\documentclass[a4,12pt]{article}

\usepackage{amsmath,amssymb,array,ascmac,bm,color,colortbl,fancyhdr,geometry,graphicx,here,listings,lscape,makeidx,mathptm,multirow,subfigure,tabularx,theorem,times,txfonts,url,verbatim}
\usepackage[square,comma,numbers,sort&compress]{natbib}

{\theorembodyfont{\normalfont}
\newtheorem{definition}{Definition}
\newtheorem{propanddef}[definition]{Proposition and definition}
\newtheorem{example}{Example}

}

\newtheorem{proposition}{Proposition}
\newtheorem{lemma}{Lemma}
\newtheorem{theorem}{Theorem}
\newtheorem{corollary}{Corollary}

\newcommand{\qed}{\hfill$\blacksquare$\par}

\makeatletter
\long\def\@makecaption#1#2{
  \vskip\abovecaptionskip
  \iftdir\sbox\@tempboxa{#1\hskip1zw#2}
    \else\sbox\@tempboxa{#1.~ #2}
  \fi
  \ifdim \wd\@tempboxa >\hsize
    \iftdir #1\hskip1zw#2\relax\par
      \else #1.~ #2\relax\par\fi
  \else
    \global \@minipagefalse
    \hbox to\hsize{\hfil\box\@tempboxa\hfil}
  \fi
  \vskip\belowcaptionskip}
\makeatother

\definecolor{orangered}{cmyk}{0,1,0.5,0}

\definecolor{navyblue}{cmyk}{0.94,0.54,0,0}

\definecolor{pinegreen}{cmyk}{0.92,0,0.59,0.25}

\newcommand{\asconv}{\stackrel{\mathrm{\footnotesize a.s.}}{\longrightarrow}}

\newcommand{\lto}{\longrightarrow}

\newcommand{\E}{\mathrm{E}}

\newcommand{\Var}{\mathrm{Var}}

\newcommand{\LA}{L_{A^{\ast}}}

\newcommand{\as}{\quad \mathrm{a.s.}}

\newcommand{\bs}[1]{\boldsymbol{#1}}

\title{\bf Optimal string clustering based on a Laplace-like mixture and EM algorithm on a set of strings}

\author{
Hitoshi Koyano,${}^{1 \ast}$ Morihiro Hayashida,${}^{2}$ and Tatsuya Akutsu${}^{2}$\\ \\
\normalsize{${}^{1}$Laboratory of Biostatistics and Bioinformatics,}\\
\normalsize{Graduate School of Medicine, Kyoto University,}\\
\normalsize{54 Kawahara-cho, Shogoin, Sakyo-ku, Kyoto 606-8507,
Japan}\\
\normalsize{${}^{2}$Laboratory of Mathematical Bioinformatics,}\\
\normalsize{Institute for Chemical Research, Kyoto University,}\\
\normalsize{Gokasho, Uji, Kyoto 611-0011, Japan}\\
\normalsize{$^\ast$Corresponding author. E-mail: koyano@kuhp.kyoto-u.ac.jp}
}

\date{}

\makeindex

\begin{document}

\maketitle

\begin{abstract}
In this study, we address the problem of clustering string data in an
unsupervised manner by developing a theory of a mixture model and an
EM algorithm for string data based on probability theory on a
topological monoid of strings developed in our previous studies.
We first construct a parametric distribution on a set of strings in
the motif of the Laplace distribution on a set of real numbers and
reveal its basic properties.
This Laplace-like distribution has two parameters: a string that
represents the location of the distribution and a positive real number
that represents the dispersion.
It is difficult to explicitly write maximum likelihood estimators of
the parameters because their log likelihood function is a complex
function, the variables of which include a string;
however, we construct estimators that almost surely converge to the
maximum likelihood estimators as the number of observed strings
increases and demonstrate that the estimators strongly consistently
estimate the parameters.
Next, we develop an iteration algorithm for estimating the parameters
of the mixture model of the Laplace-like distributions and demonstrate
that the algorithm almost surely converges to the EM algorithm for the
Laplace-like mixture and strongly consistently estimates its
parameters as the numbers of observed strings and iterations increase.
Finally, we derive a procedure for unsupervised string clustering from
the Laplace-like mixture that is asymptotically optimal in the sense
that the posterior probability of making correct classifications is
maximized.
\end{abstract}


\section{Introduction}
\label{section:Introduction}

Numbers and numerical vectors account for a large portion of data.
However, the amount of string data generated has increased
dramatically in recent years.
For example, large amounts of text data have been produced on the Web.
In the life sciences, large amounts of data regarding genes, RNAs, and
proteins have been generated.
These data are nucleotide or amino acid sequences and are represented
as strings.
Consequently, methods for analyzing string data are required in many
fields, including computer science and the life sciences.
Many studies have been conducted, especially on classification and
clustering methods for analyzing string data.
At present, a procedure of converting the strings into numerical
vectors using string
kernels~\cite{Haussler_1999,Watkins_1999a,Lodhi_2001,Leslie_2002,Paass_2002,Leslie_2003,Leslie_2004b,Vishwanathan_2004,Zien_2000,Vert_2002,Saigo_2004,Li_2005}
and subsequently applying a support vector machine that works in a
numerical vector space (see, for
example,~\cite{Aizerman_1964,Boser_1992,Cortes_1995,Drucker_1997,Vapnik_1998})
to the vectors is frequently used to classify strings.
A widely used approach to clustering strings is calculating a distance
matrix for the strings and then applying the $k$-means or nearest
neighbor method to the matrix.

However, the conversion using string kernels is not one-to-one and
thus entails a loss of information.
Furthermore, in analyzing string data, it is not necessarily natural
to convert the strings into numerical vectors because, as described in
detail in
Subsections~\ref{subsection:Summary_of_the_theory_of_randoms_strings}
and~\ref{subsection:Basic_properties} of the Appendix, a set of
strings and a numerical vector space are completely different in the
mathematical structure and the structures of a numerical vector space
that a set of strings does not have are implicitly used in analyzing
the converted data.
Thus, it is natural that the observed elements of a set of strings,
i.e., string data, be analyzed as elements of a set of strings
instead of as elements of a numerical vector space.
The most serious problem common to existing methods for analyzing
string data is that their performance has never been evaluated using
probability theory in a theoretical manner, considering that a given
string data set is part of a population generated according to a
probability law.
The performance of the methods has been generally evaluated in a
numerical manner in which they are applied to certain data sets and
cross-validation is repeated.
However, the results of performance evaluations conducted in this
manner frequently vary greatly depending on the data sets used.

Statistical methods for numerical data were rigorously constructed
based on probability theory on a set of real numbers and a numerical
vector space to make it possible to analyze data, considering that we
make an inference of a population based on a part of it observed
according to a probability law.
Similarly, statistical methods for string data should also be
constructed based on probability theory on a set of strings.
In this study, we address the problem of clustering string data in an
unsupervised manner by applying probability theory on a set of strings
developed in~\cite{Koyano_2010,Koyano_2014,Koyano_2014b}.
Our approach to unsupervised string clustering in this study is based
on a mixture model of distributions on a set of strings.
We put special emphasis on evaluating the performance of a proposed
method using probability theory in a theoretical manner, not repeating
cross-validation in a numerical manner.

Introducing a parametric probability distribution to be used as
components of a mixture model on a set of strings and constructing an
EM algorithm~\cite{Dempster_1977,Mclachlan_1997} for the mixture
model~\cite{Pearson_1894,McLachlan_2000} are critical components of
the aforementioned approach.
No parametric distributions have been introduced on a set of strings
other than an analogy of the Poisson distribution and its extension,
introduced in~\cite{Koyano_2010}.
Therefore, we first introduce an analogy of the Laplace distribution
on a set of strings that has two parameters, a string that represents
the location of the distribution and a positive real number that
represents the dispersion in
Section~\ref{section:Laplace-like_distribution_on_a_set_of_strings}.
This Laplace-like distribution on a set of strings is designed to
represent the distributions of orthologous genes.
Its basic properties are examined in detail in
Subsection~\ref{subsection:Basic_properties}.
In this preliminary study, we also reconsider the conventional
definition of a median string and assert that it is not valid under
all distance functions on a set of strings because various distance
functions are defined on a set of strings (see the second paragraph of
the next section) and each provides the set with a complex metric
structure that is completely different, for example, from a numerical
vector space.
Next, we consider the problem of estimating the parameters of the
Laplace-like distribution in
Section~\ref{section:Estimation_of_the_location_and_dispersion_parameters}.
It is difficult to explicitly write a maximum likelihood estimator of
the location parameter, and thus of the dispersion parameter, using
analytic expressions or algorithms because the maximum likelihood
estimation problem is a maximization problem of a complex function
whose arguments include a string.
We construct estimators that almost surely converge to the maximum
likelihood estimators as the number of observed strings increases and
demonstrate that the estimators strongly consistently estimate the
parameters.
Next, we develop an iteration algorithm for estimating the parameters
of the mixture model of the Laplace-like distributions on a set of
strings in
Section~\ref{section:EM_algorithm_for_the_Laplace-like_mixture_on_A^ast}.
Although the EM algorithm cannot be explicitly written in this
estimation problem, we demonstrate that the composed algorithm almost
surely converges to the EM algorithm for the Laplace-like mixture and
strongly consistently estimates its parameters as the numbers of
observed strings and iterations increase, choosing an initial value
that satisfies a given condition.
Finally, in
Section~\ref{section:Procedure_for_clustering_strings__based_on_the_Laplace-like_mixture},
we derive a procedure for clustering strings in an unsupervised manner
from the Laplace-like mixture and describe that it is asymptotically
optimal in the sense that the posterior probability of making correct
classifications is maximized.
A summary of probability theory on a set of strings that was developed
in~\cite{Koyano_2010,Koyano_2014,Koyano_2014b} and is used in this
paper is provided in
Subsection~\ref{subsection:Summary_of_the_theory_of_randoms_strings}.
In Subsection~\ref{subsection:Strong_consistency_of_MLE_in
the_case_where_the_parameter_space_is_A^ast}, we show that a result
on the strong consistency of maximum likelihood estimators in the case
where the parameter space is the set of real numbers (see, for
example,~\cite{Wald_1949,Perlman_1972}) can be extended to the case
where the parameter space is the direct product of the sets of strings
and real numbers.
The result on the strong consistency in this case is applied to
demonstrate Theorems~\ref{theorem:asymptotic_MLE}
and~\ref{theorem:strong_consistency_of_the_EM_algorithm} in
Sections~\ref{section:Estimation_of_the_location_and_dispersion_parameters}
and~\ref{section:EM_algorithm_for_the_Laplace-like_mixture_on_A^ast},
respectively.
Proofs of all results are provided in
Subsection~\ref{subsection:Proofs_of_the_results}.

\section{Laplace-like distribution on a set of strings}
\label{section:Laplace-like_distribution_on_a_set_of_strings}

Let $A = \{ a_{1}, \cdots, a_{z - 1} \}$ be an alphabet composed of $z
- 1$ letters.
We set $a_{z} = e$ for an empty letter $e$ and refer to $\bar{A} = \{
a_{1}, \cdots, a_{z} \} = A \cup \{ e \}$ as an extended alphabet.
In this study, we define a string on $A$ as a finite sequence of
elements of $A$ to the end of which an infinite sequence $o = e
\cdots$ of empty letters is appended.
Defining a string in this manner, a random string is naturally defined
in a manner in which it can realize strings of varying lengths (see
the third and fourth paragraphs of
Subsection~\ref{subsection:Summary_of_the_theory_of_randoms_strings}).
We denote the set of all strings on $A$ by $A^{\ast}$.
Our objects in this study are sequences of random variables that take
values in $A^{\ast}$ (random strings), their distributions, and their
realizations, i.e., string data.
In the following sections, we use a fundamental framework of a
probability theory on $A^{\ast}$ that was proposed
in~\cite{Koyano_2010} and extended in~\cite{Koyano_2014,Koyano_2014b}.
A summary of this theory is provided in
Subsection~\ref{subsection:Summary_of_the_theory_of_randoms_strings}.
See also the supplemental material of~\cite{Koyano_2010} for details
on the theory of random strings.

In this study, we consider distance functions that take values in
$\mathbb{N}$ (the set of natural numbers including zero) as distance
functions on $A^{\ast}$ and denote a set of them by $D$.
$D$ includes the longest common subsequence distance (see, for
example,~\cite{Bergroth_2000}), the Levenshtein
distance~\cite{Levenshtein_1965} (denoted by $d_{L}$), and the
Damerau--Levenshtein distance~\cite{Damerau_1964}.
See, for example,~\cite{Navarro_2001} for a review of the distances on
$A^{\ast}$ (the Jaro--Winkler distance~\cite{Jaro_1989,Winkler_1990}
is not a distance on $A^{\ast}$ in a mathematical sense because it
does not obey the triangle inequality).
In this paper, we regard the deletion of consecutive letters at the
end of a string as the substitution of them into empty letters and the
insertion of letters to the end of a string as the substitution of
empty letters appended to the end of the string into the letters.
We refer to the minimum number of substitutions in this sense required
to transform one string into another as the extended Hamming distance
and denote it by $d_{H^{\prime}}$.
The ordinary Hamming distance~\cite{Hamming_1950} (denoted by $d_{H}$)
is not a mathematical distance on $A^{\ast}$ because it is defined
only between strings with equal lengths, but $d_{H^{\prime}}$ is a
mathematical distance on $A^{\ast}$.
$d_{H^{\prime}}$ is regarded as the distance on $A^{\ast}$ in which
the fewest types of edit operation are allowed.
The analysis of biological sequences is started from multiple
alignment (see, for example,~\cite{Waterman_1995}) in many cases.
It is reasonable that the distance between multiple aligned sequences
is measured in terms of $d_{H^{\prime}}$ and, therefore,
$d_{H^{\prime}}$ is of special importance in biological sequence
analysis.
If, for example, $d_{L}$ is used as a distance function on $A^{\ast}$,
we write $d = d_{L}$.

We set $U_{d}(s, r) = \{ t \in A^{\ast}: d(s, t) \leq r \}$ and
$\partial U_{d}(s, r) = \{ t \in A^{\ast}: d(s, t) = r \}$ for $s \in
A^{\ast}$, $r \in \mathbb{N}$, and $d \in D$.
We denote the number of elements of a set $S$ by $|S|$ and the power
set of $S$ by $2^{S}$.
The length of $s \in A^{\ast}$ is the number of elements of $A$ that
compose $s$ and is represented by $|s|$ (see the third paragraph of
Subsection~\ref{subsection:Summary_of_the_theory_of_randoms_strings}).
We begin with introducing a parametric probability distribution on
$A^{\ast}$.

\begin{propanddef}
\label{propanddef:Laplace-like_distribution_on_A^ast}
We define the function $q_{d}( \ \cdot \ ; \lambda, \rho): A^{\ast}
\to [0, 1]$ as
\begin{equation}
q_{d}(s; \lambda, \rho) = \frac{1}{(\rho + 1) \bigl| \partial
U_{d}(\lambda, d(s, \lambda)) \bigr|} \left( \frac{\rho}{\rho + 1}
\right)^{d(s, \lambda)} \label{eq:prob_fun_Laplace-like_dist}
\end{equation}
for any $\lambda \in A^{\ast}$, $\rho \in (0, \infty)$, and $d \in D$.
Then, $q_{d}( \ \cdot \ ; \lambda, \rho)$ is a probability function on
$A^{\ast}$.
Thus, we define the set function $Q_{d}( \ \cdot \ ; \lambda, \rho):
2^{A^{\ast}} \to [0, 1]$ as
\[
Q_{d}(E; \lambda, \rho) = \sum_{s \in E} q_{d}(s; \lambda, \rho)
\]
and refer to $Q_{d}( \ \cdot \ ; \lambda, \rho)$
as a Laplace-like distribution on $A^{\ast}$ with parameters $\lambda$
and $\rho$ (denoted by $\LA(\lambda, \rho)$).
\end{propanddef}

We refer to a string that a random string $\sigma$ generates with the
highest probability as a mode string of $\sigma$ and denote it by
$\bs{M}_{m}(\sigma)$.
After reconsidering a conventional definition of a median string of
$n$ strings, we introduce a median string $\bs{M}(\sigma)$ of $\sigma$
in Subsection~\ref{subsection:Basic_properties}.
We also use a consensus sequence $\bs{M}_{c}(\sigma)$ of $\sigma = \{
\alpha_{j}: j \in \mathbb{Z}^{+} \}$, a sequence of letters to which a
marginal distribution of $\alpha_{j}$ assigns the highest probability
for each $j \in \mathbb{Z}^{+}$, as a measure of the location of a
distribution on $A^{\ast}$ in addition to $\bs{M}_{m}(\sigma)$ and
$\bs{M}(\sigma)$ and use mean absolute deviations around
$\bs{M}_{m}(\sigma), \bs{M}(\sigma)$, and $\bs{M}_{c}(\sigma)$ as
measures of the dispersion ($\mathbb{Z}^{+}$ represents the set of
positive integers).
Their definitions are provided in
Subsection~\ref{subsection:Summary_of_the_theory_of_randoms_strings}.
Although the symbol $\bs{M}_{c}(\sigma)$ is used in
Subsection~\ref{subsection:Summary_of_the_theory_of_randoms_strings}
only if a consensus sequence of $\sigma$ is unique, in the other
sections, $\bs{M}_{c}(\sigma)$ represents one of the consensus
sequences if it is not unique.

As shown in Subsection~\ref{subsection:Basic_properties} and
Section~\ref{section:Estimation_of_the_location_and_dispersion_parameters},
the Laplace-like distribution on $A^{\ast}$ has properties similar to
those of the Laplace distribution~\cite{Laplace_1774,Kotz_2001} on
$\mathbb{R}$ (the set of real numbers) in the following respects.
(i) It has two parameters, $\lambda$ and $\rho$, that represent the
location and dispersion of the distribution, respectively.
(ii) Its probability function $q_{d}(s; \lambda, \rho)$ monotonically
decreases as $d(s, \lambda)$ becomes larger (therefore, it is unimodal
at $\lambda$) and assigns an equal probability to strings with equal
distances to $\lambda$ (therefore, it is symmetric with respect to
$\lambda$).
(iii) $q_{d}(s; \lambda, \rho)$ decreases exponentially as $d(s,
\lambda)$ becomes larger and does not have inflection points, in
contrast to the normal distribution on $\mathbb{R}$.
(iv) If a random string $\sigma$ is distributed according to
$\LA(\lambda, \rho)$, we have $\bs{M}_{m}(\sigma) = \lambda$ for
any $d \in D$.
In contrast to the Laplace distribution on $\mathbb{R}$,
$\bs{M}(\sigma) = \lambda$ generally does not hold, but if $d =
d_{H^{\prime}}$, there exists $\ell \in \mathbb{N}$ such that for any
$t \in A^{\ast}$ satisfying $|t| = \ell$, we have $\bs{M}(\sigma) =
\bs{M}_{c}(\sigma) = \bs{M}_{m}(\sigma) \cdot t$, where $\cdot$
represents the concatenation (see
Subsection~\ref{subsection:Summary_of_the_theory_of_randoms_strings}
for the definition).
(v) Furthermore, the mean absolute deviation of $\sigma$ around
$\lambda$ is equal to $\rho$ for any $d \in D$.
(vi) It has the maximum entropy among all distributions on $A^{\ast}$
that satisfy the condition that the mean absolute deviation around
some fixed string is equal to a given positive real number.
(vii) If $d = d_{H^{\prime}}$, the maximum likelihood estimators of
parameters $\lambda$ and $\rho$ are asymptotically equal to a
truncated consensus sequence (introduced in the next section) and a
mean absolute deviation around it, respectively.
In the following sections, we drop $d$ from $U_{d}$, $\partial U_{d}$,
$q_{d}$, and $Q_{d}$.

\section{Estimation of the parameters of the Laplace-like distribution on $A^{\ast}$}
\label{section:Estimation_of_the_location_and_dispersion_parameters}

In this section, we consider the problem of estimating the location
and dispersion parameters $\lambda \in A^{\ast}$ and $\rho \in (0,
\infty)$ of $\LA(\lambda, \rho)$.
We suppose that observed strings $s_{1}, \cdots, s_{n} \in A^{\ast}$
are given.
We denote the relative frequency of $a_{h}$ at the $j$-th site of
$s_{1}, \cdots, s_{n}$ by $f_{j h}$ for each $j \in \mathbb{Z}^{+}$
and $h \in \{ 1, \cdots, z \}$ and set
\begin{equation}
j^{\ast} = \min \left\{ j \in \mathbb{Z}^{+}: f_{j z} = \max_{1 \leq h
\leq z} f_{j h} \right\} - 1.
\label{eq:def_j^ast}
\end{equation}
$j^{\ast}$ is the number of the last site at which it is guaranteed
that a letter that has a maximum relative frequency is nonempty for
$s_{1}, \cdots, s_{n}$.
If $j^{\ast} \geq 1$, we introduce the following condition
U$^{(\epsilon)}$ for $\epsilon > 0$.
U$^{(\epsilon)}$: There exists $j \in \{ 1, \cdots, j^{\ast} \}$ such
that for any $j^{\prime} \in \{ j, \cdots, j^{\ast} \}$,
\[
\max_{1 \leq h \leq z - 1} f_{j^{\prime} h} - \min_{1 \leq h \leq z -
1} f_{j^{\prime} h} < \epsilon
\]
holds.
If $\epsilon$ is sufficiently small, this inequality means that the
relative frequencies of all nonempty letters are almost uniform at the
$j^{\prime}$-th site.
$\lnot$U$^{(\epsilon)}$ represents the negation of U$^{(\epsilon)}$.
If $\{ j \in \{ 1, \cdots, j^{\ast} \}: \mbox{U}^{(\epsilon)} \mbox{
holds} \} \neq \emptyset$, we set
\begin{equation}
j^{(\epsilon)} = \min \left\{ j \in \left\{ 1, \cdots, j^{\ast} \right\}:
\mbox{U}^{(\epsilon)} \mbox{ holds} \right\} - 1. \label{eq:def_j^dagger}
\end{equation}
We put $\eta(j) = \arg \max_{1 \leq h \leq z} f_{j h}$.
If the maximizer of the right-hand side is not unique, we set
$\eta(j)$ as an arbitrary one of them.
We define an estimator $\bs{m}_{c}^{(\epsilon)}(s_{1}, \cdots, s_{n})$
of $\lambda$ (referred to as a truncated consensus sequence of $s_{1},
\cdots, s_{n}$) as
\begin{equation}
\bs{m}_{c}^{(\epsilon)}(s_{1}, \cdots, s_{n}) = \left\{
\begin{array}{ll}
o & \mbox{if } [j^{\ast} = 0] \lor [j^{\ast} \geq 1 \land
\mbox{U}^{(\epsilon)} \land j^{(\epsilon)} = 0], \\
a_{\eta(1)} \cdots a_{\eta(j^{\ast})} e \cdots & \mbox{if } j^{\ast} \geq 1
\land \lnot\mbox{U}^{(\epsilon)}, \label{eq:def_m_c^epsilon} \\
a_{\eta(1)} \cdots a_{\eta(j^{(\epsilon)})} e \cdots & \mbox{if } j^{\ast} \geq
1 \land \mbox{U}^{(\epsilon)} \land j^{(\epsilon)} \geq 1,
\end{array}
\right.
\end{equation}
where $\land$ and $\lor$ represent conjunction and disjunction,
respectively.
Noting Proposition~\ref{prop:consensus_seq_of_Laplace-like_dist} in
Subsection~\ref{subsection:Basic_properties}, if $d = d_{H^{\prime}}$,
$\bs{m}_{c}^{(\epsilon)}(s_{1}, \cdots, s_{n})$ becomes a reasonable
estimator of $\lambda$ as $\epsilon \lto 0$.
As described in
Proposition~\ref{prop:mean_abs_dev_of_Laplace-like_dist} in
Subsection~\ref{subsection:Basic_properties}, parameter $\rho$ is
equal to the mean absolute deviation around $\lambda$ for any $d \in
D$, and it is therefore reasonable to estimate $\rho$ by $\sum_{i =
1}^{n} d(s_{i}, \lambda) / n$ for any $d \in D$ if $\lambda$ is
known.
If $\lambda$ is unknown, $\sum_{i = 1}^{n} d(s_{i},
\bs{m}_{c}^{(\epsilon)}(s_{1}, \cdots, s_{n})) / n$ is a reasonable
estimator of $\rho$ for $d = d_{H^{\prime}}$.
In this section, we describe results in regard to the accuracy of
these estimators and their relation with the maximum likelihood
estimators of $\lambda$ and $\rho$.

Let $\{ \sigma_{i} = \{ \alpha_{i j}: j \in \mathbb{Z}^{+} \}: i \in
\mathbb{Z}^{+} \} \subset \mathcal{M}(\Omega, A^{\ast})$.
We set
\[
p(i, j, h) = P \bigl( \{ \omega \in \Omega: \alpha_{i j}(\omega) =
a_{h} \} \bigr), \quad \bar{p}(j, h, n) = \frac{1}{n} \sum_{i = 1}^{n}
p(i, j, h)
\]
for each $h = 1, \cdots, z$.
$p(i, j, h)$ represents the probability that the $j$-th letter of the
$i$-th random string realizes the $h$-th letter in the extended
alphabet $\bar{A}$, and $\bar{p}(j, h, n)$ represents the average
probability that the $h$-th letter in $\bar{A}$ is observed when $n$
observations are made.
The definitions of $\mathcal{M}(\Omega, A^{\ast}),
[\mathcal{M}(\Omega, A^{\ast})], [\mathcal{M}(\Omega, A^{\ast})^{n}]$
for $n \in \mathbb{Z}^{+}$, and $[\mathcal{M}(\Omega,
A^{\ast})^{n}]_{j}$ for $j \in \mathbb{N}$ are provided in
Subsection~\ref{subsection:Summary_of_the_theory_of_randoms_strings}.
Let S a.s. for any statement S and $\asconv$ represent that S holds
with probability one and the almost sure convergence, respectively.
First, we demonstrate the following lemma.

\begin{lemma}
\label{lemma:estimating_a_mode_string_of_a_unimodal_symmetric_distribution}
We suppose that (i) $d = d_{H^{\prime}}$ and that (ii) there exist $m
= \{ m_{1}, \cdots, m_{|m|}, e, \cdots \} \in A^{\ast}$ and $\ell \in
\mathbb{N}$ such that for any $t \in A^{\ast}$ satisfying $|t| =
\ell$, a consensus sequence of a distribution on $A^{\ast}$ is
represented by $m \cdot t$.
Let $s_{1}, \cdots, s_{n}$ be realizations of random strings
$\sigma_{1} = \{ \alpha_{1 j}: j \in \mathbb{Z}^{+} \}, \cdots,
\sigma_{n} = \{ \alpha_{n j}: j \in \mathbb{Z}^{+} \}$.
If (iii) $\alpha_{1 j}, \cdots, \alpha_{n j}$ are independent for each
$j \in \mathbb{Z}^{+}$ and (iv) there exists $N_{0} \in
\mathbb{Z}^{+}$ such that if $n \geq N_{0}$, then $(\sigma_{1},
\cdots, \sigma_{n}) \in [\mathcal{M}(\Omega, A^{\ast})^{n}]_{|m|}$
a.s. holds, $\iota(j) = \arg \max_{1 \leq h \leq z - 1} \bar{p}(j, h,
n)$ is uniquely determined independently of $n$ for each $j = 1,
\cdots, |m|$, and $\{ a_{\iota(1)}, \cdots, a_{\iota(|m|)} \} = \{
m_{1}, \cdots, m_{|m|} \}$ holds,
then there exists $\epsilon_{0} > 0$ such that if $n \geq N_{0}$ and
$\epsilon \leq \epsilon_{0}$, we have
\[
\bs{m}_{c}^{(\epsilon)}(s_{1}, \cdots, s_{n}) = m \as
\]
\end{lemma}

If $\sigma_{1}, \cdots, \sigma_{n}$ are independent, then $\alpha_{1 j},
\cdots, \alpha_{n j}$ are also independent for each $j \in
\mathbb{Z}^{+}$, but the converse is not true.
In
Lemma~\ref{lemma:estimating_a_mode_string_of_a_unimodal_symmetric_distribution},
the independence of $\alpha_{1 j}, \cdots, \alpha_{n j}$ is supposed
for each $j \in \mathbb{Z}^{+}$ (condition~(iii)), but the
independence of $\sigma_{1}, \cdots, \sigma_{n}$ is not.
Furthermore, in
Lemma~\ref{lemma:estimating_a_mode_string_of_a_unimodal_symmetric_distribution},
it is not supposed that $\sigma_{1}, \cdots, \sigma_{n}$ have an
identical distribution on $A^{\ast}$ a consensus sequence of which is
represented by $m \cdot t$.
If $\sigma_{1}, \cdots, \sigma_{n}$ have a distribution that is
unimodal at $m$ and symmetric about $m$, conditions~(ii) and~(iv) are
satisfied.
Condition~(ii) can hold even if $\sigma_{1}, \cdots, \sigma_{n}$
have different distributions if they are unimodal at $m$ and symmetric
about $m$ (see
Lemma~\ref{lemma:median_string_of_unimodal_and_symmetric_dist} in
Subsection~\ref{subsection:Basic_properties}).
More importantly, for condition~(iv) to hold, it is not necessary
that $\sigma_{1}, \cdots, \sigma_{n}$ have the identical mode
string $m$ as well as that $\sigma_{1}, \cdots, \sigma_{n}$ have
distributions unimodal at $m$ and symmetric about $m$.

From
Lemma~\ref{lemma:estimating_a_mode_string_of_a_unimodal_symmetric_distribution},
we immediately obtain the following result with respect to the
strongly consistent estimation of parameter $\lambda$ of
$\LA(\lambda, \rho)$.

\begin{proposition}
\label{prop:MLE_location_parameter}
We consider the problem of estimating parameter $\lambda$ of
$\LA(\lambda, \rho)$ with $d = d_{H^{\prime}}$ based on realizations
$s_{1}, \cdots, s_{n}$ of random strings $\sigma_{1} = \{ \alpha_{1
j}: j \in \mathbb{Z}^{+} \}, \cdots, \sigma_{n} = \{ \alpha_{n j}: j
\in \mathbb{Z}^{+} \}$.
If conditions (iii) and (iv) of
Lemma~\ref{lemma:estimating_a_mode_string_of_a_unimodal_symmetric_distribution},
in which $m$ is replaced with $\lambda$, are satisfied, there exist
$N_{0} \in \mathbb{Z}^{+}$ and $\epsilon_{0} > 0$ such that if $n \geq
N_{0}$ and $\epsilon \leq \epsilon_{0}$, we have
\[
\bs{m}_{c}^{(\epsilon)}(s_{1}, \cdots, s_{n}) = \lambda \as,
\]
i.e., $\bs{m}_{c}^{(\epsilon)}(s_{1}, \cdots, s_{n})$ strongly
consistently estimates $\lambda$ for sufficiently small $\epsilon$.
\end{proposition}

Proposition~\ref{prop:MLE_location_parameter} holds when $\rho$ is
known or unknown.

\begin{proposition}
\label{prop:MLE_dispersion_parameter}
We consider the problem of estimating parameter $\rho$ based on
realizations $s_{1}, \cdots, s_{n}$ of random strings $\sigma_{1} = \{
\alpha_{1 j}: j \in \mathbb{Z}^{+} \}, \cdots, \sigma_{n} = \{
\alpha_{n j}: j \in \mathbb{Z}^{+} \}$ (i) that are distributed
according to $\LA(\lambda, \rho)$.
(a) If $\lambda$ is known, (ii) $d \in D$ is arbitrary, and (iii)
$\sigma_{1}, \cdots, \sigma_{n}$ are independent, we have
\[
\frac{1}{n} \sum_{i = 1}^{n} d (s_{i}, \lambda) \asconv \rho
\]
as $n \lto \infty$.
Furthermore, (b) if $\lambda$ is unknown, (ii$^{\prime}$) $d =
d_{H^{\prime}}$, (iii$^{\prime}$) $\alpha_{1 j}, \cdots, \alpha_{n j}$
are independent for each $j \in \mathbb{Z}^{+}$, we have
\[
\frac{1}{n} \sum_{i = 1}^{n} d_{H^{\prime}} (s_{i},
\bs{m}_{c}^{(\epsilon)}(s_{1}, \cdots, s_{n})) \asconv \rho
\]
as $n \lto \infty$ and $\epsilon \lto 0$.
\end{proposition}

Therefore, we obtained the strong consistent estimators of $\lambda$
and $\rho$.
Are these the maximum likelihood estimators?
We first describe the relation between $\sum_{i = 1}^{n} d(s_{i},
\lambda) / n$ and the maximum likelihood estimator of $\rho$.

\begin{proposition}
\label{prop:mean_abs_dev_and_MLE_dispersion}
For any $d \in D$, the maximum likelihood estimator of parameter
$\rho$ of $\LA(\lambda, \rho)$ is given by
\[
\check{\rho}(s_{1}, \cdots, s_{n}) = \frac{1}{n} \sum_{i = 1}^{n}
d(s_{i}, \lambda)
\]
if another parameter $\lambda$ is known.
If $\lambda$ is unknown, the maximum likelihood estimator of $\rho$ is
obtained by replacing $\lambda$ on the right-hand side of the above
equation with its maximum likelihood estimator.
\end{proposition}

Next, we consider the maximum likelihood estimation of $\lambda$.
Noting that $|\partial U(\lambda, d(s_{i}, \lambda))| \geq 1$ and
$\rho / (\rho + 1) < 1$, we need to find $\lambda$ that maximizes
\begin{equation}
F(\lambda, \rho) = - \sum_{i = 1}^{n} \log \bigl| \partial U(\lambda,
d(s_{i}, \lambda)) \bigr| + \log \left( \frac{\rho}{\rho + 1} \right)
\sum_{i = 1}^{n} d(s_{i}, \lambda), \label{eq:object_function}
\end{equation}
given $\rho$.
The function $F(\lambda, \rho)$ that determines an estimate of
$\lambda$ depends on $\rho$, whereas the estimator of $\rho$ depends
on $\lambda$ from
Proposition~\ref{prop:mean_abs_dev_and_MLE_dispersion}.
Furthermore, seeking a formula for the size of a sphere of strings is
an open problem.
Therefore, it is difficult to solve the maximization problem with
respect to $\lambda$, considering both terms of the right-hand side of
Equation~(\ref{eq:object_function}) simultaneously.
However, in the case of $d = d_{L}$, several approximation algorithms
for seeking a minimizer of $\sum_{i = 1}^{n} d_{L}(s_{i}, \lambda)$
have been proposed.
Thus, in this case, it would be natural to seek an approximate
solution of the maximization problem of
Equation~(\ref{eq:object_function}) according to a procedure provided
in
Subsection~\ref{subsection:Estimation_procedure_under_the_Levenshtein_distance}.

In the case of $d = d_{H^{\prime}}$, there exists an interesting
relation between $\bs{m}_{c}^{(\epsilon)}(s_{1}, \cdots, s_{n})$ and
the maximum likelihood estimator of $\lambda$.
In this case, given a sufficiently large number of observed strings,
we can explicitly write the maximum likelihood estimators
$\check{\lambda}(s_{1}, \cdots, s_{n})$ and $\check{\rho}(s_{1},
\cdots, s_{n})$ of parameters $\lambda$ and $\rho$.

\begin{theorem}
\label{theorem:asymptotic_MLE}
We consider the problem of estimating parameters $\lambda$ and $\rho$
based on realizations $s_{1}, \cdots, s_{n}$ of random strings
$\sigma_{1}, \cdots, \sigma_{n}$ that are (i) independent and (ii)
distributed according to $\LA(\lambda, \rho)$ (iii) with $d =
d_{H^{\prime}}$.
There exist $N_{0} \in \mathbb{Z}^{+}$ and $\epsilon_{0} \geq 0$ such
that if $n \geq N_{0}$ and $\epsilon \leq \epsilon_{0}$, we have
\begin{eqnarray*}
\check{\lambda}(s_{1}, \cdots, s_{n}) &=& \bs{m}_{c}^{(\epsilon)}(s_{1},
\cdots, s_{n}) \as, \\
\check{\rho}(s_{1}, \cdots, s_{n}) &=& \frac{1}{n} \sum_{i = 1}^{n}
d_{H^{\prime}} (s_{i}, \bs{m}_{c}^{(\epsilon)}(s_{1}, \cdots, s_{n})) \as
\end{eqnarray*}
\end{theorem}

\section{Estimation algorithm for the Laplace-like mixture on $A^{\ast}$}
\label{section:EM_algorithm_for_the_Laplace-like_mixture_on_A^ast}

Let $s_{1}, \cdots, s_{n}$ be $n$ observed strings from a population
distributed according to the mixture model
\[
\bs{q}(s; \bs{\theta}) = \sum_{g = 1}^{k} \pi_{g}
q(s; \lambda_{g}, \rho_{g})
\]
of $k$ Laplace-like distributions on $A^{\ast}$ with the unknown
parameter $\bs{\theta} = (\pi_{1}, \cdots, \pi_{k}, \lambda_{1},
\cdots, \lambda_{k}, \rho_{1}, \cdots,\\ \rho_{k})$.
The parameter space of this model is $\Theta = (0, 1)^{k} \times
(A^{\ast})^{k} \times (0, \infty)^{k}$.
In this section, we develop an iteration algorithm for estimating
$\bs{\theta}$ based on $s_{1}, \cdots, s_{n}$ and then investigate its
accuracy and relation with the EM algorithm for the Laplace-like
mixture on $A^{\ast}$.

\subsection{General form of the EM algorithm for the Laplace-like mixture on $A^{\ast}$}
\label{subsection:General_form_of_the_EM_algorithm}

We denote the $i$-th observed string by $s_{i} = \{ x_{i j} \in
\bar{A}: j \in \mathbb{Z}^{+} \}$ for each $i = 1, \cdots, n$.
We suppose that $s_{i}$ is a realization of a random string
$\sigma_{i}$.
For each $g = 1, \cdots, k$, we define a $k$-dimensional real vector
$\bs{w}_{g} = (w_{g 1}, \cdots, w_{g k})$ by $w_{g g} = 1$ and $w_{g
g^{\prime}} = 0$ for $g^{\prime} \neq g$, and we set $W = \{
\bs{w}_{1}, \cdots, \bs{w}_{k} \}$.
Let $\bs{Z}_{i} = (Z_{i 1}, \cdots, Z_{i k})$ be a $k$-dimensional
latent random vector that takes values in $W$.
We define the probability function of the distribution of $\bs{Z}_{i}$
as
\begin{equation}
P(\bs{Z}_{i} = \bs{w}_{g}) = \prod_{g^{\prime} = 1}^{k}
\pi_{g^{\prime}}^{w_{g g^{\prime}}}. \label{eq:distribution_of_Z_i}
\end{equation}
Because $P(\bs{Z}_{i} = \bs{w}_{g}) = \pi_{g}$ holds, $\bs{Z}_{i} =
\bs{w}_{g}$ and $P(\bs{Z}_{i} = \bs{w}_{g})$ are interpreted to
represent the event that the $i$-th string is collected from the
$g$-th subpopulation and the probability that this event occurs,
respectively.
The probability function of the conditional distribution of
$\bs{Z}_{i}$ given $\sigma_{1}(\omega) = s_{1}, \cdots,
\sigma_{n}(\omega) = s_{n}$ is calculated as
\[
P_{\bs{\theta}}(\bs{Z}_{i} = \bs{w}_{g}| \sigma_{1}(\omega) = s_{1},
\cdots, \sigma_{n}(\omega) = s_{n}) = \frac{\pi_{g} q(s_{i}|
\lambda_{g}, \rho_{g})}{\sum_{g^{\prime} = 1}^{k} \pi_{g^{\prime}}
q(s_{i}| \lambda_{g^{\prime}}, \rho_{g^{\prime}})}.
\]
We set
\begin{equation}
\zeta_{i g} = \E_{\bs{\theta}}[Z_{i g}| \sigma_{1}(\omega) = s_{1}, \cdots,
\sigma_{n}(\omega) = s_{n}], \quad
\hat{\zeta}_{i g} = \E_{\hat{\bs{\theta}}}[Z_{i g}| \sigma_{1}(\omega) =
s_{1}, \cdots, \sigma_{n}(\omega) = s_{n}]
\label{eq:def_zeta_and_zeta_hat}
\end{equation}
for some estimator $\hat{\bs{\theta}}$ of $\bs{\theta}$.
For each $i = 1, \cdots, n$, we introduce a $k$-dimensional real vector
$\bs{z}_{i} = (z_{i 1}, \cdots, z_{i k})$ defined by
\[
z_{i g} = 1 \mbox{ and } z_{i g^{\prime}} = 0 \mbox{ for } g^{\prime}
\neq g \Longleftrightarrow s_{i} \mbox{ was collected from the } g
\mbox{-th subpopulation}.
\]
$\bs{z}_{i}$ is an unknown constant vector that is defined after the
$i$-th string was observed.
We suppose that $\sigma_{1}, \cdots, \sigma_{n}$ and $\bs{Z}_{1},
\cdots, \bs{Z}_{n}$ are independent.
We first demonstrate the following lemma that holds for any $d \in D$:

\begin{lemma}
\label{lemma:EM_algorithm_for_the_Laplace-like_mixture_general_distance}
For any $d \in D$, the EM algorithm for the Laplace-like mixture on
$A^{\ast}$ has the following form.

{\bf 1} \ Choose arbitrary initial values $\hat{\pi}_{g}^{(0)},
\hat{\lambda}_{g}^{(0)}$, and $\hat{\rho}_{g}^{(0)}$ of the parameters
for each $g = 1, \cdots, k$.

{\bf 2} \ For $t = 1, 2, \cdots$,

\quad \quad {\bf 2.1} \ Compute
\begin{equation}
\hat{\zeta}_{i g}^{(t)}
= \frac{\hat{\pi}_{g}^{(t - 1)} q(s_{i}| \hat{\lambda}_{g}^{(t - 1)},
\hat{\rho}_{g}^{(t - 1)})}{\sum_{g^{\prime} = 1}^{k}
\hat{\pi}_{g^{\prime}}^{(t - 1)} q(s_{i}|
\hat{\lambda}_{g^{\prime}}^{(t - 1)}, \hat{\rho}_{g^{\prime}}^{(t - 1)})}
\label{eq:zeta_hat_ig_t}
\end{equation}

\quad \quad \quad \quad for each $i = 1, \cdots, n$ and $g = 1, \cdots, k$.

\quad \quad {\bf 2.2} \ Compute
\begin{eqnarray}
\hat{\pi}_{g}^{(t)} &=& \frac{1}{n} \sum_{i = 1}^{n}
\hat{\zeta}_{i g}^{(t)}, \label{eq:pi_hat_g_t} \\
\hat{\lambda}_{g}^{(t)} &=& \arg \max_{\lambda_{g} \in A^{\ast}}
\sum_{i = 1}^{n} \hat{\zeta}_{i g}^{(t)} \left\{ -\log \left| \partial
U(\lambda_{g}, d(s_{i}, \lambda_{g})) \right| + d(s_{i}, \lambda_{g})
\log \left( \frac{\hat{\rho}_{g}^{(t - 1)}}{\hat{\rho}_{g}^{(t - 1)}
+ 1} \right) \right\}, \label{eq:lambda_hat_g_t} \\
\hat{\rho}_{g}^{(t)} &=& \frac{1}{\sum_{i = 1}^{n} \hat{\zeta}_{i
g}^{(t)}} \sum_{i = 1}^{n} \hat{\zeta}_{i g}^{(t)}
d(s_{i}, \hat{\lambda}_{g}^{(t)})
\label{eq:update_v_j}
\end{eqnarray}

\quad \quad \quad \quad for each $g = 1, \cdots, k$.

\quad \quad {\bf 2.3} \ If $\hat{\pi}_{g}^{(t)}$,
$\hat{\lambda}_{g}^{(t)}$, and $\hat{\rho}_{g}^{(t)}$ are sufficiently
close to $\hat{\pi}_{g}^{(t - 1)}$, $\hat{\lambda}_{g}^{(t - 1)}$, and
$\hat{\rho}_{g}^{(t - 1)}$, respectively, for

\quad \quad \quad \quad each $g = 1, \cdots, k$, terminate the
iteration and return $\hat{\pi}_{g}^{(t)}, \hat{\lambda}_{g}^{(t)}$,
and $\hat{\rho}_{g}^{(t)}$. Otherwise,

\quad \quad \quad \quad increment $t$ by one and return to Step 2.1.
\end{lemma}

In the case of $d = d_{L}$, it would be a natural method of updating
$\hat{\lambda}^{(t)}_{g}$ in Step 2.2 using the procedure described in
Subsection~\ref{subsection:Estimation_procedure_under_the_Levenshtein_distance}.

\subsection{Estimation algorithm under the extended Hamming distance}
\label{subsection:In_the_case_of_the_Hamming_distance}

In this subsection, we investigate the algorithm of
Lemma~\ref{lemma:EM_algorithm_for_the_Laplace-like_mixture_general_distance}
in the case of $d = d_{H^{\prime}}$ in detail.
We prespecify $\epsilon > 0$.
We put
\begin{equation}
f_{g j h} = \frac{1}{\sum_{i = 1}^{n} z_{i g}} \sum_{i \in \{
i^{\prime} \in \{ 1, \cdots, n \}: x_{i^{\prime} j} = a_{h} \}} z_{i
g}, \quad
\varphi_{g j h} = \frac{1}{\sum_{i = 1}^{n} \hat{\zeta}_{i g}}
\sum_{i \in \{ i^{\prime} \in \{ 1, \cdots, n \}:
x_{i^{\prime} j} = a_{h} \}} \hat{\zeta}_{i g} \label{eq:f_g_j_h}
\end{equation}
for each $g = 1, \cdots, k$, $j \in \mathbb{Z}^{+}$, and $h = 1,
\cdots, z$.
$f_{g j h}$ is the relative frequency of $a_{h}$ at the $j$-th site of
strings collected from the $g$-th subpopulation, and $\varphi_{g j h}$
is an estimator of the probability that the $j$-th letter of a string
from the $g$-th subpopulation is equal to $a_{h}$.
We set $j_{g}^{\ast}$ and $\eta(g, j)$ by replacing $f_{j h}$ in
Equation~(\ref{eq:def_j^ast}) and in the definition of $\eta(j)$,
respectively, with $\varphi_{g j h}$.
Let U$_{g}^{(\epsilon)}$ be the condition obtained by replacing
$j^{\ast}$ and $f_{j h}$ in the condition U$^{(\epsilon)}$ with
$j^{\ast}_{g}$ and $\varphi_{g j h}$, respectively.
$j^{(\epsilon)}_{g}$ is defined by replacing $j^{\ast}$ and
U$^{(\epsilon)}$ in Equation~(\ref{eq:def_j^dagger}) with
$j^{\ast}_{g}$ and U$_{g}^{(\epsilon)}$, respectively.
We introduce an estimator $\hat{\lambda}_{g}^{(\epsilon)}$ of
$\lambda_{g}$ as
\begin{equation}
\hat{\lambda}_{g}^{(\epsilon)} = \hat{\lambda}_{g}^{(\epsilon)}(s_{1},
\cdots, s_{n}) =
\left\{
\begin{array}{ll}
o & \mbox{if } [j_{g}^{\ast} = 0] \lor [j_{g}^{\ast} \geq 1 \land
\mbox{U}_{g}^{(\epsilon)} \land j_{g}^{(\epsilon)} = 0], \\
a_{\eta(g, 1)} \cdots a_{\eta(g, j_{g}^{\ast})} e \cdots & \mbox{if }
j_{g}^{\ast} \geq 1 \land \lnot\mbox{U}_{g}^{(\epsilon)},
\label{eq:lambda_hat_g} \\
a_{\eta(g, 1)} \cdots a_{\eta(g, j_{g}^{(\epsilon)})} e \cdots & \mbox{if }
j_{g}^{\ast} \geq 1 \land \mbox{U}_{g}^{(\epsilon)} \land
j_{g}^{(\epsilon)} \geq 1.
\end{array}
\right.
\end{equation}

$\hat{\lambda}_{g}^{(\epsilon)}$ is regarded as a probabilistic
extension of the truncated consensus sequence
$\bs{m}_{c}^{(\epsilon)}(s_{1}, \cdots, s_{n})$ provided by
Equation~(\ref{eq:def_m_c^epsilon}) to the case where it is unknown
from which subpopulation each observation was collected.
We denote $\hat{\lambda}_{g}^{(\epsilon)}$ obtained by replacing
$\hat{\zeta}_{i g}$ in Equation~(\ref{eq:f_g_j_h}) with
$\hat{\zeta}_{i g}^{(t)}$ by $\hat{\lambda}_{g}^{(t, \epsilon)}$.
We abbreviate the algorithm of
Lemma~\ref{lemma:EM_algorithm_for_the_Laplace-like_mixture_general_distance}
that uses $\hat{\lambda}_{g}^{(t, \epsilon)}$ as an estimate of
$\lambda_{g}$ at iteration step $t$ as Algorithm $H^{\prime}$.
We investigate the asymptotic property of Algorithm $H^{\prime}$ in
the following.
Because we develop an asymptotic theory with respect to $n$ and $t$,
we denote estimates of $\zeta_{i g}, \pi_{g}, \lambda_{g}$, and
$\rho_{g}$ from Algorithm $H^{\prime}$ by $\hat{\zeta}_{i g}^{(n, t,
\epsilon)}, \hat{\pi}^{(n, t, \epsilon)}_{g}, \hat{\lambda}^{(n, t,
\epsilon)}_{g}$, and $\hat{\rho}^{(n, t, \epsilon)}_{g}$, respectively.
Let $n_{g}$ be the number of strings collected from the $g$-th
subpopulation, and set $n^{\ast} = \min \{ n_{1}, \cdots, n_{k} \}$.
We denote strings from the $g$-th subpopulation among observed strings
$s_{1}, \cdots, s_{n}$ by $s_{g 1}, \cdots, s_{g n_{g}}$ and a random
string that generates $s_{g i}$ by $\sigma_{g i}$ for each $i = 1,
\cdots, n_{g}$ ($\sigma_{g i}$ is one of $\sigma_{1}, \cdots,
\sigma_{n}$ but which is $\sigma_{g i}$ is unknown).
We denote the true value of the parameter by $\bs{\theta}^{\ast} =
(\pi_{1}^{\ast}, \cdots, \pi_{k}^{\ast}, \lambda_{1}^{\ast}, \cdots,
\lambda_{k}^{\ast}, \rho_{1}^{\ast}, \cdots, \rho_{k}^{\ast})$ and set
$\zeta_{i g}^{\ast} = \E_{\bs{\theta}^{\ast}} [Z_{i g}|
\sigma_{1}(\omega) = s_{1}, \cdots, \sigma_{n}(\omega) = s_{n}]$
(see Equation~(\ref{eq:def_zeta_and_zeta_hat})).
In the setting described above, we first obtain the following
proposition.

\begin{proposition}
\label{prop:asymptotic_analysis_of_the_EM_algorithm}
(a) Let $\hat{\pi}_{g}^{(n, t)}, \hat{\lambda}_{g}^{(n, t)}$, and
$\hat{\rho}_{g}^{(n, t)}$ be estimators of $\pi_{g}, \lambda_{g}$, and
$\rho_{g}$, respectively, with parameters $n$ and $t$, and let
$\hat{\zeta}_{i g}^{(n, t)}$ be an estimator of $\zeta_{i g}$
calculated according to Equation~(\ref{eq:zeta_hat_ig_t}) based on
$\hat{\pi}_{g}^{(n, t - 1)}, \hat{\lambda}_{g}^{(n, t - 1)}$, and
$\hat{\rho}_{g}^{(n, t - 1)}$ for each $i = 1, \cdots, n$ and $g = 1,
\cdots, k$.
For any $d \in D$, if $\hat{\pi}_{g}^{(n, t)} \asconv
\pi_{g}^{\ast}, \hat{\lambda}_{g}^{(n, t)} = \lambda_{g}^{\ast}$ a.s.,
and $\hat{\rho}_{g}^{(n, t)} \asconv \rho_{g}^{\ast}$ hold as $n_{g},
t \lto \infty$ for each $g = 1, \cdots, k$, we have
\[
\hat{\zeta}_{i g}^{(n, t)} \asconv \zeta_{i g}^{\ast}
\]
as $n_{g}, t \lto \infty$.
Conversely, (b) in the case of $d = d_{H^{\prime}}$, if
$\hat{\zeta}_{i g}^{(n, t, \epsilon)} \asconv \zeta_{i g}^{\ast}$
holds as $n_{g}, t \lto \infty$ and $\epsilon \lto 0$, we have
\[
\hat{\pi}_{g}^{(n, t, \epsilon)} \asconv \pi_{g}^{\ast}, \quad
\hat{\lambda}_{g}^{(n, t, \epsilon)} = \lambda_{g}^{\ast} \as, \quad
\hat{\rho}_{g}^{(n, t, \epsilon)} \asconv \rho_{g}^{\ast}
\]
as $n_{g}, t \lto \infty$ and $\epsilon \lto 0$.
\end{proposition}

We note that
Proposition~\ref{prop:asymptotic_analysis_of_the_EM_algorithm} means
that in the case of $d = d_{H^{\prime}}$, $(\hat{\zeta}_{1 1}^{(n, t,
\epsilon)}, \cdots, \hat{\zeta}_{n k}^{(n, t, \epsilon)}) \asconv
(\zeta_{1 1}^{\ast}, \cdots, \zeta_{n k}^{\ast})\\ \Longleftrightarrow$
[the estimate $\hat{\bs{\theta}}^{(n, t, \epsilon)} =
(\hat{\pi}_{1}^{(n, t, \epsilon)}, \cdots, \hat{\pi}_{k}^{(n, t,
\epsilon)}, \hat{\lambda}_{1}^{(n, t, \epsilon)}, \cdots,
\hat{\lambda}_{k}^{(n, t, \epsilon)}, \hat{\rho}_{1}^{(n, t,
\epsilon)}, \cdots, \hat{\rho}_{k}^{(n, t, \epsilon)})$ from
Algorithm $H^{\prime}$ converges to the true value
$\bs{\theta}^{\ast}$ of the parameter of the Laplace-like mixture],
not that $(\hat{\zeta}_{1 1}^{(n, t, \epsilon)}, \cdots,\\
\hat{\zeta}_{n k}^{(n, t, \epsilon)}) \asconv (\zeta_{1 1}^{\ast},
\cdots, \zeta_{n k}^{\ast}) \Longleftrightarrow$ [the estimates
$\hat{\lambda}_{1}^{(n, t, \epsilon)}, \cdots,
\hat{\lambda}_{k}^{(n, t, \epsilon)}, \hat{\rho}_{1}^{(n, t,
\epsilon)}, \cdots, \hat{\rho}_{k}^{(n, t, \epsilon)}$ from
Algorithm $H^{\prime}$ converge to the maximum likelihood estimates
$\check{\lambda}_{1}, \cdots, \check{\lambda}_{k}, \check{\rho}_{1},
\cdots, \check{\rho}_{k}$ of the parameters of the subpopulation
distributions], as $n^{\ast}, t \lto \infty$ and $\epsilon \lto 0$.
From Proposition~\ref{prop:asymptotic_analysis_of_the_EM_algorithm},
Algorithm $H^{\prime}$ strongly consistently estimates the parameter
$\bs{\theta}$ of the Laplace-like mixture on $A^{\ast}$ as $n, t \lto
\infty$ and $\epsilon \lto 0$ if the approach of $(\hat{\zeta}_{1
1}^{(n, t, \epsilon)}, \cdots, \hat{\zeta}_{n k}^{(n, t, \epsilon)})
$ to $(\zeta_{1 1}^{\ast}, \cdots, \zeta_{n k}^{\ast})$ in Step~2.1
and the approach of $\hat{\bs{\theta}}^{(n, t, \epsilon)}$ to
$\bs{\theta}^{\ast}$ in Step~2.2 are alternately repeated through
iteration steps.

Next, using
Proposition~\ref{prop:asymptotic_analysis_of_the_EM_algorithm}, we
consider whether Algorithm $H^{\prime}$ can estimate the true value
$\bs{\theta}^{\ast}$ of the parameter of the Laplace-like mixture on
$A^{\ast}$ with $d = d_{H^{\prime}}$.
We approach this problem by clarifying under what conditions and with
what initial value Algorithm $H^{\prime}$ estimates the true value
$\bs{\theta}^{\ast}$.
We introduce the following two conditions:

C$_{1}$: For any $(z_{1 1}^{\prime}, \cdots, z_{n k}^{\prime}) \in [0,
1]^{n k}$ that satisfies $\sum_{g = 1}^{k} z_{i g}^{\prime} = 1$ for
each $i = 1, \cdots, n$ and $(z_{1 1}^{\prime}, \cdots, z_{n
k}^{\prime})\\ \neq (z_{1 1}, \cdots, z_{n k})$,
\[
\max_{(\lambda_{1}, \cdots, \lambda_{k}, \rho_{1}, \cdots, \rho_{k})}
\frac{1}{n} \sum_{g = 1}^{k} \sum_{i = 1}^{n} z_{i g}^{\prime} \log
q(s_{i}; \lambda_{g}, \rho_{g})
< \max_{(\lambda_{1}, \cdots, \lambda_{k}, \rho_{1}, \cdots, \rho_{k})}
\frac{1}{n} \sum_{g = 1}^{k} \sum_{i = 1}^{n} z_{i g} \log q(s_{i};
\lambda_{g}, \rho_{g}) \as
\]
holds as $n^{\ast} \lto \infty$.

Condition C$_{1}$ is satisfied, for example, if a sufficiently
large number of observed strings are given and there does not exist a
$g$-th subpopulation distribution that has a log likelihood greater
than or equal to $\ell_{g}^{\ast}$ based on a portion of strings
collected from the $g$-th subpopulation and/or including strings
collected from other subpopulations, where $\ell_{g}^{\ast}$
represents the log likelihood of a $g$-th subpopulation distribution
that has the maximum log likelihood based on all strings collected
from the $g$-th subpopulation for each $g = 1, \cdots, k$.
By Proposition~\ref{prop:strong_consistency_of_MLE_on_A^ast_times_R},
maximizing the likelihood leads to estimating the true population
distribution for the Laplace-like distribution on $A^{\ast}$, and
C$_{1}$ is therefore a natural condition to exclude pathological
situations.

C$_{2}$: The solution $(\zeta_{1 1}^{\dagger}, \cdots, \zeta_{n
k}^{\dagger}, \lambda_{1}^{\dagger}, \cdots, \lambda_{k}^{\dagger},
\rho_{1}^{\dagger}, \cdots, \rho_{k}^{\dagger})$ of the maximization
problem of
\begin{equation}
\frac{1}{n} \sum_{g = 1}^{k} \sum_{i = 1}^{n} \zeta_{i g} \log
q(s_{i}; \lambda_{g}, \rho_{g}) \label{eq:fun_of_zeta_lambda_rho}
\end{equation}
with respect to $(\zeta_{1 1}, \cdots, \zeta_{n k}, \lambda_{1},
\cdots, \lambda_{k}, \rho_{1}, \cdots, \rho_{k}) \in (0, 1)^{n k}
\times (A^{\ast})^{k} \times (0, \infty)^{k}$ is unique for given
$s_{1}, \cdots, s_{n} \in A^{\ast}$.

We set
\begin{equation}
\hat{\bs{\theta}}^{(n, 0, \tau, \epsilon)} = \arg
\max_{\hat{\bs{\theta}}^{(n, 0)} \in \Theta} \frac{1}{n} \sum_{g =
1}^{k} \sum_{i = 1}^{n} \hat{\zeta}_{i g}^{(n, \tau, \epsilon)} \log
q(s_{i}; \hat{\lambda}_{g}^{(n, \tau, \epsilon)}, \hat{\rho}_{g}^{(n,
\tau, \epsilon)})
\label{eq:def_theta_hat_dagger_n_zero}
\end{equation}
for each $n, \tau \in \mathbb{Z}^{+}$, where $\hat{\zeta}_{i g}^{(n,
\tau, \epsilon)}, \hat{\lambda}_{g}^{(n, \tau, \epsilon)}$, and
$\hat{\rho}_{g}^{(n, \tau, \epsilon)}$ represent estimates of
$\zeta_{i g}, \lambda_{g}$, and $\rho_{g}$, respectively, that
Algorithm $H^{\prime}$ with the initial value $\hat{\bs{\theta}}^{(n,
0)}$ returns at iteration step $\tau$.
$\hat{\bs{\theta}}^{(n, 0, \tau, \epsilon)}$ is an initial value with
which Algorithm $H^{\prime}$ returns estimates of $\zeta_{1 1},
\cdots, \zeta_{n k}, \lambda_{1}, \cdots, \lambda_{k}, \rho_{1},
\cdots, \rho_{k}$ that maximize
Equation~(\ref{eq:def_theta_hat_dagger_n_zero}) at iteration step
$\tau$ of all possible initial values.  The following theorem answers
the question described above.

\begin{theorem}
\label{theorem:strong_consistency_of_the_EM_algorithm}
If conditions C$_{1}$ and C$_{2}$ hold and (i) there exists an
initial value $\tilde{\bs{\theta}}^{(n, 0)} = (\tilde{\pi}_{1}^{(n,
0)}, \cdots,\\ \tilde{\pi}_{k}^{(n, 0)}, \tilde{\lambda}_{1}^{(n, 0)},
\cdots, \tilde{\lambda}_{k }^{(n, 0)}, \tilde{\rho}_{1}^{(n, 0)},
\cdots, \tilde{\rho}_{k}^{(n, 0)}) \in \Theta$ such that the estimate
$\tilde{\bs{\theta}}^{(n, t, \epsilon)}$ from Algorithm $H^{\prime}$
strongly consistently estimates $\bs{\theta}$ as $n^{\ast}, t \lto
\infty$ and $\epsilon \lto 0$, then $\bs{\theta}$ is strongly
consistently estimated by the estimator $\hat{\bs{\theta}}^{(n, t,
\tau, \epsilon)}$ from Algorithm $H^{\prime}$ with the initial
value $\hat{\bs{\theta}}^{(n, 0, \tau, \epsilon)}$ given by
Equation~(\ref{eq:def_theta_hat_dagger_n_zero}) as $n^{\ast}, t,
\tau \lto \infty$ and $\epsilon \lto 0$.
\end{theorem}

Generally, in estimating a parameter of a model by using an iteration
algorithm in which an initial value is arbitrarily chosen, the most
reliable estimate is chosen after several initial values are provided
and the behavior of the sequence of estimates from each of them is
examined.
From Theorem~\ref{theorem:strong_consistency_of_the_EM_algorithm}, in
practical data analysis using the Laplace-like mixture on $A^{\ast}$,
if a sufficiently large number of observed strings are given, choosing
several initial values and adopting the estimate that maximizes
Equation~(\ref{eq:def_theta_hat_dagger_n_zero}) for a
sufficiently large $t$ would be a realistic approach, especially when
sequences of estimates from the initial values appear to converge to
different points of the parameter space.

Applying Theorem~\ref{theorem:strong_consistency_of_the_EM_algorithm},
we can demonstrate that Algorithm $H^{\prime}$ converges to the EM
algorithm for the Laplace-like mixture on $A^{\ast}$ that we seek but
cannot write in an explicit manner.
The definition of the convergence of a sequence of algorithms to an
algorithm is introduced in
Subsection~\ref{subsection:convergence_of_a_sequence_of_algorithms}.

\begin{theorem}
\label{theorem:asymptotic_EM}
If the conditions of
Theorem~\ref{theorem:strong_consistency_of_the_EM_algorithm} are
satisfied and the convergence testing constants in Step 2.3 are
sufficiently small, Algorithm $H^{\prime}$ with the initial value
$\hat{\bs{\theta}}^{(n, 0, \tau, \epsilon)}$ that satisfies
Equation~(\ref{eq:def_theta_hat_dagger_n_zero}) almost surely
converges to the EM algorithm with the initial value
$\hat{\bs{\theta}}^{(n, 0, \tau, \epsilon)}$ for the
Laplace-like mixture on $A^{\ast}$ as $n^{\ast}, \tau \lto
\infty$ and $\epsilon \lto 0$.
\end{theorem}

\section{String clustering procedure based on the Laplace-like mixture on $A^{\ast}$}
\label{section:Procedure_for_clustering_strings__based_on_the_Laplace-like_mixture}

In this section, we derive a procedure for clustering strings based on
the results obtained in the previous sections.
We consider the problem of clustering $n$ strings $s_{1}, \cdots,
s_{n} \in A^{\ast}$ into $k$ classes.
We assume a mixture model
\[
\sum_{g = 1}^{k} \pi_{g} q(s; \lambda_{g}, \rho_{g}) = \sum_{g =
1}^{k} \frac{\pi_{g}}{(\rho_{g} + 1) \bigl| \partial U(\lambda_{g},
d(s, \lambda_{g})) \bigr|} \left( \frac{\rho_{g}}{\rho_{g} + 1}
\right)^{d(s, \lambda_{g})}
\]
of $k$ Laplace-like distributions $\LA(\lambda_{1}, \rho_{1}), \cdots,
\LA(\lambda_{k}, \rho_{k})$ on $A^{\ast}$ with mixture coefficients
$\pi_{1}, \cdots, \pi_{k}$ as a model generating $s_{1}, \cdots,
s_{n}$.
From Bayes' theorem, the posterior probability given $s_{1}, \cdots,
s_{n}$ that $s_{i}$ belongs to the $g$-th class is provided by
\[
\pi_{\bs{\theta}}(\bs{Z}_{i} = \bs{z}_{g}| s_{1}, \cdots, s_{n}) =
\frac{\pi_{g} q(s_{i}; \lambda_{g}, \rho_{g})}{\sum_{g^{\prime} = 1}^{k}
\pi_{g^{\prime}} q(s_{i}; \lambda_{g^{\prime}}, \rho_{g^{\prime}})}
\]
for each $i = 1, \cdots, n$ and $g = 1, \cdots, k$.

\begin{corollary}
\label{corollary:optimality_of_clustering_procedure}
We suppose that the conditions of
Theorem~\ref{theorem:strong_consistency_of_the_EM_algorithm} are
satisfied and denote an estimate from Algorithm $H^{\prime}$ with the
initial value $\hat{\bs{\theta}}^{(n, 0, \tau, \epsilon)} \in
\Theta$ satisfying Equation~(\ref{eq:def_theta_hat_dagger_n_zero}) by
$\hat{\bs{\theta}}^{(n, t, \tau, \epsilon)}$.
Then, the clustering procedure that
\[
\mbox{if } g^{\ast} = \arg \max_{1 \leq g \leq k}
\pi_{\hat{\bs{\theta}}^{(n, t, \tau, \epsilon)}}(\bs{Z}_{i} = \bs{z}_{g}|
s_{1}, \cdots, s_{n}),\mbox{ then classify } s_{i} \mbox{ into the }
g^{\ast}\mbox{-th class}
\]
for each $i = 1, \cdots, n$ is asymptotically optimal in the sense
that the posterior probability of making correct classifications is
maximized as $n^{\ast}, t, \tau \lto \infty$ and $\epsilon \lto
0$.
\end{corollary}

\section{Concluding remarks}
\label{section:Concluding_remarks}

In this study, using the probability theory on $A^{\ast}$ developed
in~\cite{Koyano_2010,Koyano_2014,Koyano_2014b}, we constructed the
theory of the mixture model and the EM algorithm on $A^{\ast}$ and
derived the optimal procedure for unsupervised string clustering based
on the theory.
We encountered the interesting phenomenon that an EM algorithm for the
Laplace-like mixture on $A^{\ast}$, which we sought, could not be
written in an explicit manner because of the complex metric structure
of $A^{\ast}$, but a sequence of algorithms (i.e., a sequence of
sequences of computations) that strongly consistently estimates the
parameters of the Laplace-like mixture and converges to the EM
algorithm with probability one was obtained explicitly.
This is different from the phenomena that an algorithm halts and that
a sequence of approximate solutions from an iterative algorithm
converges.
The authors addressed the problem of supervised string classification
by constructing a theory of a statistical learning machine that works
in $A^{\ast}$ in~\cite{Koyano_2014b}.
Recently, the amount of string data has increased exponentially.
The noncommutative topological monoid $A^{\ast}$ of strings has
interesting structures that are completely different, for example,
from the Euclidean space $\mathbb{R}^{p}$ and Hilbert space $L^{2}$
and will become one of important spaces on which probability theory
and methods of statistics and machine learning should be developed in
the future.

\section*{Appendix}
\label{section:Appendix}

\renewcommand{\thesubsection}{A\arabic{subsection}}

\subsection{Summary of the theory of random strings}
\label{subsection:Summary_of_the_theory_of_randoms_strings}

In this subsection of the Appendix, we describe the definitions of
several concepts in probability theory on a set of strings used in the
main text.
See the supplemental material of~\cite{Koyano_2010} for details.
In the following, we refer to a set of a finite number of letters
\[
A = \{ a_{1}, \cdots, a_{z - 1} \}
\]
as the alphabet.
For example, $A = \{ \mathtt{a}, \mathtt{c}, \mathtt{g}, \mathtt{t} \}$
is the alphabet for gene sequences.
We denote an empty letter by $e$ and set $\bar{A} = A \cup \{ e \}$.
We denote a set of $(x_{1}, \cdots, x_{n}) \in \bar{A}^{n}$ of which a
letter with the maximum frequency is uniquely determined by
$[\bar{A}^{n}]$.
A mapping $m_{c}: [\bar{A}^{n}] \to \bar{A}$ is defined as
\[
m_{c}(x_{1}, \cdots, x_{n}) = \mbox{a letter with the maximum
frequency of } x_{1}, \cdots, x_{n}
\]
and is called a consensus letter on $[\bar{A}^{n}]$.

Let $(\Omega, \mathfrak{F}, P)$ be a probability space.
We call an $\bar{A}$-valued random variable on $\Omega$ a random letter
and denote the set of all random letters by $\mathcal{M}(\Omega,
\bar{A})$.
For the mapping $\varepsilon: \Omega \to \bar{A}$, which is defined as
$\varepsilon(\omega) = e$ for all $\omega \in \Omega$, we have
$\varepsilon \in \mathcal{M}(\Omega, \bar{A})$.
The independence of $\{ \alpha_{i}: i \in \mathbb{Z}^{+} \} \subset
\mathcal{M}(\Omega, \bar{A})$ is defined in the same manner as that of
ordinary random variables.
We denote a set of $\alpha \in \mathcal{M}(\Omega, \bar{A})$ for which
there exists $x \in \bar{A}$ such that for any $y \in \bar{A}
\smallsetminus \{ x \}$, $q(x) > q(y)$ holds by $[\mathcal{M}(\Omega,
\bar{A})]$, where $q$ is a probability function of a distribution of
$\alpha$.
A mapping $M_{c}: [\mathcal{M}(\Omega, \bar{A})] \to \bar{A}$ is
defined as
\[
M_{c}(\alpha) = x \in \bar{A} \mbox{ such that } [q(x) > q(y), \forall
y \in \bar{A} \smallsetminus \{ x \}]
\]
and is called a consensus letter on $[\mathcal{M}(\Omega, \bar{A})]$.
We denote a set of $(\alpha_{1}, \cdots, \alpha_{n}) \in
\mathcal{M}(\Omega, \bar{A})^{n}$ for which a consensus letter of
$\alpha_{1}(\omega), \cdots, \alpha_{n}(\omega)$ is uniquely determined
for any $\omega \in \Omega$ by $[\mathcal{M}(\Omega, \bar{A})^{n}]$.
A mapping $\mu_{c} :[\mathcal{M}(\Omega, \bar{A})^{n}] \to
\mathcal{M}(\Omega, \bar{A})$ is defined as
\[
\mu_{c}(\alpha_{1}, \cdots, \alpha_{n})(\omega) =
m_{c}(\alpha_{1}(\omega), \cdots, \alpha_{n}(\omega))
\]
and called a consensus letter on $[\mathcal{M}(\Omega, \bar{A})^{n}]$.

In common usage in computer science, a string on the alphabet $A = \{
a_{1}, \cdots, a_{z - 1} \}$ is a finite sequence of elements of $A$.
However, in this study, we define a string as follows, although both
definitions are essentially identical:
A sequence $s = \{ x_{j} \in \bar{A}: j \in \mathbb{Z}^{+} \}$ of
elements of $\bar{A}$ is a string on $A$ if it satisfies the following
conditions:
\[
(\mbox{i}) \mbox{ there exists } h \in \mathbb{Z}^{+} \mbox{ such that }
x_{h} = e, \mbox{ and } (\mbox{ii}) \ x_{j} = e \mbox{ implies } x_{j + 1}
= e.
\]
In other words, we define a string on $A$ as a finite sequence of
elements of $A$ to which the infinite sequence $\{ e, \cdots \}$ of
empty letters is appended.
In the following, by naturally extending the above definition of a
string, we define a random string in a manner in which it can
realize strings of varying lengths.
We denote the set of all strings on $A$ by $A^{\ast}$.
Let $d$ be a distance function on $A^{\ast}$.
If $x_{k} \neq e$ and $x_{k + 1} = e$ hold for $s = \{ x_{j}: j \in
\mathbb{Z}^{+} \} \in A^{\ast}$, we set
\[
s \cdot t = \{ x_{1}, \cdots, x_{k}, y_{1}, y_{2}, \cdots \}
\]
for any $t = \{ y_{j}: j \in \mathbb{Z}^{+} \} \in A^{\ast}$
and call $s \cdot t$ the concatenation of $s$ and $t$.
A function $| \cdot |: A^{\ast} \to \mathbb{N}$ is defined as
\[
|s| = \min \{ h \in \mathbb{Z}^{+}: x_{h} = e \} - 1, \ s = \{ x_{j}:
j \in \mathbb{Z}^{+} \}
\]
and called the length on $A^{\ast}$.
Let $(s_{1}, \cdots, s_{n}) \in (A^{\ast})^{n}$ and $s_{i} = \{ x_{i
j}:j \in \mathbb{Z}^{+} \}$ for each $i = 1, \cdots, n$.
We denote a set of $(s_{1}, \cdots, s_{n})$ for which a consensus
letter of $x_{1 j}, \cdots, x_{n j}$ is uniquely determined for any $j
\in \mathbb{Z}^{+}$ by $[(A^{\ast})^{n}]$.
A mapping $\bs{m}_{c}: [(A^{\ast})^{n}] \to A^{\ast}$ is defined
as
\[
\bs{m}_{c}(s_{1}, \cdots, s_{n}) = \{ m_{c}(x_{1 j}, \cdots, x_{n
j}): j \in \mathbb{Z}^{+} \}
\]
and is called a consensus sequence on $[(A^{\ast})^{n}]$.
A function $\bs{v}_{c}: [(A^{\ast})^{n}] \to [0, \infty)$ is defined
as
\[
\bs{v}_{c}(s_{1}, \cdots, s_{n}) = \frac{1}{n} \sum_{i = 1}^{n}
d(s_{i}, \bs{m}_{c}(s_{1}, \cdots, s_{n}))
\]
and called a mean absolute deviation around $\bs{m}_{c}(s_{1},
\cdots, s_{n})$ on $[(A^{\ast})^{n}]$ (this quantity was simply called
a variance on $[(A^{\ast})^{n}]$ in~\cite{Koyano_2010}).

Next, we introduce a random string.
A sequence of random letters $\sigma = \{ \alpha_{j} \in
\mathcal{M}(\Omega, \bar{A}): j \in \mathbb{Z}^{+} \}$ is a random
string if it satisfies the following conditions:
\begin{eqnarray*}
& (\mbox{i}) & \mbox{for any } \omega \in \Omega \mbox{ there exists } h
\in \mathbb{Z}^{+} \mbox{ such that } \alpha_{h}(\omega) = e,
\mbox{ and} \\
& (\mbox{ii}) & \alpha_{j}(\omega) = e \mbox{ for } \omega \in \Omega
\mbox{ implies } \alpha_{j + 1}(\omega) = e.
\end{eqnarray*}
We denote the set of all random strings by $\mathcal{M}(\Omega,
A^{\ast})$.
\noindent A function $| \cdot |: \mathcal{M}(\Omega, A^{\ast}) \to
\mathbb{N}$ is defined as
\[
| \sigma | = \min \{ h \in \mathbb{Z}^{+}: \alpha_{h} = \varepsilon \}
- 1, \ \sigma = \{ \alpha_{j}: j \in \mathbb{Z}^{+} \}
\]
and is called the length on $\mathcal{M}(\Omega, A^{\ast})$.
The random string defined above can be regarded as a special case of a
discrete stochastic process.
Therefore, a distribution of a random string can be defined as follows:
Let $\sigma = \{ \alpha_{j}: j \in \mathbb{Z}^{+} \} \in
\mathcal{M}(\Omega, A^{\ast})$.
A set function $Q_{\sigma; j_{1}, \cdots, j_{k}}: 2^{\bar{A}^{k}} \to
[0, 1]$ is defined as
\[
Q_{\sigma; j_{1}, \cdots, j_{k}}(E) = P \bigl( \bigl\{\omega \in
\Omega:(\alpha_{j_{1}}(\omega), \cdots, \alpha_{j_{k}}(\omega)) \in E
\bigr\} \bigr)
\]
for any $k \in \mathbb{Z}^{+}$ and $j_{1}, \cdots, j_{k} \in
\mathbb{Z}^{+}$ that satisfy $j_{1} < \cdots < j_{k}$.
$Q_{\sigma; j_{1}, \cdots, j_{k}}$ is a probability measure on
$2^{\bar{A}^{k}}$ and is called a finite-dimensional distribution of
$\sigma$ at sites $j_{1}, \cdots, j_{k}$.
A function $q_{\sigma; j_{1}, \cdots, j_{k}}: \bar{A}^{k} \to [0, 1]$ is
defined as
\[
q_{\sigma; j_{1}, \cdots, j_{k}}(x_{1}, \cdots, x_{k}) = Q_{\sigma;
j_{1}, \cdots, j_{k}} \bigl( \{ (x_{1}, \cdots, x_{k}) \} \bigr)
\]
and is called a probability function of $Q_{\sigma; j_{1}, \cdots,
j_{k}}$.
For the probability function $q_{\sigma; 1, \cdots, |\sigma|}$ of the
finite-dimensional distribution at sites $1, \cdots, |\sigma|$ of
$\sigma \in \mathcal{M}(\Omega, A^{\ast})$, we define the function
$q_{\sigma}: A^{\ast} \to [0, 1]$ as
\[
q_{\sigma}(s) = \left\{
\begin{array}{ll}
q_{\sigma; 1, \cdots, | \sigma |}(x_{1}, \cdots,
x_{| \sigma |}) & (\mbox{for } x_{1}, \cdots, x_{|\sigma|} \in \bar{A}
\mbox{ such that } \\
 & \ s = (x_{1}, \cdots, x_{|\sigma|}, e, \cdots)
\mbox{ if } |\sigma| \geq |s|) \\
0 & (\mbox{if } |\sigma| < |s|).
\end{array}
\right.
\]
$q_{\sigma}$ is a probability function on $A^{\ast}$.
The independence of the random strings is defined in the following
manner:
(1) For the finite case, $\sigma_{1} = \{ \alpha_{1 j}: j \in
\mathbb{Z}^{+} \}, \cdots, \sigma_{n} = \{ \alpha_{n j}: j \in
\mathbb{Z}^{+} \} \in \mathcal{M}(\Omega, A^{\ast})$ are independent if
$\{ \alpha_{1 j}: j \in J_{1} \}, \cdots, \{ \alpha_{n j}: j \in J_{n}
\}$ are independent for any nonempty finite set $J_{1}, \cdots, J_{n}
\subset \mathbb{Z}^{+}$.
(2) For the countably infinite case, $\{ \sigma_{i}: i \in
\mathbb{Z}^{+} \} \subset \mathcal{M}(\Omega, A^{\ast})$ are independent
if $\sigma_{i_{1}}, \cdots, \sigma_{i_{k}}$ are independent for any $k
\in \mathbb{Z}^{+}$ and $i_{1}, \cdots, i_{k} \in \mathbb{Z}^{+}$.

We denote a set of $\sigma = \{ \alpha_{j}: j \in \mathbb{Z}^{+} \}
\in \mathcal{M}(\Omega, A^{\ast})$ for which a consensus letter of
$\alpha_{j}$ is uniquely determined for any $j \in \mathbb{Z}^{+}$ by
$[\mathcal{M}(\Omega, A^{\ast})]$.
A mapping $\bs{M}_{c}: [\mathcal{M}(\Omega, A^{\ast})] \to
A^{\ast}$ is defined as
\[
\bs{M}_{c}(\sigma) = \{ M_{c}(\alpha_{j}): j \in
\mathbb{Z}^{+} \}, \ \sigma = \{ \alpha_{j}: j \in \mathbb{Z}^{+} \}
\]
and is called a consensus sequence on $[\mathcal{M}(\Omega,
A^{\ast})]$.
A function $\bs{\Upsilon}_{c}: [\mathcal{M}(\Omega, A^{\ast})] \to
[0, \infty)$ is defined as
\[
\bs{\Upsilon}_{c}(\sigma) = \sum_{s \in A^{\ast}} d(s,
\bs{M}_{c}(\sigma)) q_{\sigma}(s)
\]
and called a mean absolute deviation around
$\bs{M}_{c}(\sigma)$ on $[\mathcal{M}(\Omega, A^{\ast})]$
(this quantity was simply called a variance on $[\mathcal{M}(\Omega,
A^{\ast})]$ in~\cite{Koyano_2010}).
We denote $\bs{\Upsilon}_{c}(\sigma)$ obtained by replacing
$\bs{M}_{c}(\sigma)$ on the right-hand side of the above equation with
a median string $\bs{M}(\sigma)$ introduced in
Definition~\ref{def:tentative} in
Subsection~\ref{subsection:Basic_properties} and a mode string
$\bs{M}_{m}(\sigma)$ introduced in
Section~\ref{section:Laplace-like_distribution_on_a_set_of_strings} by
$\bs{\Upsilon}(\sigma)$ and $\bs{\Upsilon}_{m}(\sigma)$, respectively.
We refer to $\bs{\Upsilon}(\sigma)$ and $\bs{\Upsilon}_{m}(\sigma)$ as
mean absolute deviations around $\bs{M}(\sigma)$ and around
$\bs{M}_{m}(\sigma)$, respectively.

Let $(\sigma_{1}, \cdots, \sigma_{n}) \in \mathcal{M}(\Omega,
A^{\ast})^{n}$ and $\sigma_{i} = \{ \alpha_{i j}: j \in \mathbb{Z}^{+}
\}$ for each $i = 1, \cdots, n$.
$[\mathcal{M}(\Omega, A^{\ast})^{n}]$ and $[\mathcal{M}(\Omega,
A^{\ast})^{n}]_{\ell}$ represent sets of $(\sigma_{1}, \cdots,
\sigma_{n})$ for which a consensus letter of $\alpha_{1 j}(\omega),
\cdots, \alpha_{n j}(\omega)$ is uniquely determined for any $j \in
\mathbb{Z}^{+}$ and $\omega \in \Omega$ and for any $j \in \{ 1,
\cdots, \ell \}$ and $\omega \in \Omega$, respectively.
A mapping $\bs{\mu}_{c}: [\mathcal{M}(\Omega, A^{\ast})^{n}] \to
\mathcal{M}(\Omega, A^{\ast})$ is defined as
\[
\bs{\mu}_{c}(\sigma_{1}, \cdots, \sigma_{n})(\omega) = \{
\mu_{c}(\alpha_{1 j}, \cdots, \alpha_{n j})(\omega): j \in
\mathbb{Z}^{+} \}
\]
and is called a consensus sequence on $[\mathcal{M}(\Omega,
A^{\ast})^{n}]$.
A mapping $\bs{\upsilon}: [\mathcal{M}(\Omega, A^{\ast})^{n}]
\to \mathcal{M}(\Omega, [0, \infty))$ is defined as
\[
\bs{\upsilon}_{c}(\sigma_{1}, \cdots, \sigma_{n})(\omega) =
\frac{1}{n} \sum_{i = 1}^{n} d(\sigma_{i}(\omega),
\bs{\mu}_{c}(\sigma_{1}, \cdots, \sigma_{n})(\omega))
\]
and called a mean absolute deviation around $\bs{\mu}_{c}(\sigma_{1},
\cdots, \sigma_{n})(\omega))$ on $[\mathcal{M}(\Omega, A^{\ast})^{n}]$
(this quantity was simply called a variance on $[\mathcal{M}(\Omega,
A^{\ast})^{n}]$ in~\cite{Koyano_2010}).

\subsection{Basic properties of the Laplace-like distribution on $A^{\ast}$}
\label{subsection:Basic_properties}

In this subsection, we describe the basic properties of $\LA(\lambda,
\rho)$ introduced in
Section~\ref{section:Laplace-like_distribution_on_a_set_of_strings}.
We consider a distribution on $A^{\ast}$ that has a parameter $m \in
A^{\ast}$.
If there exists $d \in D$ such that the probability function $q(s; m)$
of the distribution satisfies the conditions that $d(s, m) <
d(s^{\prime}, m)$ for $s, s^{\prime} \in A^{\ast}$ implies $q(s; m) >
q(s^{\prime}; m)$ and that $d(s, m) = d(s^{\prime}, m)$ implies $q(s;
m) = q(s^{\prime}; m)$, we say that
$m$ is a location parameter of the distribution with respect to $d$.
In other words, if $q(s; m)$ is unimodal at $m$ and symmetric with
respect to $m$, then $m$ is a location parameter.
We next consider a distribution on $A^{\ast}$ that has a parameter $v
\in (0, \infty)$.
We say that $v$ is a dispersion parameter if the distribution
approaches the uniform distribution on $A^{\ast}$ as $v \lto \infty$.

\begin{proposition}
\label{prop:shape_dist}
The parameters $\lambda$ and $\rho$ of $\LA(\lambda, \rho)$ are
location and dispersion parameters, respectively.
\end{proposition}

{\bf Proof.}
$q(s; \lambda, \rho)$ monotonically decreases from $1 / (\rho + 1)
> 0$ to zero as $d(s, \lambda)$ increases from zero because $(\rho /
(\rho + 1))^{d(s, \lambda)}$ and $1/|\partial U(\lambda, d(s,
\lambda))|$ monotonically decrease from one to zero as $d(s, \lambda)$
increases from zero.
In addition, $q(s; \lambda, \rho)$ depends on $s$ only through $d(s,
\lambda)$, and therefore, $q(s; \lambda, \rho) = q(s^{\prime}; \lambda,
\rho)$ holds if $d(s, \lambda) = d(s^{\prime}, \lambda)$ for $s,
s^{\prime} \in A^{\ast}$.
Thus, $\lambda$ is a location parameter.
$q(s; \lambda, \rho)$ approaches the uniform distribution on
$A^{\ast}$ as $\rho$ increases if and only if $q(\lambda; \lambda,
\rho)$ decreases as $\rho$ increases, which clearly holds from
$q(\lambda; \lambda, \rho) = 1/(\rho + 1)$.
Hence, $\rho$ is a dispersion parameter.
\qed

Here, we reconsider the conventional definition of a median string of
strings to describe a result on the relation of parameter $\lambda$ of
$\LA(\lambda, \rho)$ with a median string and a consensus sequence of
a random string $\sigma$ that is distributed according to
$\LA(\lambda, \rho)$.
A median string of $s_{1}, \cdots, s_{n} \in A^{\ast}$ cannot be
defined as a middle string after arranging $s_{1}, \cdots, s_{n}$ in
ascending order because $A^{\ast}$ is not a totally ordered set.
Therefore, \cite{Kohonen_1985} introduced the median string and set
median string of $S \subset A^{\ast}$ as
\begin{equation}
\mathrm{med}(S) = \arg \min_{s \in A^{\ast}} \sum_{t \in S} d(s, t), \quad
\mathrm{med}^{\prime}(S) = \arg \min_{s \in S} \sum_{t \in S} d(s, t),
\label{eq:Median_string_of_strings}
\end{equation}
respectively,
because a median $\mathrm{med}(x_{1}, \cdots, x_{n})$ of $x_{1},
\cdots, x_{n} \in \mathbb{R}$ is characterized as
\[
\mathrm{med}(x_{1}, \cdots, x_{n}) = \arg \min_{y \in \mathbb{R}}
\sum_{i = 1}^{n} |x_{i} - y|, \mbox{ or equivalently, }
\mathrm{med}(x_{1}, \cdots, x_{n}) = \arg \min_{y \in \mathbb{R}}
\frac{1}{n} \sum_{i = 1}^{n} |x_{i} - y|.
\]
In addition to the characterization as a minimizer of the first-order
absolute moment, a median has another characterization
by~\cite{Tukey_1975} as the deepest point in the sample or the
distribution (Tukey's depth median).
These two characterizations of a median are equivalent in $\mathbb{R}$
but not in $\mathbb{R}^{p}$ (the $p$-dimensional real vector space)
for $p \geq 2$, and they define different multidimensional medians
(see, for example,~\cite{Brown_1983,Oja_1983,Oja_1985,Donoho_1992}).
However, the characterization of a median as the deepest point in the
sample or the distribution cannot be used to define a median string
because the deepest point is defined using the concepts of the total
order and projection.
A median string is not necessarily unique, similar to an ordinary
median on $\mathbb{R}$.
A median string and a consensus sequence play an important role as a
measure of the center of strings in computer science.
See, for
example,~\cite{delaHigure_2000,Martinez_2003,Jiang_2003,Nicolas_2003,Nicolas_2005,Olivares_2008,Jiang_2012}
for theoretical results and applications of a median string.

We consider the problem of introducing a median string of a
probability distribution on $A^{\ast}$ or a random string.
One natural, but tentative, definition, which extends
Equation~(\ref{eq:Median_string_of_strings}) to a probabilistic
version, is as follows:

\begin{definition}
\label{def:tentative}
We suppose that a distance $d$ on $A^{\ast}$ is given and that a random
string $\sigma$ has a distribution on $A^{\ast}$ with a probability
function $q(s)$.
We define a median string of $\sigma$ with respect to the distance $d$
as
\[
\bs{M}(\sigma) = \arg \min_{s \in A^{\ast}} \sum_{t \in A^{\ast}}
d(s, t) q(t).
\]
\end{definition}

The definition of a median string of strings by
Equation~(\ref{eq:Median_string_of_strings}) is apparently reasonable,
but in Definition~\ref{def:tentative}, which is a natural extension of
this definition to a random string, a problem in defining a median
string as a minimizer of the first-order absolute moment arises.
We consider $A^{\ast}$ with the Levenshtein distance $d_{L}$.
In $\mathbb{R}^{p}$, spheres with centers of different points and
equal radii have equal sizes, whereas in $A^{\ast}$, spheres with
centers of different strings and equal radii do not have necessarily
equal sizes.
There exist more strings near a longer string in $A^{\ast}$.
In other words, spaces $\mathbb{R}^{p}$ and $A^{\ast}$ have different
metric structures, but the definition of a median string of strings by
Equation~(\ref{eq:Median_string_of_strings}) and
Definition~\ref{def:tentative} do not consider the differences between
the metric structures of these spaces.
We consider a unimodal and symmetric distribution on $A^{\ast}$ with
respect to $m \in A^{\ast}$ (for example, $\LA(m, \rho)$).
Choosing $m^{\prime} \in A^{\ast}$ such that $|m^{\prime}| > |m|$, we
observe that
\[
\sum_{s \in A^{\ast}} d_{L}(s, m) q(s) < \sum_{s \in A^{\ast}} d_{L}(s,
m^{\prime}) q(s)
\]
is not guaranteed because $|U(m, n)| < |U(m^{\prime}, n)|$ holds for
any $n \in \mathbb{Z}^{+}$, as illustrated in the following example.
Therefore, $m$ is not necessarily a median string according to the
above definition.

\begin{example}
\label{example:mean_absolute_deviation}
We set $A = \{ 0, 1 \}$ and consider the distribution on $A^{\ast}$
with the probability function $q(o) = 0.2, q(0) = q(1) = 0.15, q(00) =
q(01) = q(10) = q(11) = 0.125$, and $q(s) = 0$ for other $s \in
A^{\ast}$, where an infinite sequence $e \cdots$ of empty letters
connected to the end of each string was dropped.
This distribution is unimodal at $o$ (thus, $o$ is a unique mode
string) and symmetric with respect to $o$.
The mean absolute deviations around $o$, $0$, and $00$ with
respect to $d_{L}$ are equal to
\[
\sum_{s \in A^{\ast}} d_{L}(s, o) q(s) = 1.3, \
\sum_{s \in A^{\ast}} d_{L}(s, 0) q(s) = 0.975, \mbox{ and }
\sum_{s \in A^{\ast}} d_{L}(s, 00) q(s) = 1.35,
\]
respectively, and therefore, it is not minimized around $o$. \qed
\end{example}

A median string of a random string is expected to work as a measure of
the location of its distribution.
For example, for the normal and Laplace distributions on $\mathbb{R}$,
we have ``the mode $=$ the median ($=$ the expected value).''
It is desirable that a median string be defined such that the similar
relation of ``the mode string = the median string'' holds for unimodal
and symmetric distributions on $A^{\ast}$.
Considering the above-mentioned difference in the metric structure
between $A^{\ast}$ and $\mathbb{R}^{p}$,
Definition~\ref{def:tentative} can be modified as follows.

\begin{definition}
\label{def:Median_string_of_a_random_string}
Let $q(s)$ be the probability function of the distribution of a random
string $\sigma$ and $\varphi: \mathbb{N} \times \mathbb{N} \to
\mathbb{N}$ be a monotonically increasing function of both variables.
We define the median string $\bs{M}(\sigma; \varphi)$ of
$\sigma$ with respect to $\varphi$ as
\begin{equation}
\bs{M}(\sigma; \varphi) = \arg \min_{m \in A^{\ast}} \sum_{s \in
A^{\ast}} \varphi(d(s, m), |m|) q(s).
\label{eq:Median_string_of_a_random_string}
\end{equation}
\end{definition}

Under the definition of a median string of a random string by
Equation~(\ref{eq:Median_string_of_a_random_string}), if two strings
minimize the sum $\sum_{t \in A^{\ast}} d(s, t) q(t)$, the shorter one
is chosen as a median string.
Natural examples of $\varphi(d(s, m), |m|)$ include $|U(m, d(s, m))|$
and $|\partial U(m, d(s, m))|$.

\begin{example}
\label{example:mean_size_of_a_sphere}
We consider the same probability function on $A^{\ast}$ as in
Example~\ref{example:mean_absolute_deviation} under $\varphi(d(s, m),
|m|) = |\partial U(m, d_{L}(s, m))|$.
The expected size of the sphere of strings with respect to $d_{L}$ is
minimized to $\sum_{s \in A^{\ast}} |\partial U(o, d_{L}(s, o))| q(s)
= 2.8$ when $o$ is the center of the sphere.
For example, it is equal to $\sum_{s \in A^{\ast}} |\partial U(0,
d_{L}(s, 0))| q(s) = 4.775$ when $0$ is the center and, noting
$|\partial U(00, 0)| = 1$, $|\partial U(00, 1)| = 7$, and $|\partial
U(00, 2)| = 17$, it is calculated as $\sum_{s \in A^{\ast}} |\partial
U(00, d_{L}(s, 00))| q(s) = 11.0$ when $00$ is the center. \qed
\end{example}

Definition~\ref{def:Median_string_of_a_random_string} was introduced
by modifying Definition~\ref{def:tentative} to consider the difference
in the metric structure between $\mathbb{R}^{p}$ and $A^{\ast}$.
In contrast to $\mathbb{R}^{p}$, there exist several intrinsic
distance functions on $A^{\ast}$, as described in
Section~\ref{section:Laplace-like_distribution_on_a_set_of_strings}.
Consequently, it is not guaranteed that the median string provided by
Definition~\ref{def:Median_string_of_a_random_string} is equal to the
mode string for all distance functions and unimodal and symmetric
distributions on $A^{\ast}$.
However, the point here is that
Examples~\ref{example:mean_absolute_deviation}
and~\ref{example:mean_size_of_a_sphere} indicate that if we translate
a definition in $\mathbb{R}^{p}$ into $A^{\ast}$ without considering
the difference in the structure between $\mathbb{R}^{p}$ and
$A^{\ast}$, the translated definition in $A^{\ast}$ may not have some
of properties of the original definition in $\mathbb{R}^{p}$.
Furthermore, noting that even a median string as a minimizer of the
sum of distances is approximately computed and that the question of
finding a general formula for assessing the volume of a sphere of
strings is open, the above discussion also means that the complex
situation became more complex.
However, the following lemma indicates that if $d = d_{H^{\prime}}$
and a random string has a unimodal and symmetric distribution, there
exists an explicit relation among its mode string, median string
provided in Definition~\ref{def:tentative}, and consensus sequence,
and consequently, the situation is somewhat tractable.
Although in
Subsection~\ref{subsection:Summary_of_the_theory_of_randoms_strings},
the symbol $\bs{M}_{c}(\sigma)$ is used when a consensus sequence of a
random string $\sigma = \{ \alpha_{j}: j \in \mathbb{Z}^{+} \}$ is
unique, in
Lemma~\ref{lemma:median_string_of_unimodal_and_symmetric_dist} and
Proposition~\ref{prop:consensus_seq_of_Laplace-like_dist} below,
$\bs{M}_{c}(\sigma)$ represents a sequence obtained by ordering one of
letters generated by $\alpha_{j}$ with the highest probability with
respect to $j \in \mathbb{Z}^{+}$.

\begin{lemma}
\label{lemma:median_string_of_unimodal_and_symmetric_dist}
If (i) $d = d_{H^{\prime}}$ and a random string $\sigma$ has a
distribution on $A^{\ast}$ whose probability function $q(s)$ satisfies
the conditions that (ii) $d_{H^{\prime}}(s, m) <
d_{H^{\prime}}(s^{\prime}, m)$ implies $q(s) > q(s^{\prime})$ for $m
\in A^{\ast}$ and (iii) $d_{H^{\prime}}(s, m) =
d_{H^{\prime}}(s^{\prime}, m)$ implies $q(s) = q(s^{\prime})$, there
exists $\ell \in \mathbb{N}$ such that for any $t \in A^{\ast}$
satisfying $|t| = \ell$, $\bs{M}(\sigma) = \bs{M}_{c}(\sigma) = m
\cdot t$ holds.
\end{lemma}

{\bf Proof.}
{\bf (Step 1)}
Let $\sigma = \{ \alpha_{j}: j \in \mathbb{Z}^{+} \}$ and $m = \{
m_{j}: j \in \mathbb{Z}^{+} \}$.
We arbitrarily choose $j \in \{ 1, \cdots, |m| \}$.
Let $A(j, m_{j})$ be a set of strings obtained by deleting $m_{j}$
from a string in $A^{\ast}$, the $j$-th letter of which is equal to
$m_{j}$, i.e., if $\{ x_{1}, \cdots, x_{j - 1}, m_{j}, x_{j + 1},
\cdots \} \in A^{\ast}$, then $\{ x_{1}, \cdots, x_{j - 1}, x_{j + 1},
\cdots \} \in A(j, m_{j})$.
All strings in $A^{\ast}$ with the $j$-th letter equal to $m_{j}$ can
be created by inserting $m_{j}$ between the $(j - 1)$-th and $j$-th
letters of strings in $A(j, m_{j})$.
The marginal probability of $\alpha_{j}(\omega) = m_{j}$ is given
by
\begin{equation}
q_{j}(m_{j}) = \sum_{\{x_{1}, \cdots, x_{j - 1}, x_{j + 1}, \cdots
\} \in A(j, m_{j})} q(\{ x_{1}, \cdots, x_{j - 1},
m_{j}, x_{j + 1}, \cdots \}).
\label{eq:marginal_m_i}
\end{equation}
We arbitrarily choose $y_{j} \in A \smallsetminus \{ m_{j} \}$.
Defining $A(j, y_{j})$ in the same manner as $A(j, m_{j})$, we have
$A(j, m_{j}) = A(j, y_{j})$.
Therefore, the marginal probability of $\alpha_{j}(\omega) = y_{j}$ is
equal to
\begin{equation}
q_{j}(y_{j}) = \sum_{\{x_{1}, \cdots, x_{j - 1}, x_{j + 1}, \cdots \}
\in A(j, m_{j})} q(\{ x_{1}, \cdots, x_{j - 1}, y_{j}, x_{j +
1}, \cdots \}). \label{eq:marginal_y_i}
\end{equation}
We have
\[
d_{H^{\prime}}( \{ x_{1}, \cdots, x_{j - 1}, m_{j}, x_{j + 1},
\cdots \}, m) < d_{H^{\prime}}( \{ x_{1}, \cdots, x_{j - 1}, y_{j},
x_{j + 1}, \cdots \}, m)
\]
for any $\{x_{1}, \cdots, x_{j - 1}, x_{j + 1}, \cdots \} \in A(j,
m_{j})$.
Thus,
\begin{equation}
q(\{ x_{1}, \cdots, x_{j - 1}, m_{j}, x_{j + 1}, \cdots \}) >
q(\{ x_{1}, \cdots, x_{j - 1}, y_{j}, x_{j + 1}, \cdots \})
\label{eq:q_m_j_greater_q_y_j}
\end{equation}
holds from condition~(ii).
Combining Equations~(\ref{eq:marginal_m_i}),~(\ref{eq:marginal_y_i}),
and~(\ref{eq:q_m_j_greater_q_y_j}) provides
\begin{equation}
q_{j}(m_{j}) > q_{j}(y_{j}). \label{eq:q_j_m_j_<_q_j_y_j}
\end{equation}

Let $B(j)$ be a set of strings in $A^{\ast}$ with the $j$-th letter
equal to $e$.
We set $B(j, 0) = \{ o \}$ and denote a set of strings in $B(j)$ with
a nonempty $k$-th letter and a $(k + 1)$-th letter of $e$ by $B(j, k)$
for each $k = 1, \cdots, j - 1$.
We have $B(j) = \bigcup_{k = 0}^{j - 1} B(j, k)$ and $B(j, k) \cap
B(j, k^{\prime}) = \emptyset$ if $k \neq k^{\prime}$ for $k,
k^{\prime} \in \{ 0, \cdots, j - 1 \}$.
Thus,
\begin{equation}
q_{j}(e) = \sum_{s \in B(j)} q(s) = \sum_{k = 0}^{j - 1} \sum_{s \in
B(j, k)} q(s) \label{eq:marginal_prob_of_e_at_the_j-th_site}
\end{equation}
holds.
For each $k = 0, \cdots, j - 2$, $B^{\prime}(j, k)$ represents a set
of strings obtained by substituting (1) $e$s at the $(k + 1)$-th to
$(j - 1)$-th sites of each $s \in B(j, k)$ with arbitrary nonempty
letters, (2) $e$ at the $j$-th site with $m_{j}$, and (3) $e$s at the
$(j + 1)$-th to $(2j - 1 - k)$-th sites with arbitrary nonempty
letters.
We denote a set of strings obtained by substituting $e$ at the $j$-th
site of each $s \in B(j, j - 1)$ with $m_{j}$ by $B^{\prime}(j, j -
1)$ and then set $B^{\prime}(j) = \bigcup_{k = 0}^{j - 1}
B^{\prime}(j, k)$.
$B^{\prime}(j)$ is a subset of strings with the $j$-letter equal to
$m_{j}$, and we have $B^{\prime}(j, k) \cap B^{\prime}(j, k^{\prime}) =
\emptyset$ for $k \neq k^{\prime}$.
Therefore,
\begin{equation}
q_{j}(m_{j}) > \sum_{s \in B^{\prime}(j)} q(s) = \sum_{k = 0}^{j - 1}
\sum_{s \in B^{\prime}(j, k)} q(s)
\label{eq:marginal_prob_of_m_j_at_the_j-th_site}
\end{equation}
holds.
Noting that $|B(j, k)| = |B^{\prime}(j, k)|$ and $d_{H^{\prime}}(s, m)
> d_{H^{\prime}}(s^{\prime}, m)$ for any $s \in B(j, k)$ and
$s^{\prime} \in B^{\prime}(j, k)$, we obtain
\begin{equation}
\sum_{k = 0}^{j - 1} \sum_{s \in B(j, k)} q(s) < \sum_{k = 0}^{j - 1}
\sum_{s \in B^{\prime}(j, k)} q(s)
\label{eq:ineq_over_B_j_k_and_B_prime_j_k}
\end{equation}
from condition (ii).
Combining
Equations~(\ref{eq:marginal_prob_of_e_at_the_j-th_site}),~(\ref{eq:marginal_prob_of_m_j_at_the_j-th_site}),
and~(\ref{eq:ineq_over_B_j_k_and_B_prime_j_k}) leads to
\begin{equation}
q_{j}(e) < q_{j}(m_{j}). \label{eq:q_j_e_<_q_j_m_j}
\end{equation}
From Equations~(\ref{eq:q_j_m_j_<_q_j_y_j})
and~(\ref{eq:q_j_e_<_q_j_m_j}), the consensus letter of the marginal
distribution of $\alpha_{j}$ is equal to $m_{j}$ for each $j = 1,
\cdots, |m|$.

{\bf (Step 2)}
We arbitrarily choose $j \geq |m| + 1$.
For any $x \in A$, let $C_{1}(j, x)$ be a set of strings that have
$x$ at the $j$-th site,
nonempty letters at all sites before the $j$-th,
and $e$s at all sites after the $j$-th and
$C_{2}(j, x)$ be a set of strings that have
$x$ at the $j$-th site,
nonempty letters at all sites before the $j$-th,
nonempty letters at a finite number of consecutive sites more
than or equal to one from the $(j + 1)$-th,
and $e$s at all of the following sites.
Because a set of strings that have $x$ at the $j$-th site is
represented by $C_{1}(j, x) \cup C_{2}(j, x)$, we have
\begin{equation}
q_{j}(x) = \sum_{s \in C_{1}(j, x) \cup C_{2}(j, x)} q(s).
\label{eq:rep_q_j_x}
\end{equation}
We set $C(j, x, r) = \{ s \in C_{1}(j, x) \cup C_{2}(j, x):
d_{H^{\prime}}(s, m) = r \}$ for each $r \in \mathbb{N}$.
For any $x \in A$, we have $C_{1}(j, x) \cup C_{2}(j, x) = \bigcup_{0
\leq r < \infty} C(j, x, r)$ and $C(j, x, r) \cap C(j, x,
r^{\prime}) = \emptyset$ if $r \neq r^{\prime}$.
Thus,
\[
q_{j}(x) = \sum_{0 \leq r < \infty} \sum_{s \in C(j, x, r)} q(s).
\]
holds from Equation~(\ref{eq:rep_q_j_x}).
Therefore, noting $|C(j, x, r)| = |C(j, x^{\prime}, r)|$ for any $x,
x^{\prime} \in A$ and condition (iii), we have
\begin{equation}
q_{j}(x) = q_{j}(x^{\prime}). \label{eq:marginal_prob_nonempty_letter}
\end{equation}
Equation~(\ref{eq:marginal_prob_nonempty_letter}) means that all
nonempty letters are consensus letters (i.e., have a maximum marginal
probability) at the $j$-th site if $q_{j}(x) > q_{j}(e)$.

For a nonempty letter to have a maximum marginal probability at each
site, it is necessary that there exist strings that have a nonempty
letter at each site.
However, $A^{\ast}$ is a set of sequences obtained by concatenating a
finite sequence of nonempty letters with an infinite sequence of $e$s,
and a probability function of any distribution on $A^{\ast}$ assigns
probabilities only to strings in $A^{\ast}$.
Hence, there exists $J_{0} \in \mathbb{Z}^{+}$ such that if $j \geq
J_{0}$, then
\begin{equation}
q_{j}(x) < q_{j}(e) \label{eq:existence_J_0}
\end{equation}
holds for any $x \in A$.
From the result of Step 1, we have $J_{0} \geq |m| + 1$.
We set $J_{0}^{\ast} = \min \{ J_{0} \in \mathbb{Z}^{+}: \mbox{if } j
\geq J_{0}, \mbox{ Equation~(\ref{eq:existence_J_0}) holds for any } x
\in A \}$.
From Equation~(\ref{eq:marginal_prob_nonempty_letter}), to conclude
$\bs{M}_{c}(\sigma) = m \cdot t$, it suffices to demonstrate that
\begin{equation}
q_{j}(x) > q_{j}(e) \label{eq:q_j_x_>_q_j_e}
\end{equation}
holds for any $j$ satisfying $|m| + 1 \leq j \leq J_{0}^{\ast} - 1$ if
$J_{0}^{\ast} \geq |m| + 2$.

{\bf (Step 3)}
Let $D_{1}(j, x)$ be a set of strings that have $x \in A$ at the
$j$-th site, nonempty letters at the first to the $(j - 1)$-th sites,
and $e$s at all sites after the $(j + 1)$-th.
We denote a set of strings that have $x$ at the $j$-th site, nonempty
letters at the first to the $(j - 1)$-th sites and at the $(j +
1)$-th, and $e$s at all sites after the $(j + 1)$-th by $D_{2}(j, x)$.
$D_{3}(j, x)$ represents a set of strings that have $x$ at the $j$-th
site, nonempty letters at the first to the $(j - 1)$-th sites and at a
finite number of consecutive sites more than or equal to two from the
$(j + 1)$-th, and $e$s at all of the following sites.
Because $D_{i}(j, x) \cap D_{i^{\prime}}(j, x) = \emptyset$ holds if
$i \neq i^{\prime}$ for $i, i^{\prime} \in \{ 1, 2, 3 \}$ and a set of
strings that have $x$ at the $j$-th site is represented by $D_{1}(j,
x) \cup D_{2}(j, x) \cup D_{3}(j, x)$, we have
\begin{equation}
q_{j}(x) = \sum_{s \in D_{1}(j, x)} q(s) + \sum_{s \in D_{2}(j, x)}
q(s) + \sum_{s \in D_{3}(j, x)} q(s). \label{eq:q_j_x}
\end{equation}
Let $D_{4}(j, x)$ be a set of strings that have $x$ at the $(j +
1)$-th site, nonempty letters at the first to the $j$-th sites, and
$e$s at all sites after the $(j + 1)$-th and $D_{5}(j, x)$ be a set of
strings that have $x$ at the $(j + 1)$-th site, nonempty letters at
the first to the $j$-th sites and at a finite number of consecutive
sites more than or equal to one from the $(j + 2)$-th, and $e$s at all
of the following sites.
$D_{4}(j, x) \cap D_{5}(j, x) = \emptyset$ holds and a set of strings
that have $x$ at the $(j + 1)$-th site is written as $D_{4}(j, x) \cup
D_{5}(j, x)$.
Consequently,
\begin{equation}
q_{j + 1}(x) = \sum_{s \in D_{4}(j, x)} q(s) + \sum_{s \in D_{5}(j, x)}
q(s) \label{eq:q_j_plus_1_x}
\end{equation}
holds.
We set $D_{i}(j, x, r) = \{ s \in D_{i}(j, x): d_{H^{\prime}}(s, m) =
r \}$ for each $i = 2, \cdots, 5$ and $r \in \mathbb{N}$.
We have $B_{i}(j, x) = \bigcup_{0 \leq r < \infty} B_{i}(j, x, r)$ and
$B_{i}(j, x, r) \cap B_{i}(j, x, r^{\prime}) = \emptyset$ if $r \neq
r^{\prime}$.
Furthermore, $|D_{2}(j, x, r)| = |D_{4}(j, x, r)|$ and $|D_{3}(j, x,
r)| = |D_{5}(j, x, r)|$ hold.
Thus, noting condition~(iii), we obtain
\[
\sum_{s \in D_{2}(j, x)} q(s)
= \sum_{0 \leq r < \infty}\sum_{s \in B_{2}(j, x, r)} q(s)
= \sum_{0 \leq r < \infty}\sum_{s \in B_{4}(j, x, r)} q(s)
= \sum_{s \in D_{4}(j, x)} q(s)
\]
and $\sum_{s \in D_{3}(j, x)} q(s) = \sum_{s \in D_{5}(j, x)} q(s)$.
Therefore, from Equations~(\ref{eq:q_j_x})
and~(\ref{eq:q_j_plus_1_x}), we have
\begin{equation}
q_{j}(x) > q_{j + 1}(x), \label{eq:monotonically_decrease}
\end{equation}
i.e., $q_{j}(x)$ is monotonically decreasing with respect to $j \geq
|m| + 1$.

Let $D(j, e)$ be a set of strings that have nonempty letters at
consecutive sites more than or equal to zero and fewer than or equal
to $j - 1$ from the first site and $e$s at all of the following sites
and $D_{6}(j + 1, e)$ be a set of strings that have nonempty letters
at the first to the $j$-th sites and $e$s at all of the following
sites.
$D(j + 1, e) \cap D_{6}(j + 1, e) = \emptyset$ holds, and sets of
strings that have $e$ at the $j$-th site and at the $(j + 1)$-th site
are written as $D(j, e)$ and $D(j, e) \cup D_{6}(j + 1, e)$,
respectively.
Therefore, we have
\begin{equation}
q_{j}(e)
= \sum_{s \in D(j, e)} q(s) < \sum_{s \in D(j, e)} q(s)
+ \sum_{s \in D_{6}(j + 1, e)} q(s) = q_{j + 1}(e),
\label{eq:monotonically_increase}
\end{equation}
i.e., $q_{j}(e)$ is monotonically increasing with respect to $j \geq
|m| + 1$.
From Equations~(\ref{eq:monotonically_decrease})
and~(\ref{eq:monotonically_increase}), we obtain
Equation~(\ref{eq:q_j_x_>_q_j_e}).
Hence, $\bs{M}_{c}(\sigma) = m \cdot t$ has been demonstrated.

{\bf (Step 4)}
From the results of Steps~1 to~3, putting $\bs{M}_{c}(\sigma) = \{
m_{1}, \cdots, m_{|m|}, t_{1}, \cdots, t_{\ell}, e, \cdots \} = \{
m_{j}^{\prime}: j \in \mathbb{Z}^{+} \}$, we have
\begin{eqnarray*}
\sum_{s \in A^{\ast}} d_{H}(s_{j}, m_{j}^{\prime}) q(s) &<&
\sum_{s \in A^{\ast}} d_{H}(s_{j}, y_{j}) q(s),
 \ j = 1, \cdots, |m|, |m| + \ell + 1, \cdots, \\
\sum_{s \in A^{\ast}} d_{H}(s_{j}, m_{j}^{\prime}) q(s)
&\leq& \sum_{s \in A^{\ast}} d_{H}(s_{j}, y_{j}) q(s),
 \ j = |m| + 1, \cdots, |m| + \ell
\end{eqnarray*}
for any $y_{j} \in \bar{A}$ satisfying $y_{j} \neq m_{j}^{\prime}$,
where $s = \{ s_{j}: j \in \mathbb{Z}^{+} \}$.
Consequently, we obtain
\begin{equation}
\sum_{s \in A^{\ast}} \sum_{j \in \mathbb{Z}^{+}} d_{H}(s_{j},
m_{j}^{\prime}) q(s) < \sum_{s \in A^{\ast}} \sum_{j \in
\mathbb{Z}^{+}} d_{H}(s_{j}, y_{j}) q(s).
\label{eq:inequality_on_the_weighted_sum_of_Hamming_distances}
\end{equation}
If $d = d_{H^{\prime}}$, then
\begin{equation}
\bs{M}(\sigma) = \arg \min_{s^{\prime} \in A^{\ast}} \sum_{s \in A^{\ast}}
d_{H^{\prime}}(s, s^{\prime}) q(s)
= \arg \min_{s^{\prime} \in A^{\ast}} \sum_{s \in A^{\ast}} \sum_{j
\in \mathbb{Z}^{+}} d_{H}(s_{j}, s_{j}^{\prime}) q(s)
\label{eq:def_median_string}
\end{equation}
holds, where $s^{\prime} = \{ s_{j}^{\prime}: j \in \mathbb{Z}^{+} \}$.
Combining
Equations~(\ref{eq:inequality_on_the_weighted_sum_of_Hamming_distances})
and~(\ref{eq:def_median_string}) provides $\bs{M}(\sigma) =
\bs{M}_{c}(\sigma)$.
The proof is completed.
\qed

Lemma~\ref{lemma:median_string_of_unimodal_and_symmetric_dist} states
that the mode string is not equal to the median string (and the
consensus sequence) for unimodal and symmetric distributions on
$A^{\ast}$ in contrast to $\mathbb{R}$, which implies that
distributions on $A^{\ast}$ are intractable compared with those on
$\mathbb{R}$.
Noting Step 3 in the above proof, we observe that the larger the
dispersion of the distribution ($\rho$ in the case of $\LA(\lambda,
\rho)$) is, the larger $\ell$ is.
The relation between the mode string and the consensus sequence for
unimodal and symmetric distributions on $A^{\ast}$ described in
Lemma~\ref{lemma:median_string_of_unimodal_and_symmetric_dist},
especially the fact that the probabilities of all nonempty letters are
uniform at the $(|m| + 1)$-th to $(|m| + \ell)$-th sites, underlies
the theoretical results in the case of $d = d_{H^{\prime}}$ described
in the main text.
We immediately obtain the following proposition on parameter $\lambda$
of $\LA(\lambda, \rho)$ from
Lemma~\ref{lemma:median_string_of_unimodal_and_symmetric_dist}.

\begin{proposition}
\label{prop:consensus_seq_of_Laplace-like_dist}
We suppose that a random string $\sigma$ is distributed according to
$\LA(\lambda, \rho)$.
(a) For any $d \in D$, we have $\bs{M}_{m}(\sigma) = \lambda$.
(b) If $d = d_{H^{\prime}}$, there exists $\ell \in \mathbb{N}$ such
that for any $t \in A^{\ast}$ satisfying $|t| = \ell$, we have
$\bs{M}(\sigma) = \bs{M}_{c}(\sigma) = \lambda \cdot t$.
\end{proposition}

{\bf Proof.}
Obvious from Proposition~\ref{prop:shape_dist} and
Lemma~\ref{lemma:median_string_of_unimodal_and_symmetric_dist}.
\qed

With respect to parameter $\rho$ of $\LA(\lambda, \rho)$, we have the
following proposition:

\begin{proposition}
\label{prop:mean_abs_dev_of_Laplace-like_dist}
If a random string $\sigma$ is distributed according to $\LA(\lambda,
\rho)$, we have $\bs{\Upsilon}_{m}(\sigma) = \rho$ for any $d \in D$.
\end{proposition}

{\bf Proof.}
Setting $r = d(s, \lambda)$ and noting that the power series $\sum_{r
= 0}^{\infty} r (\rho / (\rho + 1))^{r}$ converges and its sum is
equal to $\rho (\rho + 1)$ from $\rho / (\rho + 1) < 1$, we obtain
\begin{eqnarray*}
& & \sum_{s \in A^{\ast}} d(s, \lambda) q(s; \lambda, \rho)
 \ = \ \frac{1}{\rho + 1} \sum_{s \in A^{\ast}} \frac{d(s, \lambda)}{\bigl|
\partial U(\lambda, d(s, \lambda)) \bigr|} \left( \frac{\rho}{\rho +
1} \right)^{d(s, \lambda)} \\
&=& \frac{1}{\rho + 1} \sum_{r = 0}^{\infty} \frac{r}{|\partial U(\lambda,
r)|} \left( \frac{\rho}{\rho + 1} \right)^{r} |\partial U(\lambda, r)|
 \ = \ \frac{1}{\rho + 1} \rho (\rho + 1)
 \ = \ \rho
\end{eqnarray*}
for any $d \in D$.
Thus, from Part~(a) of
Proposition~\ref{prop:consensus_seq_of_Laplace-like_dist}, we have
$\bs{\Upsilon}_{m}(\sigma) = \rho$ for any $d \in D$.
\qed

It is well known that among all continuous distributions with the
support $(0, \infty)$ the mean of which is equal to a given positive
real number, the exponential distribution has the maximum entropy.
Among all continuous distributions with the support $\mathbb{R}$ that
have a given mean and variance, the normal distribution has the
maximum entropy.
Similarly, the Laplace distribution maximizes the entropy among all
continuous distributions with the support $\mathbb{R}$ that satisfy the
condition that the first absolute moment about some fixed point is
equal to a given positive real number~\cite{Kagan_1973}.
The following proposition states that $\LA(\lambda, \rho)$ has a
similar property.

\begin{proposition}
\label{prop:maximum_entropy_A_ast}
Among all distributions on $A^{\ast}$ satisfying the condition that
the first absolute moment about some fixed string $m \in A^{\ast}$ is
equal to a given positive real number $v$, $\LA(m, v)$ maximizes the
entropy.
\end{proposition}

{\bf Proof.}
In this proof, we denote a value of a function $q$ on $A^{\ast}$ at
$s$ by $q_{s}$.
Although the constraints are
\[
(\mathrm{i}) \ q_{s} > 0, \forall s \in A^{\ast}, \ \ 
(\mathrm{ii}) \ \displaystyle{\sum_{s \in A^{\ast}}} q_{s} = 1, \ \ 
(\mathrm{iii}) \ \displaystyle{\sum_{s \in A^{\ast}}} d(s, m) q_{s} =
v,
\]
we first seek a function that maximizes the entropy among functions on
$A^{\ast}$ that satisfy constraints (ii) and (iii).
The Lagrangian is
\[
L = - \sum_{s \in A^{\ast}} q_{s} \log q_{s} - c_{1} \left( \sum_{s
\in A^{\ast}} q_{s} - 1 \right) - c_{2} \left( \sum_{s \in A^{\ast}} d(s,
m) q_{s} - v \right)
\]
for undetermined multipliers $c_{1}, c_{2} \neq 0$, and therefore, we
have $\partial L/\partial q_{t} = - \log q_{t} - 1 - c_{1} - c_{2}
d(t, m)$ for each $t \in A^{\ast}$.
Thus, the necessary condition to maximize the entropy under
constraints (ii) and (iii) is given by
\begin{equation}
q_{t} = \exp (- 1 - c_{1}) \exp (- c_{2} d(t, m)).
\label{eq:necessary-2}
\end{equation}
We set $d(t, m) = r$.
Noting that $c_{2} > 0$ from $q_{t} \leq 1$ and
Equation~(\ref{eq:necessary-2}) and that $\log((v + 1) / v) > 0$ from
$v > 0$ and making the parametrization of $c_{2} = \log ((v + 1) / v)$
provides
\begin{equation}
q_{t} = \exp (- 1 - c_{1}) \left( \frac{v}{v + 1} \right)^{r}.
\label{eq:parametrization_of_c_2}
\end{equation}
From $\sum_{r = 0}^{\infty} (v/v + 1)^{r} = 1/(1 - v / (v + 1)) = v +
1$, we have
\begin{equation}
\sum_{r = 0}^{\infty} \frac{1}{v + 1} \left( \frac{v}{v + 1}
\right)^{r} = 1.
\label{eq:sum_with_respect_to_r_is_equal_to_one}
\end{equation}
Equation~(\ref{eq:necessary-2}) holds for any $t \in A^{\ast}$ if and
only if $d(t, m) = d(t^{\prime}, m)$ implies $q_{t} = q_{t^{\prime}}$
for $t, t^{\prime} \in A^{\ast}$ because $c_{1}$ and $c_{2}$ are
constants and $q_{t}$ depends on $t$ only through $d(t, m)$.
Moreover, the number of $t^{\prime} \in A^{\ast}$ such that
$d(t^{\prime}, m) = r$ is equal to $|\partial U(m, r)|$ for $r \in
\mathbb{N}$.
Hence, noting Equations~(\ref{eq:parametrization_of_c_2})
and~(\ref{eq:sum_with_respect_to_r_is_equal_to_one}) and constraint
(ii), we obtain
\[
q_{t} = \frac{1}{(v + 1) \bigl| \partial U(m, d(t, m)) \bigr|} \left(
\frac{v}{v + 1} \right)^{d(t, m)}.
\]
The above $q_{t}$ also satisfies constraint (i).
Because the entropy is a concave function, its maximization subject to
linear constraints by Lagrange's method provides a global maximum.
The proof is completed.
\qed

\subsection{Estimation procedure of parameter $\lambda$ under the Levenshtein distance}
\label{subsection:Estimation_procedure_under_the_Levenshtein_distance}

In this subsection, we describe an estimation procedure of the
location parameter $\lambda$ of $\LA(\lambda, \rho)$ in the case of $d
= d_{L}$.

{\bf 1} \ Seek
\[
\check{\lambda}^{(0)} = \arg \min_{\lambda \in A^{\ast}} \sum_{i =
1}^{n} d_{L}(s_{i}, \lambda)
\]

\quad using an existing algorithm (for
example,~\cite{Kohonen_1985,Martinez_2000,Martinez_2003}).

{\bf 2} \ Compute
\[
\check{\rho}^{(0)} = \frac{1}{n} \sum_{i = 1}^{n} d_{L}(s_{i},
\check{\lambda}^{(0)}), \quad
F_{\ast}^{(0)} = F(\check{\lambda}^{(0)}, \check{\rho}^{(0)})
\]

\quad (see Equation~(\ref{eq:object_function}) for the definition of $F$).

{\bf 3} \ For $t = 1, 2, \cdots$,

\quad \quad {\bf 3.1} \ Set $\{ s_{\gamma} \in A^{\ast}: \gamma \in \Gamma^{(t - 1)} \} = \partial U(\check{\lambda}^{(t - 1)}, 1)$ and compute
\[
v_{\gamma} = \frac{1}{n} \sum_{i = 1}^{n} d_{L}(s_{i},
s_{\gamma}), \quad F(s_{\gamma}, v_{\gamma})
\]

\quad \quad \quad \quad for each $\gamma \in \Gamma^{(t - 1)}$.

\quad \quad {\bf 3.2} \ If there exists $\gamma \in \Gamma^{(t - 1)}$
such that $F(s_{\gamma}, v_{\gamma}) < F_{\ast}^{(t - 1)}$ holds, set
\begin{eqnarray*}
\gamma^{\ast} = \arg \min_{\gamma \in \Gamma^{(t - 1)}} F(s_{\gamma},
v_{\gamma}), \quad
\check{\lambda}^{(t)} = s_{\gamma^{\ast}}, \quad
\check{\rho}^{(t)} = v_{\gamma^{\ast}}, \quad
F_{\ast}^{(t)} = F(\check{\lambda}^{(t)}, \check{\rho}^{(t)}), \quad
t = t + 1
\end{eqnarray*}

\quad \quad \quad \quad and return to Step 3.1.
Otherwise, terminate the iteration and return
\[
\check{\lambda} = \check{\lambda}^{(t - 1)}, \quad
\check{\rho} = \check{\rho}^{(t - 1)}.
\]

\subsection{Strong consistency of maximum likelihood estimators in cases where the parameter space is $A^{\ast}$ and $A^{\ast} \times (0, \infty)$}
\label{subsection:Strong_consistency_of_MLE_in the_case_where_the_parameter_space_is_A^ast}

Maximum likelihood estimators for string parameters are strongly
consistent under quite general conditions, as are maximum likelihood
estimators for parameters that are real
numbers~\cite{Wald_1949,Perlman_1972}.
In this subsection, we describe propositions with respect to the
strong consistency of maximum likelihood estimators in cases where the
parameter space is $A^{\ast}$ and $A^{\ast} \times (0, \infty)$.
The proposition for the latter case underlies
Theorems~\ref{theorem:asymptotic_MLE}
and~\ref{theorem:strong_consistency_of_the_EM_algorithm} and
consequently also Theorem~\ref{theorem:asymptotic_EM} and
Corollary~\ref{corollary:optimality_of_clustering_procedure}.

Let $\sigma \in \mathcal{M}(\Omega, A^{\ast})$ have a probability
function $q(s; \theta)$ with a parameter $\theta \in A^{\ast}$.
$\theta^{\ast}$ represents the true value of the parameter.
We set
\begin{equation}
\eta(\theta^{\prime}, \theta) = \sum_{s \in A^{\ast}} \log (q(s;
\theta)) q(s; \theta^{\prime}) \label{eq:entropic_divergence}
\end{equation}
for any $\theta, \theta^{\prime} \in A^{\ast}$.
It is verified that $\eta(\theta^{\prime}, \theta) \leq
\eta(\theta^{\prime}, \theta^{\prime})$ holds for any $\theta \in
A^{\ast}$ in the same manner as in the case where the sample and
parameter spaces are $\mathbb{R}$.

\begin{proposition}
\label{prop:strong_consistency_of_MLE_on_A^ast}
We suppose that $\sigma_{1}, \cdots, \sigma_{n} \in
\mathcal{M}(\Omega, A^{\ast})$ (i) are independent and (ii) have the
identical probability function $q(s; \theta)$ with a parameter $\theta
\in A^{\ast}$ and that (iii) $q(s; \theta)$ has the support
$A^{\ast}$.
We denote the realization of $\sigma_{i}$ by $s_{i}$ for each $i = 1,
\cdots, n$ and the maximum likelihood estimator of $\theta$ based on
$s_{1}, \cdots, s_{n}$ by $\check{\theta}^{(n)}$.
If (iv) $\eta(\theta^{\ast}, \theta) < \eta(\theta^{\ast},
\theta^{\ast})$ holds for any $\theta \in A^{\ast} \smallsetminus \{
\theta^{\ast} \}$, there exists $N_{0} \in \mathbb{Z}^{+}$ such
that if $n \geq N_{0}$, we have $\check{\theta}^{(n)} =
\theta^{\ast}$ a.s.
\end{proposition}

{\bf Proof.}
We denote the log likelihood function of $\theta$ based on $s_{1},
\cdots, s_{n}$ by $\ell(\theta) = \ell(\theta; s_{1}, \cdots, s_{n})$.
From condition~(ii), we have
\begin{eqnarray}
\E_{\theta^{\ast}} \left[ \frac{1}{n} \ell(\theta) \right]
&=& \frac{1}{n} \sum_{i = 1}^{n} \E_{\theta^{\ast}} [\log q(\sigma_{i};
\theta)]
= \E_{\theta^{\ast}} [\log q(\sigma_{1}; \theta)] \nonumber \\
&=& \sum_{s_{1} \in A^{\ast}} \log (q(s_{1}; \theta)) q(s_{1};
\theta^{\ast})
= \eta (\theta^{\ast}, \theta)
\label{eq:expected_mean_log_likelihood}
\end{eqnarray}
for any $\theta \in A^{\ast}$.
Noting that $0 < q(s_{i}; \theta) \leq 1$ holds from condition~(iii),
we find that
\[
-\infty < \log q(s_{i}; \theta) \leq 0, \quad -\infty < \sum_{s_{i}
\in A^{\ast}} \log (q(s_{i}; \theta)) q(s_{i}; \theta) \leq 0.
\]
Thus,
\begin{eqnarray*}
\Var_{\theta^{\ast}} [\log q(\sigma_{i}; \theta)]
&=& \E_{\theta^{\ast}} \left[ \left\{ \log q(\sigma_{i}; \theta) -
\E_{\theta^{\ast}} [\log q(\sigma_{i}; \theta) ] \right\}^{2} \right] \\
&=& \sum_{s_{i} \in A^{\ast}} \left\{ \log q(s_{i}; \theta) - \sum_{s_{i}
\in A^{\ast}} \log (q(s_{i}; \theta)) q(s_{i}; \theta) \right\}^{2}
q(s_{i}; \theta) \ < \ \infty
\end{eqnarray*}
holds.
Therefore, using the strong law of large numbers from conditions~(i)
and~(ii) and noting Equation~(\ref{eq:expected_mean_log_likelihood}),
we obtain
\begin{equation}
\frac{1}{n} \ell(\theta)
= \frac{1}{n} \sum_{i = 1}^{n} \log q(s_{i}; \theta)
\asconv \frac{1}{n} \sum_{i = 1}^{n} \E_{\theta^{\ast}} [\log q(\sigma_{i};
\theta)]
= \E_{\theta^{\ast}} \left[ \frac{1}{n} \ell (\theta) \right]
= \eta(\theta^{\ast}, \theta)
\label{eq:mean_log_likelihood_conv}
\end{equation}
as $n \lto \infty$ for any $\theta \in A^{\ast}$.
We set
\[
\delta = \min_{\theta \in A^{\ast} \smallsetminus \{ \theta^{\ast} \}}
| \eta(\theta^{\ast}, \theta) - \eta(\theta^{\ast},
\theta^{\ast}) |.
\]
We have $\delta > 0$ from condition~(iv).
Hence, for any $\theta \in A^{\ast}$ there exists $N_{\theta} \in
\mathbb{Z}^{+}$ such that if $n \geq N_{\theta}$, then
\[
\left| \frac{1}{n} \ell(\theta) - \eta(\theta^{\ast}, \theta) \right|
< \frac{\delta}{2} \as
\]
holds from Equation~(\ref{eq:mean_log_likelihood_conv}).
Thus, if $n \geq N_{0}$ for $N_{0} = \max \{ N_{\theta}: \theta \in
A^{\ast} \}$, we have
\[
\max_{\theta \in A^{\ast}} \left| \frac{1}{n} \ell(\theta) -
\eta(\theta^{\ast}, \theta) \right| < \frac{\delta}{2} \as
\]
Therefore,
\[
\left| \frac{1}{n} \ell(\theta^{\ast}) - \eta(\theta^{\ast},
\theta^{\ast}) \right| < \frac{\delta}{2}, \quad \left| \frac{1}{n}
\ell(\theta) - \eta(\theta^{\ast}, \theta) \right| < \frac{\delta}{2}
\mbox{ for any } \theta \neq \theta^{\ast}
\]
hold.
From the definition of $\delta$, the above inequalities mean that
$\ell(\theta) / n$ attains a maximum value at $\theta = \theta^{\ast}$
for any $n \geq N_{0}$.
The maximizer of $\ell(\theta) / n$ is $\check{\theta}^{(n)}$.
Hence, we obtain $\check{\theta}^{(n)} = \theta^{\ast}$ a.s. for any
$n \geq N_{0}$ from condition~(iv).
\qed

We suppose that $\sigma \in \mathcal{M}(\Omega, A^{\ast})$ has a
probability function $q(s; \bs{\theta})$ with a parameter $\bs{\theta}
= (\theta_{1}, \theta_{2}) \in A^{\ast} \times (0, \infty)$.
$\bs{\theta}^{\ast} = (\theta_{1}^{\ast}, \theta_{2}^{\ast})$
represents the true value of the parameter.
$\eta(\bs{\theta}^{\prime}, \bs{\theta})$ is defined for any
$\bs{\theta}, \bs{\theta}^{\prime} \in A^{\ast} \times (0, \infty)$ by
Equation~(\ref{eq:entropic_divergence}).
$\eta(\bs{\theta}^{\prime}, \bs{\theta}) \leq
\eta(\bs{\theta}^{\prime}, \bs{\theta}^{\prime})$ holds for any
$\bs{\theta} \in A^{\ast}$.
We introduce the following regular conditions, which are obtained by
slightly modifying the regular conditions for the strong consistency
of maximum likelihood estimators in the case where the parameter
space is $\mathbb{R}$:

{\bf 1.} If $|\theta_{2} - \theta_{2}^{\ast}| > 0$, we have
$\eta(\bs{\theta}^{\ast}, \bs{\theta}^{\ast}) -
\eta(\bs{\theta}^{\ast}, \bs{\theta}) > 0$.

{\bf 2.} Setting
\[
g_{M}(s) = \sup_{\substack{\theta_{1} \neq \theta_{1}^{\ast} \lor\\
|\theta_{2} - \theta_{2}^{\ast}| > M}} \log q(s; \bs{\theta})
\]
for $M > 0$, we have
\[
c_{g} = \E_{\bs{\theta}^{\ast}} [g_{M}(s)] < \eta(\bs{\theta}^{\ast},
\bs{\theta}^{\ast})
\]
for a sufficiently large $M$.

{\bf 3.} $q(s; \bs{\theta})$ is partially differentiable with respect
to $\theta_{2}$ for any $s \in A^{\ast}$, and setting
\[
h_{M}(s) = \sup_{|\theta_{2} - \theta_{2}^{\ast}| \leq M} \left|
\frac{\partial}{\partial \theta_{2}} \log q(s; \bs{\theta}) \right|
\]
for $M$ for which regular condition~2 holds, we have
\[
c_{h} = \E_{\bs{\theta}^{\ast}} [h_{M}(s)] < \infty.
\]

\begin{proposition}
\label{prop:strong_consistency_of_MLE_on_A^ast_times_R}
We suppose that $\sigma_{1}, \cdots, \sigma_{n} \in
\mathcal{M}(\Omega, A^{\ast})$ (i) are independent and (ii) have the
identical probability function $q(s; \bs{\theta})$
and (iii) $q(s; \bs{\theta})$ has the support $A^{\ast}$.
We denote the realization of $\sigma_{i}$ by $s_{i}$ for each $i = 1,
\cdots, n$ and the maximum likelihood estimator of $\bs{\theta}$ based
on $s_{1}, \cdots, s_{n}$ by $\check{\bs{\theta}}^{(n)} =
(\check{\theta}_{1}^{(n)}, \check{\theta}_{2}^{(n)})$.
If regular conditions~1 to~3 are satisfied, there exists $N_{0} \in
\mathbb{Z}^{+}$ such that if $n \geq N_{0}$, we have
$\check{\theta}_{1}^{(n)} = \theta_{1}^{\ast}$ a.s. and
$\check{\theta}_{2}^{(n)} \asconv \theta_{2}^{\ast}$ as $n \lto
\infty$.
\end{proposition}

{\bf Proof.}
We denote the log likelihood function of $\bs{\theta}$ based on
$s_{1}, \cdots, s_{n}$ by $\ell(\bs{\theta}) = \ell(\bs{\theta};
s_{1}, \cdots, s_{n})$.
From conditions (i) and (ii),
\begin{eqnarray}
\E_{\bs{\theta}^{\ast}} \left[ \frac{1}{n} \ell(\bs{\theta}) \right]
&=& \frac{1}{n} \sum_{i = 1}^{n} \E_{\bs{\theta}^{\ast}} [\log q(\sigma_{i};
\bs{\theta})]
 \ = \ \E_{\bs{\theta}^{\ast}} [\log q(\sigma_{1}; \bs{\theta})]
\nonumber \\
&=& \sum_{s_{1} \in A^{\ast}} \log (q(s_{1}; \bs{\theta})) q(s_{1};
\bs{\theta}^{\ast})
 \ = \ \eta (\bs{\theta}^{\ast}, \bs{\theta})
\label{eq:expected_mean_log_likelihood_conv}
\end{eqnarray}
holds for any $\bs{\theta} \in A^{\ast} \times (0, \infty)$.
We have
\[
-\infty < \log q(s_{i}; \bs{\theta}) \leq 0, \quad -\infty < \sum_{s_{i}
\in A^{\ast}} \log (q(s_{i}; \bs{\theta})) q(s_{i}; \bs{\theta}) \leq 0
\]
from condition (iii).
Thus,
\begin{eqnarray*}
\Var_{\bs{\theta}^{\ast}} [\log q(\sigma_{i}; \bs{\theta})]
&=& \E_{\bs{\theta}^{\ast}} \left[ \left\{ \log q(\sigma_{i}; \bs{\theta}) -
\E_{\bs{\theta}^{\ast}} [\log q(\sigma_{i}; \bs{\theta}) ] \right\}^{2}
\right] \\
&=& \sum_{s_{i} \in A^{\ast}} \left\{ \log q(s_{i}; \bs{\theta}) - \sum_{s_{i}
\in A^{\ast}} \log (q(s_{i}; \bs{\theta})) q(s_{i}; \bs{\theta})
\right\}^{2} q(s_{i}; \bs{\theta}) \ < \ \infty
\end{eqnarray*}
holds.
Therefore, using the strong law of large numbers from conditions~(i)
and~(ii) and noting
Equation~(\ref{eq:expected_mean_log_likelihood_conv}), we obtain
\[
\frac{1}{n} \ell(\bs{\theta})
= \frac{1}{n} \sum_{i = 1}^{n} \log q(s_{i}; \bs{\theta})
\asconv \frac{1}{n} \sum_{i = 1}^{n} \E_{\bs{\theta}^{\ast}} [\log
q(\sigma_{i}; \bs{\theta})]
= \E_{\bs{\theta}^{\ast}} \left[\frac{1}{n} \ell (\bs{\theta}) \right]
= \eta(\bs{\theta}^{\ast}, \bs{\theta})
\]
as $n \lto \infty$ for any $\bs{\theta} \in A^{\ast} \times (0, \infty)$.
Consequently,
\begin{equation}
\frac{1}{n} \ell(\bs{\theta}^{\ast}) \asconv \eta(\bs{\theta}^{\ast},
\bs{\theta}^{\ast}) \label{eq:result_of_Step_1}
\end{equation}
as $n \lto \infty$.
By the definition of $c_{g}$ and the strong law of large numbers, we
have
\begin{eqnarray*}
\frac{1}{n} \sup_{\substack{\theta_{1} \neq \theta_{1}^{\ast} \lor\\
|\theta_{2} - \theta_{2}^{\ast}| > M}} \ell(\bs{\theta})
&=& \frac{1}{n} \sup_{\substack{\theta_{1} \neq \theta_{1}^{\ast} \lor\\
|\theta_{2} - \theta_{2}^{\ast}| > M}} \sum_{i = 1}^{n} \log q(s_{i};
\bs{\theta})\\
&\leq& \frac{1}{n} \sum_{i = 1}^{n} \sup_{\substack{\theta_{1} \neq
\theta_{1}^{\ast} \lor\\ |\theta_{2} - \theta_{2}^{\ast}| > M}}
\log q(s_{i}; \bs{\theta})
 \ = \ \frac{1}{n} \sum_{i = 1}^{n} g_{M}(s)
 \ \asconv \ c_{g}
\end{eqnarray*}
as $n \lto \infty$.
Hence, from regular condition~2,
\[
\frac{1}{n} \sup_{\substack{\theta_{1} \neq \theta_{1}^{\ast}
\lor\\ |\theta_{2} - \theta_{2}^{\ast}| > M}} \ell
(\bs{\theta}) < \eta(\bs{\theta}^{\ast}, \bs{\theta}^{\ast}) \as
\]
as $n \lto \infty$.
Thus, noting Equation~(\ref{eq:result_of_Step_1}), we observe that the
maximizer of $\ell(\bs{\theta})$ in $A^{\ast} \times (0, \infty)$,
i.e., the maximum likelihood estimate $\check{\bs{\theta}}^{(n)}$ of
$\bs{\theta}$ satisfies
\[
\check{\theta}_{1}^{(n)} = \theta_{1}^{\ast} \quad \mbox{and} \quad
\check{\theta}_{2}^{(n)} \in [-M, M] \as
\]
as $n \lto \infty$.
The almost sure convergence of $\check{\theta}_{2}^{(n)}$ to
$\theta_{2}^{\ast}$ is demonstrated using regular conditions~1 and~3
in the same manner as in the proof of the strong consistency of the
maximum likelihood estimator in the case where the parameter space is
$\mathbb{R}$.
\qed

Regular conditions 1 to 3 are quite general, and $\LA(\lambda,
\rho)$ satisfies these conditions as well as condition (iii) of
Proposition~\ref{prop:strong_consistency_of_MLE_on_A^ast_times_R}.

\subsection{Convergence of a sequence of algorithms to an algorithm}
\label{subsection:convergence_of_a_sequence_of_algorithms}

In this subsection, we define the convergence of a sequence of
algorithms to an algorithm.
This concept is used in describing Theorem~\ref{theorem:asymptotic_EM}
in
Section~\ref{section:EM_algorithm_for_the_Laplace-like_mixture_on_A^ast}.

Let $X$ and $Y$ be input and output spaces, respectively, and $d_{Y}$
be a distance function on $Y$.
Here, as algorithms we consider a sequence of computations that
returns an output $y \in Y$ for each input $x \in X$ and for which the
number $\ell(x)$ of computations required to return $y$ and the
$j$-th computation $c_{x}^{(j)}$ for each $j = 1, \cdots, \ell(x)$ can
vary depending on $x$.
In other words, for a function $\ell: X \to \mathbb{Z}^{+}$ and
computations $c_{x}^{(1)}, \cdots, c_{x}^{(\ell(x))}$ determined for
each $x \in X$, we consider algorithms represented as a composition $C:
X \to Y, \ x \mapsto c_{x}^{(\ell(x))} \circ \cdots \circ
c_{x}^{(1)}(x)$ of these computations.
Let $k$ be a $\mathbb{Z}^{+}$-valued function defined on
$\mathbb{Z}^{+} \times X$ and $B_{n} = b_{n, x}^{(k(n, x))} \circ
\cdots \circ b_{n, x}^{(1)}$ be an algorithm with the input
and output spaces $X$ and $Y$ for each $n \in \mathbb{Z}^{+}$.
We consider a sequence $\{ B_{n}: n \in \mathbb{Z}^{+ }\}$ of
algorithms with the parameter $n \in \mathbb{Z}^{+}$.

\begin{definition}
\label{def:convergence_of_a_sequence_of_algorithms}
We say that $\{ B_{n}: n \in \mathbb{Z}^{+} \}$ converges to $C$ as $n
\lto \infty$ if for any $x \in X$ (i) there exists $N_{0} \in
\mathbb{Z}^{+}$ such that if $n \geq N_{0}$, then $k(n, x) =
\ell(x)$ holds and (ii) $d_{Y}(b_{n, x}^{(k(n, x))} \circ \cdots \circ
b_{n, x}^{(1)}(x), c_{x}^{(\ell(x))} \circ \cdots \circ
c_{x}^{(1)}(x)) \lto 0$ as $n \lto \infty$.
\end{definition}

In the case where the parameter space is $[0, \infty)$, the
convergence of $\{ B_{n} \}$ to $C$ is defined in a similar manner.
Furthermore, if $X$ is a probability space, the convergence in
probability and the almost sure convergence of $\{ B_{n} \}$ to $C$
are defined in a trivial manner based on
Definition~\ref{def:convergence_of_a_sequence_of_algorithms}.
We can also define the convergence of $\{ B_{n} \}$ to $C$ in a
stronger manner by replacing condition~(ii) in
Definition~\ref{def:convergence_of_a_sequence_of_algorithms} with
(ii$^{\prime}$) for each $j = 1, \cdots, \ell(x)$ there exists a
metric space $(Y_{j}, d_{j})$ such that $b_{n, x}^{(j)} \circ \cdots
\circ b_{n, x}^{(1)}(x), c_{x}^{(j)} \circ \cdots \circ c_{x}^{(1)}(x)
\in Y_{j}$ holds and $d_{j}(b_{n, x}^{(j)} \circ \cdots \circ b_{n,
x}^{(1)}(x), c_{x}^{(j)} \circ \cdots \circ c_{x}^{(1)}(x)) \lto 0$
as $n \lto \infty$.

\subsection{Proofs of the results}
\label{subsection:Proofs_of_the_results}

In this subsection, proofs of the results described in the main text
are provided.

{\bf Proof of Proposition and definition~\ref{propanddef:Laplace-like_distribution_on_A^ast}.}
The nonnegativity is trivial.
Noting that there exist $| \partial U(\lambda, r)|$ strings in
$A^{\ast}$ that satisfy $r = d(s, \lambda)$ and that $\sum_{r =
0}^{\infty}(\rho / (\rho + 1))^{r} = \rho + 1$ holds from $\rho /
(\rho + 1) < 1$, we obtain
\begin{eqnarray*}
& & \sum_{s \in A^{\ast}} q(s; \lambda, \rho)
 \ = \ \frac{1}{\rho + 1} \sum_{s \in A^{\ast}} \frac{1}{\bigl|
\partial U(\lambda, d(s, \lambda)) \bigr|} \left( \frac{\rho}{\rho
+ 1} \right)^{d(s, \lambda)} \\
&=& \frac{1}{\rho + 1} \sum_{r = 0}^{\infty} \frac{1}{| \partial
U(\lambda, r)|} \left( \frac{\rho}{\rho + 1} \right)^{r} | \partial
U(\lambda, r)|
 \ = \ \frac{1}{\rho + 1} (\rho + 1) \ = \ 1.
\end{eqnarray*}
\qed

{\bf Proof of Lemma~\ref{lemma:estimating_a_mode_string_of_a_unimodal_symmetric_distribution}.}
Corollary~4.2 in~\cite{Koyano_2014b} holds for any $d \in D$ (the
proof described in~\cite{Koyano_2014b} works for any $d \in D$).
We denote the relative frequency of $x \in \bar{A}$ at the $j$-th site
of $s_{1}, \cdots, s_{n} \in A^{\ast}$ by $f_{j}^{(n)}(x)$ for each $j
\in \mathbb{Z}^{+}$.
Let $m \cdot t = \{ m_{1}, \cdots, m_{|m|}, t_{1}, \cdots, t_{\ell},
e, \cdots \}$.

{\bf (Step 1)} We first consider the case of $|m| \geq 1$.
We set
\[
x_{j}^{\ast} = \arg \max_{x \in \bar{A} \smallsetminus \{ m_{j} \}}
q_{j}(x), \quad \delta_{j} = q_{j}(m_{j}) - q_{j}(x_{j}^{\ast})
\]
for each $j = 1, \cdots, |m|$,
\[
\delta_{j} = q_{j}(a_{1}) - q_{j}(e)
\]
for each $j = |m| + 1, \cdots, |m| + \ell$ if $\ell \geq 1$, and
\[
x_{|m| + \ell + 1}^{\ast} = \arg \max_{x \in A} q_{|m| + \ell + 1}(x),
\quad \delta_{|m| + \ell + 1} = q_{|m| + \ell + 1}(e) - q_{|m| + \ell
+ 1}(x_{|m| + \ell + 1}^{\ast}).
\]
From condition~(ii), $\delta_{j} > 0$ holds for each $j = 1, \cdots,
|m| + \ell + 1$.
We put $\epsilon_{j} = \delta_{j} / 3$.
By condition~(iii) and the strong law of large numbers, there exists
$N_{0} \in \mathbb{Z}^{+}$ such that if $n \geq N_{0}$, we have
\[
\left| f_{j}^{(n)}(m_{j}) - q_{j}(m_{j}) \right| < \epsilon_{j}, \quad
\left| f_{j}^{(n)}(e) - q_{j}(e) \right| <
\epsilon_{j} \as
\]
for each $j = 1, \cdots, |m|$,
\[
\left| f_{j}^{(n)}(a_{1}) - q_{j}(a_{1}) \right| < \epsilon_{j}, \quad
\left| f_{j}^{(n)}(e) - q_{j}(e) \right| < \epsilon_{j} \as
\]
for each $j = |m| + 1, \cdots, |m| + \ell$ if $\ell \geq 1$, and
\[
\left| f_{|m| + \ell + 1}^{(n)}(e) - q_{|m| + \ell + 1}(e) \right| <
\epsilon_{|m| + \ell + 1}, \quad \left| f_{|m| + \ell +
1}^{(n)}(x_{|m| + \ell + 1}^{\ast}) - q_{|m| + \ell + 1} (x_{|m| +
\ell + 1}^{\ast}) \right| < \epsilon_{|m| + \ell + 1} \as
\]
Thus, noting the definition of $\epsilon_{j}$, if $n \geq N_{0}$, then
$f_{j}^{(n)}(m_{j}) > f_{j}^{(n)}(e)$ a.s. for each $j = 1, \cdots,
|m|$, $f_{j}^{(n)}(a_{1}) > f_{j}^{(n)}(e)$ a.s. for each $j = |m| +
1, \cdots, |m| + \ell$ if $\ell \geq 1$, and $f_{|m| + \ell +
1}^{(n)}(e) > f_{|m| + \ell + 1}^{(n)}(x_{|m| + \ell + 1}^{\ast})$
a.s. hold.
Consequently, from Equation~(\ref{eq:def_j^ast}), we obtain
\begin{equation}
j^{\ast} = |m| + \ell \as \label{eq:j^ast}
\end{equation}
for $n \geq N_{0}$ in both cases of $\ell = 0$ and $\ell \geq 1$.

We set $\epsilon_{0} = \min_{1 \leq j \leq |m|} \epsilon_{j}$.
By condition~(iii) and the strong law of large numbers, for any
$\epsilon \leq \epsilon_{0}$ there exists $N_{1} \in \mathbb{Z}^{+}$
such that if $n \geq N_{1}$, then
\[
\left| f_{j}^{(n)}(m_{j}) - q_{j}(m_{j}) \right| < \epsilon, \quad
\left| f_{j}^{(n)}(x_{j}^{\ast}) - q_{j}(x_{j}^{\ast}) \right| <
\epsilon \as
\]
hold for each $j = 1, \cdots, |m|$.
Therefore, noting the definition of $\epsilon_{0}$, we have
\begin{equation}
\left| f_{j}^{(n)}(m_{j}) - f_{j}^{(n)}(x_{j}^{\ast}) \right|
> \epsilon \as \label{eq:f_j_m_j^prime_minux_f_j_x_j^ast_<_epsilon}
\end{equation}
for $n \geq N_{1}$.
We consider the case of $\ell \geq 1$.
Noting condition~(ii) and using the strong law of large numbers from
condition~(iii), there exists $N_{2} \in \mathbb{Z}^{+}$ such that if
$n \geq N_{2}$, then
\begin{equation}
\max_{x \in A} f_{j}^{(n)}(x) - \min_{x \in A} f_{j}^{(n)}(x) <
\epsilon \as \label{eq:max_minus_min_<_epsilon}
\end{equation}
holds for each $j = |m| + 1, \cdots, |m| + \ell$.
Hence, noting the definition of condition U$^{(\epsilon)}$ and
Equation~(\ref{eq:def_j^dagger}), we have
\begin{equation}
j^{(\epsilon)} = |m| \ (\geq 1) \as \label{eq:j^epsilon}
\end{equation}
for $n \geq \max \{ N_{0}, N_{1}, N_{2} \}$.
From Equations~(\ref{eq:j^ast}) and~(\ref{eq:j^epsilon}),
$\bs{m}_{c}^{(\epsilon)}(s_{1}, \cdots, s_{n})$ is determined
according to the lower of Equation~(\ref{eq:def_m_c^epsilon}) in the
case of $\ell \geq 1$ and, from the definition of $\eta(j)$, is a
sequence obtained by concatenating a substring composed of the first
to the $|m|$-th letters of a consensus sequence of $s_{1}, \cdots,
s_{n}$ with an infinite sequence of $e$s.
In the case of $\ell = 0$, U$^{(\epsilon)}$ does not hold from
Equation~(\ref{eq:f_j_m_j^prime_minux_f_j_x_j^ast_<_epsilon}).
Thus, $\bs{m}_{c}^{(\epsilon)}(s_{1}, \cdots, s_{n})$ is determined
according to the middle of Equation~(\ref{eq:def_m_c^epsilon}) and is
the same sequence as in the case of $\ell \geq 1$.
Therefore, noting conditions~(iii) and~(iv) and applying Corollary~4.2
in~\cite{Koyano_2014b} to only the substring composed of the first to
the $|m|$-th letters of $\bs{m}_{c}^{(\epsilon)}(s_{1}, \cdots,
s_{n})$, if $n \geq \max \{ N_{0}, N_{1}, N_{2} \}$, then
\begin{equation}
\bs{m}_{c}^{(\epsilon)}(s_{1}, \cdots, s_{n}) = \{ m_{1}, \cdots,
m_{|m|}, e, \cdots \} = m
\as \label{eq:case_of_|m|_greater_or_equal_one}
\end{equation}
holds.

{\bf (Step 2)}
We next consider the case of $|m| = 0$.
We set
\[
x_{j}^{\ast} = \arg \max_{x \in A} q_{j}(x), \quad \delta_{j} =
q_{j}(e) - q_{j}(x^{\ast}_{j}), \quad \epsilon_{j} =
\frac{\delta_{j}}{3}
\]
for each $j \geq \ell + 1$ and put $\epsilon_{0} = \min_{j \geq \ell +
1} \epsilon_{j}$.
Conducting a similar discussion to that in Step~1 noting
condition~(ii) and using the strong law of large numbers from
condition~(iii), there exists $N_{3} \in \mathbb{Z}^{+}$ such that if
$n \geq N_{3}$, then $f_{j}^{(n)}(e) > f_{j}^{(n)}(x_{j}^{\ast})$
a.s. holds for each $j \geq \ell + 1$.
Consequently, in the case of $\ell = 0$, we have $j^{\ast} = 0$
a.s. for $n \geq N_{3}$ from Equation~(\ref{eq:def_j^ast}).
On the other hand, in the case of $\ell \geq 1$, $j^{\ast} \geq 1$
a.s. holds.
Furthermore, conducting a similar discussion to that in Step~1, for
any $\epsilon \leq \epsilon_{0}$, there exists $N_{4} \in
\mathbb{Z}^{+}$ such that if $n \geq N_{4}$,
Equation~(\ref{eq:max_minus_min_<_epsilon}) holds for each $j = 1,
\cdots, \ell$.
Hence, we have $j^{(\epsilon)} = 0$ a.s. for $n \geq \max \{ N_{3},
N_{4} \}$.
Thus, in both cases of $\ell = 0$ and $\ell \geq 1$,
$\bs{m}_{c}^{(\epsilon)}(s_{1}, \cdots, s_{n})$ is determined
according to the upper of Equation~(\ref{eq:def_m_c^epsilon}), and if
$n \geq \max \{ N_{3}, N_{4} \}$ and $\epsilon \leq \epsilon_{0}$,
then
\begin{equation}
\bs{m}_{c}^{(\epsilon)}(s_{1}, \cdots, s_{n}) = o = m
\as \label{eq:case_of_|m|_equal_zero}
\end{equation}
holds.
Combining Equations~(\ref{eq:case_of_|m|_greater_or_equal_one})
and~(\ref{eq:case_of_|m|_equal_zero}) completes the proof.
\qed

{\bf Proof of Proposition~\ref{prop:MLE_location_parameter}.}
Obvious from
Lemma~\ref{lemma:estimating_a_mode_string_of_a_unimodal_symmetric_distribution}, Proposition~\ref{prop:shape_dist},
and Lemma~\ref{lemma:median_string_of_unimodal_and_symmetric_dist}.
\qed

{\bf Proof of Proposition~\ref{prop:MLE_dispersion_parameter}.}
We demonstrate Part~(b).
Condition~(iv) of
Lemma~\ref{lemma:estimating_a_mode_string_of_a_unimodal_symmetric_distribution}
is satisfied from condition~(i).
Thus, noting conditions~(ii$^{\prime}$) and~(iii$^{\prime}$) and using
Proposition~\ref{prop:MLE_location_parameter}, there exist $N_{0} \in
\mathbb{Z}^{+}$ and $\epsilon_{0} > 0$ such that if $n \geq N_{0}$ and
$\epsilon \leq \epsilon_{0}$, we have $\bs{m}_{c}^{(\epsilon)}(s_{1},
\cdots, s_{n}) = \lambda$ a.s., and therefore,
\begin{equation}
\frac{1}{n} \sum_{i = 1}^{n} d_{H^{\prime}}(s_{i},
\bs{m}_{c}^{(\epsilon)}(s_{1}, \cdots, s_{n})) = \frac{1}{n} \sum_{i =
1}^{n} d_{H^{\prime}}(s_{i}, \lambda) \as
\label{eq:truncated_consensus_seq_equal_lambda}
\end{equation}
From conditions~(ii$^{\prime}$) and~(iii$^{\prime}$),
$d_{H^{\prime}}(s_{1}, \lambda), \cdots, d_{H^{\prime}}(s_{n},
\lambda)$ are also independent.
Using Proposition~\ref{prop:mean_abs_dev_of_Laplace-like_dist} from
conditions~(i) and~(ii$^{\prime}$), we have $\E_{\lambda, \rho}
[d_{H^{\prime}}(\sigma_{i}, \lambda)] = \bs{\Upsilon}_{m}(\sigma_{i}) =
\rho$ for each $i = 1, \cdots, n$.
Furthermore,
\[
\E_{\lambda, \rho} \left[ \Bigl\{ d_{H^{\prime}}(\sigma_{i}, \lambda)
- \E_{\lambda, \rho} [ d_{H^{\prime}}(\sigma_{i}, \lambda) ]
\Bigr\}^{2} \right] = \sum_{s \in A^{\ast}} \{ d_{H^{\prime}}(s,
\lambda) - \rho \}^{2} q(s; \lambda, \rho) < \infty
\]
holds.
Hence, by the strong law of large numbers, we obtain
\begin{equation}
\frac{1}{n} \sum_{i = 1}^{n} d_{H^{\prime}}(s_{i}, \lambda) \asconv
\frac{1}{n} \sum_{i = 1}^{n} \E_{\lambda, \rho}
[d_{H^{\prime}}(\sigma_{i}, \lambda)] = \rho
\label{eq:sample_mean_absolute_dev_equal_rho}
\end{equation}
as $n \lto \infty$.
Combining Equations~(\ref{eq:truncated_consensus_seq_equal_lambda})
and~(\ref{eq:sample_mean_absolute_dev_equal_rho}) provides $\sum_{i =
1}^{n} d_{H^{\prime}}(s_{i}, \bs{m}_{c}^{(\epsilon)}(s_{1}, \cdots,
s_{n})) / n \asconv \rho$ as $n \lto \infty$ and $\epsilon \lto 0$.
Part (a) is demonstrated in a similar manner.
\qed

{\bf Proof of
Proposition~\ref{prop:mean_abs_dev_and_MLE_dispersion}.}
We suppose that $n$ strings $s_{1}, \cdots, s_{n}$ are observed
and set $r_{i} = d(s_{i}, \lambda)$ for each $i = 1, \cdots, n$.
The log likelihood function of $\lambda$ and $\rho$ is given by
\[
\ell(\lambda, \rho; s_{1}, \cdots, s_{n}) = - n \log (\rho + 1) -
\sum_{i = 1}^{n} \log \bigl| \partial U(\lambda, d(s_{i}, \lambda))
\bigr| + \log \left( \frac{\rho}{\rho + 1} \right) \sum_{i = 1}^{n}
d(s_{i}, \lambda).
\]
Solving
\[
\frac{\partial}{\partial \rho} \ell(\lambda, \rho; s_{1}, \cdots, s_{n})
= - \frac{n}{\rho + 1} + \frac{1}{\rho (\rho + 1)} \sum_{i = 1}^{n}
d(s_{i}, \lambda) = 0
\]
with respect to $\rho$, we obtain
\[
\rho^{\ast} = \frac{1}{n} \sum_{i = 1}^{n} d(s_{i}, \lambda).
\]
We have
\[
\frac{\partial^{2}}{\partial \rho^{2}} \ell(\lambda, \rho; s_{1}, \cdots,
s_{n}) = \frac{n \rho^{2} - 2 R \rho - R}{\rho^{2} (\rho + 1)^{2}}
\]
for $R = \sum_{i = 1}^{n} d(s_{i}, \lambda)$.
Noting
\[
\frac{\partial^{2}}{\partial \rho^{2}} \ell(\lambda, \rho; s_{1}, \cdots,
s_{n}) < 0 \Longleftrightarrow \frac{R}{n} -
\frac{\sqrt{R (R + n)}}{n} < \rho < \frac{R}{n} +
\frac{\sqrt{R (R + n)}}{n},
\]
$\partial^{2} \ell(\lambda, \rho^{\ast}; s_{1}, \cdots, s_{n}) /
\partial \rho^{2} < 0$ holds.
Therefore, if $\lambda$ is known, the maximum likelihood
estimator of $\rho$ is given by
\[
\check{\rho}(s_{1}, \cdots, s_{n}) = \frac{1}{n} \sum_{i = 1}^{n}
d(s_{i}, \lambda).
\]
If $\lambda$ is unknown, the maximum likelihood estimator of $\rho$ is
obtained by replacing $\lambda$ on the right-hand side of the above
equation with its maximum likelihood estimator.
\qed

{\bf Proof of Theorem~\ref{theorem:asymptotic_MLE}.}
Conditions~(i) and~(ii) imply conditions~(iii) and~(iv) of
Lemma~\ref{lemma:estimating_a_mode_string_of_a_unimodal_symmetric_distribution},
respectively.
Therefore, applying Proposition~\ref{prop:MLE_location_parameter} from
condition~(iii), there exist $N_{0} \in \mathbb{Z}^{+}$ and
$\epsilon_{0} > 0$ such that if $n \geq N_{0}$ and $\epsilon \leq
\epsilon_{0}$, then
\begin{equation}
\bs{m}_{c}^{(\epsilon)}(s_{1}, \cdots, s_{n}) = \lambda \as
\label{eq:sample_consensus_seq_eq_lambda}
\end{equation}
holds.
Because the Laplace-like distribution satisfies the conditions of
Proposition~\ref{prop:strong_consistency_of_MLE_on_A^ast_times_R},
there exists $N_{1} \in \mathbb{Z}^{+}$ such that if $n \geq N_{1}$,
we have
\begin{equation}
\check{\lambda}(s_{1}, \cdots, s_{n}) = \lambda \as
\label{eq:MLE_lambda}
\end{equation}
Combining Equations~(\ref{eq:sample_consensus_seq_eq_lambda})
and~(\ref{eq:MLE_lambda}) provides $\check{\lambda}(s_{1}, \cdots,
s_{n}) = \bs{m}_{c}^{(\epsilon)}(s_{1}, \cdots, s_{n})$ a.s. for any
$n \geq \max \{ N_{0}, N_{1} \}$.
Consequently, we obtain
\[
\check{\rho}(s_{1}, \cdots, s_{n}) = \frac{1}{n} \sum_{i = 1}^{n}
d_{H^{\prime}} (s_{i}, \bs{m}_{c}^{(\epsilon)}(s_{1}, \cdots, s_{n})) \as
\]
for any $n \geq \max \{ N_{0}, N_{1} \}$ from
Proposition~\ref{prop:mean_abs_dev_and_MLE_dispersion}.
\qed

{\bf Proof of
Lemma~\ref{lemma:EM_algorithm_for_the_Laplace-like_mixture_general_distance}.}
Examining the process of deriving an EM algorithm for a mixture model
(see, for example,~\cite{Mclachlan_1997}), we observe that
Equations~(\ref{eq:zeta_hat_ig_t}) and~(\ref{eq:pi_hat_g_t}) that
provide the formulae for updating $\hat{\zeta}_{i g}^{(t)}$ and
$\hat{\pi}_{g}^{(t)}$, respectively, are common to all distributions
on all spaces.
The object function of the maximization in the M step of the EM
algorithm is given by
\begin{equation}
\sum_{i = 1}^{n} \sum_{g = 1}^{k} \hat{\zeta}_{i g}^{(t)} \left\{ -\log
(\rho_{g} + 1) - \log \left| \partial U(\lambda_{g}, d(s_{i},
\lambda_{g})) \right| + d(s_{i}, \lambda_{g}) \log
\left( \frac{\rho_{g}}{\rho_{g} + 1} \right) \right\}.
\label{eq:object_function_EM_algorithm}
\end{equation}
Because the order of the summations with respect to $i$ and $g$ can be
interchanged, seeking $\lambda_{1}, \cdots, \lambda_{k}$ that maximize
Equation~(\ref{eq:object_function_EM_algorithm}) is equivalent to
seeking $\lambda_{g}$ that maximizes
\[
\sum_{i = 1}^{n} \hat{\zeta}_{i g}^{(t)} \left\{ -\log \left| \partial
U(\lambda_{g}, d(s_{i}, \lambda_{g})) \right| + d(s_{i}, \lambda_{g})
\log \left( \frac{\rho_{g}}{\rho_{g} + 1} \right)
\right\}
\]
for each $g = 1, \cdots, k$, noting $- \log | \partial U(\lambda_{g},
d(s_{i}, \lambda_{g}))| \leq 0$ and $\log (\rho_{g} / (\rho_{g} + 1)) <
0$.
Therefore, Equation~(\ref{eq:lambda_hat_g_t}) provides a procedure for
updating an estimate of $\lambda_{g}$ if the maximization problem can
be solved when the distance $d$ is specified.
For each $g = 1, \cdots, k$, partially differentiating
Equation~(\ref{eq:object_function_EM_algorithm}) with respect to
$\rho_{g}$ leads to
\begin{eqnarray*}
& & \frac{\partial}{\partial \rho_{g}} \sum_{i = 1}^{n}
\sum_{g^{\prime} = 1}^{k} \hat{\zeta}_{i g^{\prime}}^{(t)} \left\{ - \log
(\rho_{g^{\prime}} + 1) - \log \left| \partial U(\lambda_{g^{\prime}}, d(s_{i},
\lambda_{g^{\prime}})) \right| + d(s_{i}, \lambda_{g^{\prime}}) \log \left(
\frac{\rho_{g^{\prime}}}{\rho_{g^{\prime}} + 1} \right) \right\} \\
&=& \sum_{i = 1}^{n} \hat{\zeta}_{i g}^{(t)} \left( -\frac{1}{\rho_{g}
+ 1} + \frac{d(s_{i}, \lambda_{g})}{\rho_{g}} - \frac{d(s_{i},
\lambda_{g})}{\rho_{g} + 1} \right)
 \ = \ \frac{1}{\rho_{g} (\rho_{g} + 1)} \sum_{i = 1}^{n}
\hat{\zeta}_{i g}^{(t)} d(s_{i}, \lambda_{g}) - \frac{\rho_{g}}{\rho_{g}
(\rho_{g} + 1)} \sum_{i = 1}^{n} \hat{\zeta}_{i g}^{(t)}.
\end{eqnarray*}
Therefore, solving the equation
\[
\frac{1}{\rho_{g} (\rho_{g} + 1)} \sum_{i = 1}^{n} \hat{\zeta}_{i g}^{(t)}
d(s_{i}, \lambda_{g}) - \frac{\rho_{g}}{\rho_{g} (\rho_{g} + 1)}
\sum_{i = 1}^{n} \hat{\zeta}_{i g}^{(t)} = 0
\]
with respect to $\rho_{g}$, we obtain
\[
\rho_{g} = \frac{1}{\sum_{i = 1}^{n} \hat{\zeta}_{i g}^{(t)}} \sum_{i
= 1}^{n} \hat{\zeta}_{i g}^{(t)} d(s_{i}, \lambda_{g}).
\]
Hence, Equation~(\ref{eq:update_v_j}) provides a procedure for
updating an estimate of $\rho_{g}$.
\qed

{\bf Proof of
Proposition~\ref{prop:asymptotic_analysis_of_the_EM_algorithm}.}
The topology on $A^{\ast}$ is a discrete topology for any $d \in D$
(see the second paragraph of
Section~\ref{section:Laplace-like_distribution_on_a_set_of_strings}).
Thus, from Equations~(\ref{eq:prob_fun_Laplace-like_dist})
and~(\ref{eq:zeta_hat_ig_t}), $\hat{\zeta}_{i g}^{(n, t)}$ is a
continuous function of $\hat{\pi}_{g}^{(n, t)}, \hat{\lambda}_{g}^{(n,
t)}$, and $\hat{\rho}_{g}^{(n, t)}$ for each $i = 1, \cdots, n$ and
$g = 1, \cdots, k$.
Hence, Part~(a) is obvious.
Therefore, we demonstrate Part~(b).

{\bf (Step 1)}
From the independence of $\bs{Z}_{1}, \cdots, \bs{Z}_{n}$, $Z_{1 g},
\cdots, Z_{n g}$, and consequently $\zeta_{1 g}, \cdots, \zeta_{n g}$,
are also independent for each $g = 1, \cdots, k$.
Moreover, we have $\E_{\bs{\theta}^{\ast}}[\zeta_{i g}^{\ast}] =
\pi_{g}^{\ast}$ and
\[
\Var_{\bs{\theta}^{\ast}} [\zeta_{i g}^{\ast}]
= \sum_{(s_{1}, \cdots, s_{n}) \in (A^{\ast})^{n}} (\zeta_{i g}^{\ast} -
\pi_{g}^{\ast})^{2} \prod_{i = 1}^{n} \bs{q}(s_{i}; \bs{\theta}^{\ast})
< \infty
\]
for each $i = 1, \cdots, n$ and $g = 1, \cdots, k$.
Thus, using the strong law of large numbers provides
\begin{equation}
\frac{1}{n} \sum_{i = 1}^{n} \zeta_{i g}^{\ast} \asconv \pi_{g}^{\ast}
\label{eq:mean_of_zeta_ig_asconv_pi_g^ast}
\end{equation}
as $n_{g} \lto \infty$.
Hence, from the condition of Part~(b) and
Equation~(\ref{eq:pi_hat_g_t}), we obtain $\hat{\pi}_{g}^{(n, t,
\epsilon)} \asconv \pi_{g}^{\ast}$ as $n_{g}, t \lto \infty$ and
$\epsilon \lto \infty$.

{\bf (Step 2)}
Noting that $Z_{1 g}, \cdots, Z_{n g}$ are independent and that
$\E_{\pi_{g}^{\ast}}[Z_{i g}] = \pi^{\ast}_{g}$ and
$\Var_{\pi_{g}^{\ast}}[Z_{i g}] < \infty$ hold from
Equation~(\ref{eq:distribution_of_Z_i}) and applying the strong law of
large numbers, we have $\sum_{i = 1}^{n} z_{i g} / n \asconv
\pi_{g}^{\ast}$ as $n_{g} \lto \infty$.
Therefore, from Equation~(\ref{eq:mean_of_zeta_ig_asconv_pi_g^ast}),
\begin{equation}
\frac{1}{n} \sum_{i = 1}^{n} z_{i g} \asconv \frac{1}{n} \sum_{i =
1}^{n} \zeta_{i g}^{\ast}
\label{eq:mean_of_z_i_g_to_mean_of_zeta_i_g^ast}
\end{equation}
holds as $n_{g} \lto \infty$.
We obtain
\begin{equation}
\frac{1}{n} \sum_{i \in \{ i^{\prime} \in \{ 1, \cdots, n \}:
x_{i^{\prime} j} = a_{h} \}} z_{i g} \asconv \frac{1}{n} \sum_{i \in
\{ i^{\prime} \in \{ 1, \cdots, n \}: x_{i^{\prime} j} = a_{h} \}}
\zeta_{i g}^{\ast}
\label{eq:var_mean_of_z_ig_asconv_var_mean_of_zeta_ig^ast}
\end{equation}
as $n_{g} \lto \infty$ in the same manner.
Combining Equations~(\ref{eq:mean_of_z_i_g_to_mean_of_zeta_i_g^ast})
and~(\ref{eq:var_mean_of_z_ig_asconv_var_mean_of_zeta_ig^ast}) with
the condition of Part~(b) gives
\begin{equation}
\frac{1}{n} \sum_{i = 1}^{n} \hat{\zeta}_{i g}^{(n, t, \epsilon)}
\asconv \frac{1}{n} \sum_{i = 1}^{n} z_{i g}, \quad
\frac{1}{n} \sum_{i \in \{ i^{\prime} \in \{ 1, \cdots, n \}: x_{i^{\prime} j}
= a_{h} \}} \hat{\zeta}_{i g}^{(n, t, \epsilon)}
\asconv \frac{1}{n} \sum_{i \in \{ i^{\prime} \in \{ 1, \cdots, n \}: x_{i^{\prime} j}
= a_{h} \}} z_{i g} \label{eq:mean_of_zeta_hat_ig_mean_of_conv_z_ig}
\end{equation}
as $n_{g}, t \lto \infty$ and $\epsilon \lto 0$.
Thus, from Equation~(\ref{eq:f_g_j_h}), we obtain $\varphi_{g j h}
\asconv f_{g j h}$ as $n_{g}, t \lto \infty$ and $\epsilon \lto 0$.
Hence, using Equations~(\ref{eq:def_m_c^epsilon})
and~(\ref{eq:lambda_hat_g}), there exist $N_{0}, T_{0} \in
\mathbb{Z}^{+}$ and $\epsilon_{0} > 0$ such that if $n_{g} \geq N_{0},
t \geq T_{0}$, and $\epsilon \leq \epsilon_{0}$, then
\begin{equation}
\hat{\lambda}_{g}^{(n, t, \epsilon)} = \bs{m}_{c}^{(\epsilon)}(s_{g 1},
\cdots, s_{g n_{g}}) \as \label{eq:lambda_hat_g_equal_m_c^epsilon_as}
\end{equation}
Applying Proposition~\ref{prop:MLE_location_parameter} to $\sigma_{g
1}, \cdots, \sigma_{g n_{g}}$, there exist $N_{1} \in
\mathbb{Z}^{+}$ and $\epsilon_{1} > 0$ such that if $n_{g} \geq N_{1}$
and $\epsilon \leq \epsilon_{1}$, then
\begin{equation}
\bs{m}_{c}^{(\epsilon)}(s_{g 1}, \cdots, s_{g n_{g}}) = \lambda_{g}^{\ast}
\as \label{eq:m_c^epsilon_equal_lambda_g^ast_as}
\end{equation}
Noting Equations~(\ref{eq:lambda_hat_g_equal_m_c^epsilon_as})
and~(\ref{eq:m_c^epsilon_equal_lambda_g^ast_as}), we see that there
exist $N_{2}, T_{2} \in \mathbb{Z}^{+}$ and $\epsilon_{2} > 0$ such
that if $n_{g} \geq N_{2}, t \geq T_{2}$, and $\epsilon \leq
\epsilon_{2}$, then $\hat{\lambda}_{g}^{(n, t, \epsilon)} =
\lambda_{g}^{\ast}$ a.s. holds.

{\bf (Step 3)}
From Equation~(\ref{eq:mean_of_z_i_g_to_mean_of_zeta_i_g^ast}), we
have
\[
\frac{1}{n} \sum_{i = 1}^{n} z_{i g} d_{H^{\prime}}(s_{i},
\lambda_{g}^{\ast}) \asconv \frac{1}{n} \sum_{i = 1}^{n} \zeta_{i
g}^{\ast} d_{H^{\prime}}(s_{i}, \lambda_{g}^{\ast})
\]
as $n_{g} \lto \infty$.
Thus, from the condition of Part~(b),
\begin{equation}
\frac{1}{n} \sum_{i = 1}^{n} \hat{\zeta}_{i g}^{(n, t)}
d_{H^{\prime}}(s_{i}, \lambda_{g}^{\ast}) \asconv \frac{1}{n} \sum_{i =
1}^{n} z_{i g} d_{H^{\prime}}(s_{i}, \lambda_{g}^{\ast})
\label{eq:sum_zeta_hat_d_H_prime_to_sum_z_d_H_prime}
\end{equation}
holds as $n_{g}, t \lto \infty$ and $\epsilon \lto 0$.
Using the left of
Equation~(\ref{eq:mean_of_zeta_hat_ig_mean_of_conv_z_ig}) and
Equation~(\ref{eq:sum_zeta_hat_d_H_prime_to_sum_z_d_H_prime}), we have
\begin{equation}
\hat{\rho}_{g}^{(n, t, \epsilon)}
= \frac{1}{\sum_{i = 1}^{n} \hat{\zeta}_{i g}^{(n, t, \epsilon)}}
\sum_{i = 1}^{n} \hat{\zeta}_{i g}^{(n, t, \epsilon)} d_{H^{\prime}}(s_{i},
\hat{\lambda}_{g}^{(n, t, \epsilon)})
\asconv \frac{1}{\sum_{i = 1}^{n} z_{i g}} \sum_{i = 1}^{n} z_{i g}
d_{H^{\prime}}(s_{i}, \hat{\lambda}_{g}^{(n, t, \epsilon)})
\label{eq:rho_hat_g_conv_1}
\end{equation}
as $n_{g}, t \lto \infty$ and $\epsilon \lto \infty$.
Noting the result of Step 2 and applying
Proposition~\ref{prop:MLE_dispersion_parameter} to $\sigma_{g
1}, \cdots, \sigma_{g n_{g}}$, we obtain
\begin{equation}
\frac{1}{\sum_{i = 1}^{n} z_{i g}} \sum_{i = 1}^{n} z_{i g}
d_{H^{\prime}}(s_{i}, \hat{\lambda}_{g}^{(n, t, \epsilon)})
\asconv \frac{1}{n_{g}} \sum_{i = 1}^{n_{g}} d_{H^{\prime}} (s_{g i},
\lambda_{g}^{\ast})
\asconv \rho_{g}^{\ast}
\label{eq:rho_hat_g_conv_2}
\end{equation}
as $n_{g}, t \lto \infty$ and $\epsilon \lto 0$.
From Equations~(\ref{eq:rho_hat_g_conv_1})
and~(\ref{eq:rho_hat_g_conv_2}), $\hat{\rho}_{g}^{(n, t, \epsilon)}
\asconv \rho_{g}^{\ast}$ holds as $n_{g}, t \lto \infty$ and $\epsilon
\lto \infty$.
Combining the results of Steps~1 to~3 completes the proof.
\qed

{\bf Proof of
Theorem~\ref{theorem:strong_consistency_of_the_EM_algorithm}.}
{\bf (Step 1)}
Let $\check{\lambda}_{g}^{(n_{g})}$ and $\check{\rho}_{g}^{(n_{g})}$
represent the maximum likelihood estimators of $\lambda_{g}$ and
$\rho_{g}$ based on $s_{g 1}, \cdots, s_{g n_{g}}$, respectively, for
each $g = 1, \cdots, k$.
Applying Theorem~\ref{theorem:asymptotic_MLE} to $\sigma_{g 1},
\cdots, \sigma_{g n_{g}}$, there exists $N_{0} \in \mathbb{Z}^{+}$
such that if $n^{\ast} \geq N_{0}$, then
$\check{\lambda}_{1}^{(n_{1})}, \cdots, \check{\lambda}_{k}^{(n_{k})},
\check{\rho}_{1}^{(n_{1})}, \cdots, \check{\rho}_{k}^{(n_{k})}$ are
uniquely determined with probability one.
Therefore, if $n^{\ast} \geq N_{0}$, we have
\begin{equation}
\frac{1}{n} \sum_{g = 1}^{k} \sum_{i = 1}^{n} z_{i g} \log q(s_{i};
\lambda_{g}, \rho_{g}) < \frac{1}{n} \sum_{g = 1}^{k} \sum_{i =
1}^{n} z_{i g} \log q(s_{i}; \check{\lambda}_{g}^{(n_{g})},
\check{\rho}_{g}^{(n_{g})}) \as
\label{eq:ineq_MLE}
\end{equation}
for any $(\lambda_{1}, \cdots, \lambda_{k}, \rho_{1}, \cdots,
\rho_{k}) \in (A^{\ast})^{k} \times (0, \infty)^{k} \smallsetminus \{
(\check{\lambda}_{1}^{(n_{1})}, \cdots, \check{\lambda}_{k}^{(n_{k})},
\check{\rho}_{1}^{(n_{1})}, \cdots, \check{\rho}_{k}^{(n_{k})}) \}$.
By
Proposition~\ref{prop:strong_consistency_of_MLE_on_A^ast_times_R},
there exists $N_{1} \in \mathbb{Z}^{+}$ such that if $n^{\ast} \geq
N_{1}$, then
\begin{equation}
(\check{\lambda}_{1}^{(n_{1})}, \cdots, \check{\lambda}_{k}^{(n_{k})}) =
(\lambda_{1}^{\ast}, \cdots, \lambda_{k}^{\ast}) \as
\label{eq:lambda_check_aseq_lambda}
\end{equation}
and
\begin{equation}
(\check{\rho}_{1}^{(n_{1})}, \cdots, \check{\rho}_{k}^{(n_{k})}) \asconv
(\rho_{1}^{\ast}, \cdots, \rho_{k}^{\ast})
\label{eq:rho_check_asconv_rho}
\end{equation}
as $n^{\ast} \lto \infty$.

{\bf (Step 2)}
Let $\tilde{\zeta}_{i g}^{(n, t, \epsilon)}, \tilde{\lambda}_{g}^{(n,
t, \epsilon)}$, and $\tilde{\rho}_{g}^{(n, t, \epsilon)}$ be
estimates of $\zeta_{i g}, \lambda_{g}$, and $\rho_{g}$ obtained at
iteration step $t$ from the initial value
$\tilde{\bs{\theta}}^{(n, 0)}$ using Algorithm $H^{\prime}$,
respectively.
Noting condition~(i) and using
Equations~(\ref{eq:lambda_check_aseq_lambda})
and~(\ref{eq:rho_check_asconv_rho}), there exist $N_{2}, T_{2} \in
\mathbb{Z}^{+}$ and $\epsilon_{0} > 0$ such that if $n^{\ast} \geq
N_{2}, t \geq T_{2}$, and $\epsilon \leq \epsilon_{0}$, then
\begin{equation}
(\tilde{\lambda}_{1}^{(n, t, \epsilon)}, \cdots,
\tilde{\lambda}_{k}^{(n, t, \epsilon)}) = (\check{\lambda}_{1}^{(n_{1})},
\cdots, \check{\lambda}_{k}^{(n_{k})}) \as
\label{eq:lambda_tilde_aseq_lambda_check}
\end{equation}
and
\begin{equation}
(\tilde{\rho}_{1}^{(n, t, \epsilon)}, \cdots,
\tilde{\rho}_{k}^{(n, t, \epsilon)}) \asconv (\check{\rho}_{1}^{(n_{1})},
\cdots, \check{\rho}_{k}^{(n_{k})})
\label{eq:rho_tilde_asconv_rho_check}
\end{equation}
as $n^{\ast}, t \lto \infty$ and $\epsilon \lto 0$.
Applying Part~(a) of
Proposition~\ref{prop:asymptotic_analysis_of_the_EM_algorithm} from
condition~(i), we have
\begin{equation}
(\tilde{\zeta}_{1 1}^{(n, t, \epsilon)}, \cdots, \tilde{\zeta}_{n
k}^{(n, t, \epsilon)}) \asconv (\zeta_{1 1}^{\ast}, \cdots,
\zeta_{n k}^{\ast}) \label{eq:zeta_tilde_asconv_zeta_ast}
\end{equation}
as $n^{\ast}, t \lto \infty$ and $\epsilon \lto 0$.
Hence, we obtain
\begin{equation}
\frac{1}{n} \sum_{g = 1}^{k} \sum_{i = 1}^{n} \tilde{\zeta}_{i g}^{(n,
t, \epsilon)} \log q(s_{i}; \tilde{\lambda}_{g}^{(n, t, \epsilon)},
\tilde{\rho}_{g}^{(n, t, \epsilon)}) \asconv \frac{1}{n} \sum_{g =
1}^{k} \sum_{i = 1}^{n} z_{i g} \log q(s_{i};
\tilde{\lambda}_{g}^{(n, t, \epsilon)}, \tilde{\rho}_{g}^{(n, t,
\epsilon)}) \label{eq:zeta_tilde_log_q_asconv_z_log_q}
\end{equation}
as $n^{\ast}, t \lto \infty$ and $\epsilon \lto 0$ in a similar manner
to deriving
Equation~(\ref{eq:sum_zeta_hat_d_H_prime_to_sum_z_d_H_prime}).
Combining
Equations~(\ref{eq:ineq_MLE}),~(\ref{eq:lambda_tilde_aseq_lambda_check}),~(\ref{eq:rho_tilde_asconv_rho_check}),
and~(\ref{eq:zeta_tilde_log_q_asconv_z_log_q}) provides
\begin{equation}
\frac{1}{n} \sum_{g = 1}^{k} \sum_{i = 1}^{n} z_{i g} \log q(s_{i};
\lambda_{g}, \rho_{g}) \leq \frac{1}{n} \sum_{g = 1}^{k} \sum_{i =
1}^{n} \tilde{\zeta}_{i g}^{(n, t, \epsilon)} \log q(s_{i};
\tilde{\lambda}_{g}^{(n, t, \epsilon)}, \tilde{\rho}_{g}^{(n, t,
\epsilon)}) \as
\label{eq:key_inequality}
\end{equation}
for any $(\lambda_{1}, \cdots, \lambda_{k}, \rho_{1}, \cdots,
\rho_{k}) \in (A^{\ast})^{k} \times (0, \infty)^{k} \smallsetminus \{
(\check{\lambda}_{1}^{(n_{1})}, \cdots, \check{\lambda}_{k}^{(n_{k})},
\check{\rho}_{1}^{(n_{1})}, \cdots, \check{\rho}_{k}^{(n_{k})}) \}$ as
$n^{\ast}, t \lto \infty$ and $\epsilon \lto 0$.
From Equations~(\ref{eq:ineq_MLE}) and~(\ref{eq:key_inequality}) and
condition C$_{1}$, there exist $N_{3}, T_{3} \in \mathbb{Z}^{+}$ and
$\epsilon_{1} > 0$ such that if $n^{\ast} \geq N_{3}, t \geq T_{3}$,
and $\epsilon \leq \epsilon_{1}$, then
\begin{equation}
(\tilde{\lambda}_{1}^{(n, t, \epsilon)}, \cdots,
\tilde{\lambda}_{k}^{(n, t, \epsilon)}) = (\lambda_{1}^{\dagger},
\cdots, \lambda_{k}^{\dagger}) \as
\label{eq:lambda_tilde_aseq_lambda_ast}
\end{equation}
and
\begin{equation}
(\tilde{\zeta}_{1 1}^{(n, t, \epsilon)}, \cdots, \tilde{\zeta}_{n
k}^{(n, t, \epsilon)}) \asconv (\zeta_{1 1}^{\dagger}, \cdots,
\zeta_{n k}^{\dagger}), \quad (\tilde{\rho}_{1}^{(n, t, \epsilon)},
\cdots, \tilde{\rho}_{k}^{(n, t, \epsilon)}) \asconv
(\rho_{1}^{\dagger}, \cdots, \rho_{k}^{\dagger})
\label{eq:rho_tilde_asconv_rho_ast}
\end{equation}
as $n^{\ast}, t \lto \infty$ and $\epsilon \lto 0$.
In other words, combined with condition C$_{1}$ and
Proposition~\ref{prop:strong_consistency_of_MLE_on_A^ast_times_R},
condition~(i) means that there exists an initial value
$\tilde{\bs{\theta}}^{(n, 0)}$ with which Algorithm $H^{\prime}$
returns a sequence $\{ \tilde{\zeta}_{1 1}^{(n, t, \epsilon)}, \cdots,
\tilde{\zeta}_{n k}^{(n, t, \epsilon)}, \tilde{\lambda}_{1}^{(n, t,
\epsilon)}, \cdots, \tilde{\lambda}_{k}^{(n, t, \epsilon)},
\tilde{\rho}_{1}^{(n, t, \epsilon)}, \cdots, \tilde{\rho}_{k}^{(n, t,
\epsilon)} \}$ of estimates that almost surely converges to the
maximizer $(\zeta_{1 1}^{\dagger}, \cdots, \zeta_{n k}^{\dagger},
\lambda_{1}^{\dagger}, \cdots, \lambda_{k}^{\dagger},
\rho_{1}^{\dagger}, \cdots, \rho_{k}^{\dagger})$ of
Equation~(\ref{eq:fun_of_zeta_lambda_rho}) as $n^{\ast}, t \lto
\infty$ and $\epsilon \lto 0$.

{\bf (Step 3)}
$\hat{\bs{\theta}}^{(n, 0, \tau, \epsilon)}$ is an initial value
with which Algorithm $H^{\prime}$ returns estimates of $\zeta_{1 1},
\cdots, \zeta_{n k},\\ \lambda_{1}, \cdots, \lambda_{k}, \rho_{1},
\cdots, \rho_{k}$ that maximize
Equation~(\ref{eq:fun_of_zeta_lambda_rho}) in a set of their estimates
at iteration step $\tau$ for all possible initial values.
Therefore, from the result of Step~2, there exist $N_{4}, T_{4} \in
\mathbb{Z}^{+}$ and $\epsilon_{2} > 0$ such that if $n^{\ast} \geq
N_{4}, t, \tau \geq T_{4}$, and $\epsilon \leq \epsilon_{2}$,
then
\begin{equation}
(\hat{\lambda}_{1}^{(n, t, \tau, \epsilon)}, \cdots,
\hat{\lambda}_{k}^{(n, t, \tau, \epsilon)}) =
(\lambda_{1}^{\dagger}, \cdots, \lambda_{k}^{\dagger}) \as
\label{eq:lambda_hat_asconv_lambda_ast}
\end{equation}
and
\begin{equation}
(\hat{\zeta}_{1 1}^{(n, t, \tau, \epsilon)}, \cdots,
\hat{\zeta}_{n k}^{(n, t, \tau, \epsilon)}) \asconv
(\zeta_{1 1}^{\dagger}, \cdots, \zeta_{n k}^{\dagger}), \quad
(\hat{\rho}_{1}^{(n, t, \tau, \epsilon)}, \cdots,
\hat{\rho}_{k}^{(n, t, \tau, \epsilon)}) \asconv
(\rho_{1}^{\dagger}, \cdots, \rho_{k}^{\dagger})
\label{eq:rho_hat_asconv_rho_ast}
\end{equation}
as $n^{\ast}, t, \tau \lto \infty$ and $\epsilon \lto 0$.
Noting condition C$_{2}$ and combining
Equations~(\ref{eq:lambda_tilde_aseq_lambda_ast})
to~(\ref{eq:rho_hat_asconv_rho_ast}), there exist $N_{5}, T_{5} \in
\mathbb{Z}^{+}$ and $\epsilon_{3} > 0$ such that if $n^{\ast} \geq
N_{5}, t, \tau \geq T_{5}$, and $\epsilon \leq \epsilon_{3}$, then
\begin{equation}
(\hat{\lambda}_{1}^{(n, t, \tau, \epsilon)}, \cdots,
\hat{\lambda}_{k}^{(n, t, \tau, \epsilon)}) =
(\tilde{\lambda}_{1}^{(n, t, \epsilon)}, \cdots,
\tilde{\lambda}_{k}^{(n, t, \epsilon)}) \as
\label{eq:lambda_hat_aseq_lambda_tilde}
\end{equation}
and
\begin{equation}
(\hat{\zeta}_{1 1}^{(n, t, \tau, \epsilon)}, \cdots,
\hat{\zeta}_{n k}^{(n, t, \tau, \epsilon)}) \asconv
(\tilde{\zeta}_{1 1}^{(n, t, \epsilon)}, \cdots,
\tilde{\zeta}_{n k}^{(n, t, \epsilon)}), \quad
(\hat{\rho}_{1}^{(n, t, \tau, \epsilon)}, \cdots,
\hat{\rho}_{k}^{(n, t, \tau, \epsilon)}) \asconv
(\tilde{\rho}_{1}^{(n, t, \epsilon)}, \cdots,
\tilde{\rho}_{k}^{(n, t, \epsilon)}) \label{eq:rho_hat_asconv_rho_tilde}
\end{equation}
as $n^{\ast}, t, \tau \lto \infty$ and $\epsilon \lto 0$.
By Equation~(\ref{eq:lambda_hat_aseq_lambda_tilde}), the right of
Equation~(\ref{eq:rho_hat_asconv_rho_tilde}), and condition~(i),
$(\hat{\lambda}_{1}^{(n, t, \tau, \epsilon)}, \cdots,
\hat{\lambda}_{k}^{(n, t, \tau, \epsilon)}, \hat{\rho}_{1}^{(n,
t, \tau, \epsilon)}, \cdots, \hat{\rho}_{k}^{(n, t,
\tau, \epsilon)})$ strongly consistently estimates
$(\lambda_{1}, \cdots, \lambda_{k}, \rho_{1}, \cdots, \rho_{k})$ as
$n^{\ast}, t, \tau \lto \infty$ and $\epsilon \lto 0$.
Combining Equation~(\ref{eq:zeta_tilde_asconv_zeta_ast}) and the left
of Equation~(\ref{eq:rho_hat_asconv_rho_tilde}) provides
\[
(\hat{\zeta}_{1 1}^{(n, t, \tau, \epsilon)}, \cdots,
\hat{\zeta}_{n k}^{(n, t, \tau, \epsilon)}) \asconv (\zeta_{1
1}^{\ast}, \cdots, \zeta_{n k}^{\ast})
\]
as $n^{\ast}, t, \tau \lto \infty$ and $\epsilon \lto 0$.
Hence, using Part~(b) of
Proposition~\ref{prop:asymptotic_analysis_of_the_EM_algorithm}, we see
that $(\hat{\pi}_{1}^{(n, t, \tau, \epsilon)}, \cdots,
\hat{\pi}_{k}^{(n, t, \tau, \epsilon)})$ strongly consistently
estimates $(\pi_{1}, \cdots, \pi_{k})$ as $n^{\ast}, t, \tau
\lto \infty$ and $\epsilon \lto 0$.
\qed

{\bf Proof of Theorem~\ref{theorem:asymptotic_EM}.}
Algorithm $H^{\prime}$ forms a sequence of algorithms with respect to
the three parameters $n^{\ast}, \tau$, and $\epsilon$ in this theorem,
whereas a sequence of algorithms has one parameter $n$ in
Definition~\ref{def:convergence_of_a_sequence_of_algorithms}.
It suffices to demonstrate that there exist $N_{0}, T_{0} \in
\mathbb{Z}^{+}$ and $\epsilon_{0} > 0$ such that if $n_{g} \geq
N_{0}, t, \tau \geq T_{0}$, and $\epsilon \leq \epsilon_{0}$,
then $\hat{\lambda}_{g}^{(n, t, \tau, \epsilon)}$ is equal to
the maximizer of Equation~(\ref{eq:lambda_hat_g_t}) with $d =
d_{H^{\prime}}$ with probability one for each $g = 1, \cdots, k$.

{\bf (Step 1)}
Under the conditions of
Theorem~\ref{theorem:strong_consistency_of_the_EM_algorithm}, the
estimate $\hat{\bs{\theta}}^{(n, t, \tau, \epsilon)}$ at iteration
step $t$ from Algorithm $H^{\prime}$ with the initial value
$\hat{\bs{\theta}}^{(n, 0, \tau, \epsilon)}$ that satisfies
Equation~(\ref{eq:fun_of_zeta_lambda_rho}) almost surely converges to
the true value $\bs{\theta}^{\ast}$ of the parameter as $n^{\ast}, t,
\tau \lto \infty$ and $\epsilon \lto 0$.
Therefore, using Part~(a) of
Proposition~\ref{prop:asymptotic_analysis_of_the_EM_algorithm}, we
have
\[
(\hat{\zeta}_{1 1}^{(n, t, \tau, \epsilon)}, \cdots,
\hat{\zeta}_{n k}^{(n, t, \tau, \epsilon)}) \asconv (\zeta_{1
1}^{\ast}, \cdots, \zeta_{n k}^{\ast})
\]
as $n^{\ast}, t, \tau \lto \infty$ and $\epsilon \lto 0$.
Thus, we obtain
\begin{eqnarray}
& & \frac{1}{n} \sum_{i = 1}^{n} \hat{\zeta}_{i g}^{(n, t, \tau,
\epsilon)} \left\{ - \log \left| \partial U(\lambda_{g},
d_{H^{\prime}}(s_{i}, \lambda_{g})) \right| + d_{H^{\prime}}(s_{i},
\lambda_{g}) \log \left( \frac{\hat{\rho}_{g}^{(n, t - 1, \tau,
\epsilon)}}{\hat{\rho}_{g}^{(n, t - 1, \tau, \epsilon)} + 1} \right)
\right\} \nonumber \\
&\asconv&
\frac{1}{n} \sum_{i = 1}^{n} z_{i g} \left\{ - \log \left| \partial
U(\lambda_{g}, d_{H^{\prime}}(s_{i}, \lambda_{g})) \right| +
d_{H^{\prime}}(s_{i}, \lambda_{g}) \log \left( \frac{\hat{\rho}_{g}^{(n,
t - 1, \tau, \epsilon)}}{\hat{\rho}_{g}^{(n, t - 1, \tau,
\epsilon)} + 1} \right) \right\} \label{eq:EM_obj_fun_conv} \\
&=&
\frac{1}{n} \sum_{i = 1}^{n_{g}} \left\{ - \log \left| \partial
U(\lambda_{g}, d_{H^{\prime}}(s_{g i}, \lambda_{g})) \right| +
d_{H^{\prime}}(s_{g i}, \lambda_{g}) \log \left( \frac{\hat{\rho}_{g}^{(n,
t - 1, \tau, \epsilon)}}{\hat{\rho}_{g}^{(n, t - 1, \tau,
\epsilon)} + 1} \right) \right\} \nonumber
\end{eqnarray}
as $n_{g}, t, \tau \lto \infty$ and $\epsilon \lto 0$ in a
similar manner to obtaining
Equation~(\ref{eq:sum_zeta_hat_d_H_prime_to_sum_z_d_H_prime}).
In other words, the object function~(\ref{eq:lambda_hat_g_t}) of the
maximization almost surely converges to the log likelihood function of
$\lambda_{g}$ based on $s_{g 1}, \cdots, s_{g n_{g}}$.

{\bf (Step 2)}
Let $\check{\lambda}_{g}(s_{g 1}, \cdots, s_{g n_{g}})$ denote the
maximizer of the right-hand side of
Equation~(\ref{eq:EM_obj_fun_conv}), i.e., the maximum likelihood
estimate of $\lambda_{g}$ based on $s_{g 1}, \cdots, s_{g n_{g}}$.
By Proposition~\ref{prop:strong_consistency_of_MLE_on_A^ast_times_R},
there exists $N_{1} \in \mathbb{Z}^{+}$ such that if $n_{g} \geq
N_{1}$, then
\[
\check{\lambda}_{g}(s_{g 1}, \cdots, s_{g n_{g}}) = \lambda^{\ast}_{g}
\as
\]
holds, and therefore
$\check{\lambda}_{g}(s_{g 1}, \cdots, s_{g n_{g}})$ is uniquely
determined with probability one.
Applying Theorem~\ref{theorem:asymptotic_MLE} to $\sigma_{g 1},
\cdots, \sigma_{g n_{g}}$, there exist $N_{2} \in \mathbb{Z}^{+}$ and
$\epsilon_{2} > 0$ such that if $n_{g} \geq N_{2}$ and $\epsilon \leq
\epsilon_{2}$, we have
\begin{equation}
\bs{m}_{c}^{(\epsilon)}(s_{g 1}, \cdots, s_{g n_{g}}) =
\check{\lambda}_{g}(s_{g 1}, \cdots, s_{g n_{g}}) \as
\label{eq:m_c^epsilon_equal_lambda_check_g}
\end{equation}
We observe that there exist $N_{3}, T_{3} \in \mathbb{Z}^{+}$ and
$\epsilon_{3} > 0$ such that if $n_{g} \geq N_{3}, t, \tau \geq
T_{3}$, and $\epsilon \leq \epsilon_{3}$, then
\begin{equation}
\hat{\lambda}_{g}^{(n, t, \tau, \epsilon)} =
\bs{m}_{c}^{(\epsilon)}(s_{g 1}, \cdots, s_{g n_{g}}) \as
\label{eq:lambda_hat_g_equal_m_c^epsilon}
\end{equation}
holds in the same manner as obtaining
Equation~(\ref{eq:lambda_hat_g_equal_m_c^epsilon_as}).
Combining Equations~(\ref{eq:m_c^epsilon_equal_lambda_check_g})
and~(\ref{eq:lambda_hat_g_equal_m_c^epsilon}), there exist $N_{4},
T_{4} \in \mathbb{Z}^{+}$ and $\epsilon_{4} > 0$ such that if $n_{g}
\geq N_{4}, t, \tau \geq T_{4}$, and $\epsilon \leq \epsilon_{4}$, we
have
\[
\hat{\lambda}_{g}^{(n, t, \tau, \epsilon)} = \check{\lambda}_{g}(s_{g
1}, \cdots, s_{g n_{g}}) \as
\]
Thus, supposing that there exists $\tilde{\lambda}_{g} \in A^{\ast}$
such that
\begin{eqnarray*}
& & \frac{1}{n} \sum_{i = 1}^{n} \hat{\zeta}_{i g}^{(n, t, \tau,
\epsilon)} \left\{ - \log \left| \partial U(\hat{\lambda}_{g}^{(n,
t, \tau, \epsilon)}, d_{H^{\prime}}(s_{i}, \hat{\lambda}_{g}^{(n, t,
\tau, \epsilon)})) \right| + d_{H^{\prime}}(s_{i},
\hat{\lambda}_{g}^{(n, t, \tau, \epsilon)}) \log \left(
\frac{\hat{\rho}_{g}^{(n, t - 1, \tau, \epsilon)}}{\hat{\rho}_{g}^{(n,
t - 1, \tau, \epsilon)} + 1} \right) \right\} \\
&<& \frac{1}{n} \sum_{i = 1}^{n} \hat{\zeta}_{i g}^{(n, t, \tau,
\epsilon)} \left\{ - \log \left| \partial U(\tilde{\lambda}_{g},
d_{H^{\prime}}(s_{i}, \tilde{\lambda}_{g})) \right| + d_{H^{\prime}}(s_{i},
\tilde{\lambda}_{g}) \log \left( \frac{\hat{\rho}_{g}^{(n, t - 1,
\tau, \epsilon)}}{\hat{\rho}_{g}^{(n, t - 1, \tau, \epsilon)}
+ 1} \right) \right\} \as
\end{eqnarray*}
holds as $n_{g}, t, \tau \lto \infty$ and $\epsilon \lto \infty$ leads
to a contradiction with Equation~(\ref{eq:EM_obj_fun_conv}) because
$\check{\lambda}_{g}(s_{g 1}, \cdots, s_{g n_{g}})$ is the unique
maximizer of the right-hand side of
Equation~(\ref{eq:EM_obj_fun_conv}) for $n_{g} \geq N_{1}$ with
probability one.
Hence, there exist $N_{5}, T_{5} \in \mathbb{Z}^{+}$ and $\epsilon_{5}
> 0$ such that if $n_{g} \geq N_{5}, t, \tau \geq T_{5}$, and
$\epsilon \leq \epsilon_{5}$, then $\hat{\lambda}_{g}^{(n, t, \tau,
\epsilon)}$ is equal to the maximizer of the left-hand side of
Equation~(\ref{eq:EM_obj_fun_conv}) with probability one.
This completes the proof.
\qed

{\bf Proof of
Corollary~\ref{corollary:optimality_of_clustering_procedure}.}
Obvious from the manner of constructing the clustering procedure and
Theorem~\ref{theorem:strong_consistency_of_the_EM_algorithm}.
\qed


\end{document}